\documentclass[11pt]{article}
\usepackage{xcolor}
 \usepackage{geometry}                % See geometry.pdf to learn the layout options. There are lots.
\usepackage{graphicx}
\usepackage{subcaption}
\usepackage{authblk}

\usepackage{amssymb}
\usepackage[normalem]{ulem}
\usepackage{hyperref}
\usepackage{enumitem}
\usepackage{mathrsfs}
\usepackage{epstopdf}
\usepackage{color}
\usepackage{bbm, dsfont}
\usepackage{tikz-cd}
\usepackage{amsfonts} %mathbb{}
\usepackage{geometry}
\usepackage{amsthm}
\usepackage{amsmath}
\usepackage{titlesec}
\usepackage{fixltx2e}
\usepackage{amssymb}
\usepackage{enumitem}
\usepackage{float}
\restylefloat{table}
\usepackage{tikz}
\usetikzlibrary{tikzmark}
\usetikzlibrary{arrows}
\usetikzlibrary{cd}
\usetikzlibrary{automata, positioning, arrows}
\usepackage [english]{babel}
\usepackage [autostyle, english = american]{csquotes}
\usepackage{algorithm}
\usepackage[noend]{algpseudocode} %pseudocode
\makeatletter
\def\BState{\State\hskip-\ALG@thistlm}
\makeatother
\def\acts{\curvearrowright}

\def\X^n{\text{X}^n}

\usepackage{hyperref}

\newcommand\addtag{\refstepcounter{equation}\tag{\theequation}}
\makeatletter
\newcommand*{\rom}[1]{\expandafter\@slowromancap\romannumeral #1@}
\makeatother

\theoremstyle{theorem}
\newtheorem{thm}{Theorem}[section]
\newtheorem{prop}[thm]{Proposition}
\newtheorem{cor}[thm]{Corollary}
\newtheorem{lemma}[thm]{Lemma}

\theoremstyle{definition}
\newtheorem{defn}{Definition}[section]
\newtheorem{rmk}[thm]{Remark}
\newenvironment{defnis}[1]
  {\renewcommand{\thedefn}{\ref{#1}$'$}%
   \addtocounter{defn}{-1}%
   \begin{defn}}
  {\end{defn}}

\newtheorem{exmp}[thm]{Example}

\usepackage{etoolbox}

\makeatletter
\patchcmd{\ttlh@hang}{\parindent\z@}{\parindent\z@\leavevmode}{}{}
\patchcmd{\ttlh@hang}{\noindent}{}{}{}
\makeatother

\usepackage{listings}
\usepackage{color} %red, green, blue, yellow, cyan, magenta, black, white
\definecolor{mygreen}{RGB}{28,172,0} % color values Red, Green, Blue
\definecolor{mylilas}{RGB}{170,55,241}

\newlist{Assumptions}{enumerate}{1}
\setlist[Assumptions]{label=     \textbf{Assumption} \arabic*:}

\makeatletter

\newsavebox{\@brx}
\newcommand{\llangle}[1][]{\savebox{\@brx}{\(\m@th{#1\langle}\)}%
  \mathopen{\copy\@brx\kern-0.5\wd\@brx\usebox{\@brx}}}
\newcommand{\rrangle}[1][]{\savebox{\@brx}{\(\m@th{#1\rangle}\)}%
  \mathclose{\copy\@brx\kern-0.5\wd\@brx\usebox{\@brx}}}
\makeatother

\usepackage{lipsum} % for generating filler text
\usepackage{titlesec}
\titleformat{\subsection}[runin]% runin puts it in the same paragraph
       {\normalfont\bfseries}% formatting commands to apply to the whole heading
       {\thesubsection}% the label and number
       {0.5em}% space between label/number and subsection title
       {}% formatting commands applied just to subsection title
       [.]% punctuation or other commands following subsection title

\begin{document}
\title{Growth of groups with incompressible elements, \rom{1}}
\author{Zheng Kuang\thanks{Texas A\&M University, College Station, USA; Email: kzkzkzz@tamu.edu}}

\date{}
\maketitle

\begin{abstract}
We define the class of groups of bounded type from tile inflations. These tile inflations also determine some automata describing the groups. In the case when the automata are stationary, we show that if the set of incompressible elements of a group in this class is finite, then this group has subexponential growth with a bounded power in the exponent. Then we describe some examples with certain special structures of orbital graphs and give explicit ways to find the upper bounds for the growth functions of these groups. 
\end{abstract}
\tableofcontents
\section{Introduction}
  We study a class of groups containing many examples of groups of intermediate growth, extending the methods of the papers~\cite{BNZ} and~\cite{nekrash18} to groups with orbital graphs not quasi-isometric to a line. Similarly to the classical examples of groups of intermediate growth, such as the Grigorchuk group~\cite{gri84}, the groups in our class have growth functions bounded from above by functions of the form $\exp(CR^\alpha)$ for $\alpha\in (0, 1)$. 

  Since the groups in this class do not usually act on rooted trees, the classical methods of \emph{strong contraction} (e.g., used to prove that the Grigorchuk group has intermediate growth) cannot be applied. We replace it by the condition of having finitely many ``incompressible'' words in the group. This condition implies strong contraction for groups acting on rooted trees. The language of incompressible words was studied in the context of groups acting on rooted trees in the papers~\cite{KBRP06,lotsofauth12,Fran}.

  The general condition of finiteness of the language of incopressible words is defined, following~\cite{BNZ}, using travels in the orbital graphas and ``traverses'' of the associated tiles. It is not clear at the moment what is the corresponding notion of strong contraction associated with this condition in general.

  The main result is as follows (Theorem \ref{main}).
 \begin{thm}
    Let $\mathsf{B}=(V,E,\mathbf{s},\mathbf{r})$ be a simple and stationary Bratteli diagram. Let $G$ be a group of bounded type of homeomorphisms of $\Omega(\mathsf{B})$ described by a stationary automaton $\mathcal{A}$. Suppose $G$ acts minimally on $\Omega(\mathsf{B})$. If the set of incompressible elements of $G$ is finite, then the growth function $\gamma_G(R)\preccurlyeq \exp(R^{\alpha})$  for some $\alpha\in (0,1)$.
 \end{thm}

The Grigorchuk group is an example of a self-similar group (acting by tree automorphisms on a regular rooted tree). In \cite{gri84}, R. Grigorchuk proved that the growth function $\gamma_G(R)$ of his group satisfies \[\exp(\sqrt{R})\preccurlyeq \gamma_G(R)\preccurlyeq\exp(R^{\alpha})\addtag\label{bounddd}\] where $\alpha=\log_{32}(31)\approx 0.991$. Later in 1998, a sharper upper bound was given by L. Bartholdi in \cite{Bar98}. Bartholdi proved that the upper estimate for power $\alpha$ equal to $\alpha_0=\log(2)/\log(2/\eta)\approx 0.767$ where $\eta\approx 0.811$ is the positive real root of the polynomial $x^3+x^2+x-2$. In 2018, A. Erschler and T. Zheng \cite{EZ18} proved that the upper bound improved by Bartholdi is actually tight. More precisely, for all $\epsilon>0$ \[\exp(R^{\alpha_0-\epsilon})\preccurlyeq\gamma_G(R)\preccurlyeq\exp(R^{\alpha_0}\addtag)\] where $\alpha_0$ is as in the Bartholdi's estimate. 

Since the original Grigorchuk group, more groups of intermediate growth were constructed. These groups were mainly variants of Grigorchuk's group. See for instance \cite{BS01,BE10,BE11}. In particular, R. Grigorchuk extended his results by constructing uncountably many groups of intermediate growth. See again \cite{gri84} and \cite[Chapters 10,11]{man11}. 

In the paper \cite{nekrash18}, Nekrashevych expanded the class of groups of intermediate growth by giving the first examples of simple groups of intermediate growth (which are also periodic). These groups are obtained by \textit{fragmenting} non-free minimal actions of infinite dihedral groups $D_{\infty}$. 

The original growth estimate in \cite{nekrash18} was larger than $\exp(CR^\alpha)$, but in~\cite{BNZ}, L.~Bartholdi, V.~Nekrashevych and T.~Zheng introduced a new method of counting ``traverses'' in the orbital graphs (summarized here in Subsection \ref{subestivia} and Proposition \ref{ultesti}), and improved the result of  \cite{nekrash18} by showing that the groups from it have growth bounded by $\exp(CR^\alpha)$.

The methods of~\cite{BNZ} will be the main tool of our paper. Our aim is to describe a wide class of examples for which their method of ``traverses'' can be applied.

The groups described in the papers~\cite{nekrash18} and~\cite{BNZ} have linear orbital graphs, i.e., the orbital graphs are quasi-isometric to the real line (or half-line). This is because they are obtained from fragmentations of $D_{\infty}$ whose orbital graphs are linear, while fragmentations preserve quasi-isometry class of orbital graphs. In 2019, J.~Cantu \cite{JC19} generalized the notion of fragmentations in \cite{nekrash18} to study periodicity of groups with non-linear orbital graphs. The class of groups we define in this paper is based on the groups described by J.~Cantu, and we show that many of the groups described by him are of intermediate growth.

\subsection{Sketch of the method}
 We generalize the method developed in \cite{BNZ} to groups with non-linear orbital graphs. To do so, we first define \textit{groups of bounded type} (acting on the spaces of paths of simple and stationary Bratteli diagrams) from \textit{tile inflation processes}. It is a generalization of groups defined by bounded automata (acting on regular rooted trees). See \cite{sidki000} and \cite{Bon}. In the cases when the infinite orbital graphs obtained from embeddings of finite tiles are \textit{well-labeled}, the tile inflation process will define an inverse semigroup $G_0$ acting by partially defined transformations on $\Omega(\mathsf{B})$ (the space of infinite paths of $\mathsf{B}$).
 
 A tile inflation process determines a time-varying automaton $\mathcal{A}$ describing the inverse semigroup $G_0$. The automaton $\mathcal{A}$ also defines a family of inverse semigroups whose initial states generate $G_0$. In certain nice cases, we can define a \textit{group of bounded type} from $G_0$ by relabeling the generators of $G_0$ such that the infinite orbital graphs are \textit{perfectly labeled}. We assume that $G$ acts minimally on $\Omega(\mathsf{B})$. 

 For an inverse semigroup of bounded type $G_0$ (as well as a group of bounded type $G$), we can define the notions of \textit{return} and \textit{incompressible elements}. Using these notions, we can bound the growth function of the group $G$. Let $\mathcal{T}$ be a finite tile. We define \textit{traverses} of $\mathcal{T}$ which are \textit{trajectories} of (semi)group elements on $\mathcal{T}$ starting and ending at boundary points of $\mathcal{T}$ while all the edges in the middle do not touch boundary points. By the process of tile inflations and minimality, a traverse of a tile on a deeper level of $\mathsf{B}$ induces traverses of tiles on earlier levels. This gives rise to an \textbf{injective} map from the set of traverses of tiles on level $m$ to that of the tiles on level $n$ for some $n<m$. This map is called the \textit{last-moment map}. We show that if the set of incompressible elements of a group of bounded type with finite cycles is finite, i.e., every group element that is long enough contains a return, then the last moment map is non-surjective for all $m$ and $n=m-l$ for a fixed $l$. Then we can apply Proposition \ref{ultesti} (a modification of \cite[Proposition 2.7]{BNZ}) to obtain a subexponential growth estimate.
 
\subsection{Organization of the paper}
Section \ref{prelim} contains all the notions that will be used in the paper. Subsection \ref{tileinf} is the ingredient to define groups of bounded type, and Subsection \ref{subestivia} is the tool for obtaining subexponential growth estimates. Section \ref{inct} is a systematic treatment of inverse semigroups of bounded type. We will give the definitions of groups of bounded type and incompressible elements. Then we show polynomial growth and repetitivity of orbital (tile) graphs. Polynomial growth will follow from contractions and linear repetitivity will follow from some other assumptions. Section \ref{mmain} gives the proof of the main theorem. The notion of the last-moment map is defined in the proof. The main theorem is applied in Section \ref{exmmmmp} to provide examples of groups with bounded power in the exponent. 

\subsection*{Acknowledgements}
I would like to express my sincere gratitude to Volodymyr Nekrashevych for his patient guidance through the completion of this work. I would also like to thank Jintao Deng, Ning Ning, and Josiah Owens for valuable suggestions on improving the text. I am deeply grateful for the support from Diego Martinez and the beautiful program \cite{Mar22} he wrote for me. 

\section{Preliminaries}\label{prelim}
In this section (except for the last subsection), we use left action notations (and right action notations for the last subsection) to maintain consistency with the literature.  

\subsection{Tile inflations} \label{tileinf}

Tile inflations are used to define inverse semigroups acting on Cantor sets. In certain nice cases, these semigroups can be turned into groups, which we will discuss in Section \ref{inct}. For visualization of examples of tile inflations, the reader may refer to the Python package \cite{Mar22}.

We start with the definition of \textit{graphs with boundary}. The definition is similar to Definition 1 in  \cite[Chapter 2]{serre} with some additional structure. 

\begin{defn}\label{graphs}
A \textit{graph with boundary} $\Gamma$ consists of a set of \textit{vertices} $V(\Gamma)$, a set of \textit{edges} $E(\Gamma)$,  two partially defined maps
\[\mathsf{s}:E(\Gamma)\dashrightarrow V(\Gamma),\] 
\[\mathsf{r}:E(\Gamma)\dashrightarrow V(\Gamma),\]
and a map
\[E(\Gamma)\rightarrow E(\Gamma),\text{    } e\mapsto e^{-1},\] which satisfy the conditions: 
\begin{enumerate}
\item each $e\in E(\Gamma)$ belongs to at least one of the domains of $\mathsf{s}$ and $\mathsf{r}$; 
\item for each $e\in E(\Gamma)$ we have $(e^{-1})^{-1}=e$, $e^{-1}\neq e$ and $\mathsf{s}(e)=\mathsf{r}(e^{-1})$ or $\mathsf{r}(e)=\mathsf{s}(e^{-1})$. 
\end{enumerate}
The maps $\mathsf{s},\mathsf{r}$ are called respectively \textit{source} and \textit{range} maps, which are defined on subsets of $E(\Gamma)$. For each $e\in E(\Gamma)$, the edge $e^{-1}$ is called the \textit{inverse} edge of $e$. 
\end{defn}

We interpret an edge $e$ as an arrow from $v_0=\mathsf{s}(e)$ to $v_1=\mathsf{r}(e)$ and sometimes denote them as $v_0\rightarrow v_1$.

Informally, we allow some of the edges to be ``hanging'' with source or range to be not in the set of vertices.

Let $\Gamma$ be a graph with boundary. Let us define a class of edge-labeling.

Let $\mathcal{S}$ be a finite set of labels. We assume that $\mathcal{S}$ is symmetric, i.e., there is an involution map \[\mathcal{S}\rightarrow\mathcal{S}, \textup{ } F\mapsto F^{-1}\] mapping each $F\in\mathcal{S}$ to a unique $F^{-1}\in\mathcal{S}$ with $(F^{-1})^{-1}=F$. 

Each element in the set of edges $E(\Gamma)$ is labeled by an element of $\mathcal{S}$, where if $e\in E(\Gamma)$ is labeled by $F\in\mathcal{S}$, then $e^{-1}$ is labeled by the unique $F^{-1}\in\mathcal{S}$ such that one of the following cases will happen.
\begin{enumerate}
    \item There exist $\gamma_1,\gamma_2\in V(\Gamma)$ connected by an arrow (edge) $\gamma_1\rightarrow\gamma_2$ labeled by $F$, and by an arrow $\gamma_1\leftarrow\gamma_2$ labeled by $F^{-1}$ while no other edges connecting $\gamma_1,\gamma_2$ (if exist) are labeled by $F$ or $F^{-1}$. 
    \item\label{bdd} There exists $\gamma\in V(\Gamma)$ and $e\in E(\Gamma)$ such that $\mathsf{s}(e)=\gamma$ but $e$ is not in the domain of $\mathsf{r}$ and $e$ is labeled by $F$. Correspondingly, $e^{-1}$, labeled by $F^{-1}$, is such that $\mathsf{r}(e^{-1})=\gamma$ and it is not in the domain of $\mathsf{s}$. In other words, there exists a $\gamma\in V(\Gamma)$
    to which an outgoing edge labeled by $F$ and an incoming edge labeled by $F^{-1}$ are attached, but $\gamma$ is not connected to other vertices by these two edges. 
\end{enumerate}
If both 1 and 2 happen for the same $F$, then $\gamma_1\neq \gamma$.

\begin{defn}\label{bddpt}
   The point $\gamma$ in Case \ref{bdd} above is called a \textit{boundary vertex} of a graph. The edges labeled by $F$ and $F^{-1}$ at $\gamma$ are called \textit{boundary edges}. 
\end{defn} 

\begin{defn}\label{welllabel}
    Let $\Gamma$ be a graph whose edges are labeled by elements of the set $\mathcal{S}$. A vertex $v$ of $\Gamma$ is called \textit{well-labeled} if, for every $F\in \mathcal{S}$, there is at most one edge starting in $v$ labeled by $F$ and at most one edge ending in $v$ labeled by $F$. The graph $\Gamma$ is said to be \textit{well-labeled} if each of its vertices is well-labeled. It is said to be \textit{perfectly-labeled} if for every vertex $v$ and every $F\in \mathcal{S}$ there is exactly one edge starting in $v$ labeled by $F$ and exactly one edge ending in $v$ labeled by $F$. 
\end{defn}

\begin{defn}\label{brat} A \textit{Bratteli diagram} $\mathsf{B}$ consists of sequences $(V_1,V_2,\ldots)$ and $(E_1,E_2,\ldots)$ of finite sets and sequences of maps $\mathbf{s}_n:E_n\rightarrow V_n$ and $\mathbf{r}_n:E_n\rightarrow V_{n+1}$. The sets $\bigsqcup_{n\geq 1}V_n$ and $\bigsqcup_{n\geq 1}E_n$ are called, respectively, \textit{vertices} and \textit{edges} of the diagram. The maps $\mathbf{s}_n$ and $\mathbf{r}_n$ are called \textit{source map} and \textit{range map}, respectively, and they are assumed to be surjective. We will write $\mathbf{s}$ and $\mathbf{r}$ if no ambiguity would be caused. Denote $\mathsf{B}=((V_n)_{n=1}^{\infty},(E_n)_{n=1}^{\infty},\mathbf{s},\mathbf{r})$. Also denote by $\Omega(\mathsf{B})$ the space of all infinite paths of $\mathsf{B}$ starting in $V_1$, by $\Omega_n(\mathsf{B})$ the space of all paths starting in $V_1$ ending in $V_{n+1}$, and by $\Omega^*(\mathsf{B}):=\bigcup_{n=1}^{\infty}\Omega_n(\mathsf{B})$ the space of all finite paths. Here a \emph{path} in the diagram is a sequence $(e_1, e_2, \ldots, e_n)$ of edges such that $\mathsf{r}(e_i)=\mathsf{s}(e_{i+1})$.
\end{defn}
The space $\Omega(\mathsf{B})$ is a closed subset with the subset topology of the direct product $\prod\limits_{n=1}^{\infty}E_n$. It is compact, totally disconnected, and metrizable.

A \textit{tile} $\mathcal{T}_{v,n}$ is a graph with boundary defined by a \textbf{tile inflation} process, described below, where $V(\mathcal{T}_{v,n})$ consists of all paths of length $n$ of $\mathsf{B}$ ending in the same vertex $v\in E_{n+1}$.

Let us describe this inflation process.
%which is a generalization of \cite[Definition 23]{Bon11}.
Suppose  all tiles up to the $n$-th level of $\mathsf{B}$ are constructed. 

\begin{defn}\label{admissible}
  Let $\mathcal{T}_{1},\mathcal{T}_2$ be tiles (possibly $\mathcal{T}_1=\mathcal{T}_2$) on the $n$-th level of $\mathsf{B}$. A pair of boundary edges $e_1, e_2$ of $\mathcal{T}_1$, $\mathcal{T}_2$, respectively are said to be \textit{compatible} if  both edges are labeled by the same label $F$ and either $\mathsf{s}(e_1)$ and $\mathsf{r}(e_2)$ are undefined, or vice versa. The corresponding vertices of the edges $e_1$ and $e_2$ are called \emph{admissible}.
\end{defn}
\begin{exmp}
  If each $F=F^{-1}$, then each finite boundary point is admissible to itself. 
\end{exmp}

Let $v\in V_{n+2}$ and $\mathcal{T}_v$ be an $(n+1)$-st level tile. It is obtained as follows. 

We choose a set  $P_{n+1}$ (called the set of \textit{connectors} of level $n+1$) of triples 
$(\gamma_1e_1,\gamma_2e_2, F)$ where vertices $\gamma_1,\gamma_2$ of $n$th level tile are such that are admissible boundary points with corresponding compatible edges labeled by $F$ (the edges are unique, by our condition that the graphs are well labeled). The points $\gamma_1e_1$ and $\gamma_2e_2$ are called \textit{connecting points}. 

For each edge $e\in \mathbf{r}^{-1}_{n+1}(v)$, take a copy, denoted $\mathcal{T}_{e}$, of the tile $\mathcal{T}_{\mathbf{s}_{n+1}(e)}$. 

Consider a pair of edges $e_1,e_2\in\mathbf{r}_{n+1}^{-1}(v)$, and suppose that $\gamma_1$ and $\gamma_2$ are boundary points of $\mathcal{T}_{\mathbf{s}_{n+1}(e_1)}$ 
and $\mathcal{T}_{\mathbf{s}_{n+1}(e_2)}$ respectively, such that $(\gamma_1e_1,\gamma_2e_2, F)\in P_{n+1}$. 

Connect then the vertex $\gamma_1e_1$ of $\mathcal{T}_{e_1}$ to the vertex $\gamma_2e_2$ of $\mathcal{T}_{e_2}$ by a pair of arrows $\gamma_1e_1\rightarrow\gamma_2e_2$ and $\gamma_1e_1\leftarrow\gamma_2e_2$ labeled by $F$ and $F^{-1}$, respectively, or labeled by $F^{-1}$ and $F$, respectively. 

The tile $\mathcal{T}_v$ is obtained by applying the above process to all connectors in $P_{n+1}$. 
\begin{rmk}
    Note that a boundary point $\gamma$ on level $n$ (i.e., a path in the Bratteli diagram) must be continued to either a boundary point $\gamma e$ or a connecting point $\gamma e'$ on level $n+1$, or both happen. Meanwhile, non-boundary points cannot be continued to boundary points or connecting points. 
\end{rmk}
  Let $\gamma\in\Omega(\mathsf{B})$. Write $\gamma=e_1e_2\ldots e_n\ldots$. Denote by $\gamma_n=e_1\ldots e_n$ its $n$-th \textit{truncation} ending in a vertex $v_{n+1}\in V_{n+1}$. 
  \begin{defn}
      The \textit{infinite tile} of $\gamma$, denoted $\mathcal{T}_{\gamma}$, is the inductive limit of the embeddings of graphs
\begin{center}
    \begin{tikzcd}
\mathcal{T}_{v_2} \arrow[r, "\phi_{v_2,e_2}"] & \mathcal{T}_{v_3} \arrow [r, "\phi_{v_3,e_3}"] & \ldots \arrow [r, "\phi_{v_{n-1},e_{n-1}}"]& \mathcal{T}_{v_{n}} \arrow [r, "\phi_{v_n,e_{n}}"] & \ldots,
\end{tikzcd}
\end{center}
where $\phi_{v_i,e_{i}}:\eta\mapsto \eta e_{i}$, for $\eta\in \Omega_{i}(\mathsf{B})$ ending in $v_{i}$, is an ismorphism of $\mathcal{T}_{v_i}$ with a subgraph of $\mathcal{T}_{v_{i+1}}$. 
\end{defn}
\begin{defn}
    Let $\mathcal{T}$ be an infinite tile. An infinite path $\xi\in\mathcal{T}$ is said to be a \textit{boundary point of $\mathcal{T}$} if all its $n$-th truncations are boundary vertices in the finite tiles, for $n\in\mathbb{N}$.
\end{defn}

\begin{defn}
    Let $\textup{X}$ be a topological space. An \textit{inverse semigroup} acting on $\textup{X}$ by partial  homeomorphisms, denoted $G\acts \textup{X}$, is a collection of homeomorphisms $g:\textup{Dom}(g)\rightarrow\textup{Ran}(g)$ between open subsets of $\textup{X}$ closed under composition and taking inverses.
\end{defn}

The labels of the edges of infinite tiles define bijections between subsets of $\Omega(\mathsf{B})$, as it follows from the construction of tile inflation and Definition \ref{welllabel}). Namely, for every label $F\in\mathcal{S}$ and an infinite path $\gamma\in\Omega(\mathsf{B})$ we define $F(\gamma)$ as $\mathsf{r}(e)$, where the edge $e$ of $\mathcal{T}_\gamma$ is such that $\mathsf{s}(e)=\gamma$ and $e$ is labeled by $F$. Note that such a (non-boundary) edge may not exist. On the other hand, if it exists, then it exists for all paths $\gamma'\in\Omega(\mathsf{B})$ that have a sufficiently long common beginning with $\gamma$. It follows that the domain of $F$ is open. It follows that every tile inflation rule defines a finite collection of partial homeomorphisms of $\Omega(\mathsf{B})$, hence generates an inverse semigroup.

On the other hand, the union of the domain of $F$ with the set of boundary points of infinite tiles that have a boundary edge labeled by $F$ is clopen.

In some cases the local homeomorphism of $\Omega(\mathsf{B})$ can be extended to a homeomorphism of the closure of the domain. In other words, we can  connect the corresponding infinite tiles at the boundary points by the arrows labeled by $F$ and $F^{-1}$ to form a new well labeled graph. Then each $F\in\mathcal{S}$ defines a transformation (partial homeomorphism) with a \textbf{clopen domain} of $\Omega(\mathsf{B})$, denoted $F:\textup{Dom}(F)\rightarrow\textup{Ran}(F)$. Since each $F\in\mathcal{S}$ corresponds to a unique $F^{-1}$, the set $\mathcal{S}$ actually generates an inverse semigroup, denoted $G=\langle \mathcal{S}\rangle$. 

The following is an example of tile inflations that define an inverse semigroup (in fact, group) action on $\Omega(\mathsf{B})$. 
\begin{exmp}\label{addingmach}
    Consider the Bratteli diagram $\mathsf{B}$ in Figure \ref{biroot}. The space $\Omega^*(\mathsf{B})$ is a binary rooted tree whose vertices are identified with finite sequences of edges $0,1$. Since there is only 1 vertex on each level of $\mathsf{B}$, there is only 1 tile on each level, denoted $\mathcal{T}_n$. Let $\mathcal{S}=\{a,a^{-1}\}$. Consider the following tile inflation process. Each $\mathcal{T}_n$ is a chain of $2^n$ vertices corresponding to paths \[00\ldots 0,10\ldots 0,01\ldots 0,\ldots ,11\ldots 1\] where each vertex is connected to the right by an arrow labeled by $a$ and to the left by an arrow labeled by $a^{-1}$. The boundary points are $00\ldots 0$ and $11\ldots 1$, and they are admissible (See Definition \ref{admissible}). Indeed, there is an arrow (boundary edge) going out of $00\ldots 0$ labeled by $a^{-1}$ and an arrow going in labeled by $a$ while there is an arrow going out of $11\ldots 1$ labeled by $a$ and an arrow going in labeled by $a^{-1}$. The tile $\mathcal{T}_{n+1}$ is constructed as follows. Let $(11\ldots 10,00\ldots 01,a)\in P_{n+1}$. Take two copies of $\mathcal{T}_n$, one appending $0$ to the end of each vertex and the other appending $1$, denoted respectively $\mathcal{T}_n0$ and $\mathcal{T}_n1$. Connect the vertex $11\ldots 10$ on $\mathcal{T}_n0$ to $00\ldots 01$ on $\mathcal{T}_n1$ by an arrow labeled by $a$ and by a backward arrow labeled by $a^{-1}$. Then we get the tile $\mathcal{T}_{n+1}$ whose boundary points are $00\ldots 00$ and $11\ldots 11$ with boundary edges described the same way as for $\mathcal{T}_n$. By the construction, the infinite paths $0^{\omega}=00..$ and $1^{\omega}=11\ldots $ are boundary points of infinite tiles, with boundary edges described the same way as their $n$-th truncations. Denote these two tiles by $\mathcal{T}_{0^{\omega}}$ and $\mathcal{T}_{1^{\omega}}$. Connect $1^{\omega}$ to $0^{\omega}$ by $a$ and $0^{\omega}$ to $1^{\omega}$ by $a^{-1}$. Denote the new graph by $\Gamma_{1^{\omega}}$. It follows that all infinite graphs are well-labeled (in fact, perfectly labeled). Hence the tile inflation described above defines a $\mathbb{Z}$ action on $\Omega(\mathsf{B})$ called the \textit{binary adding machine}. 
    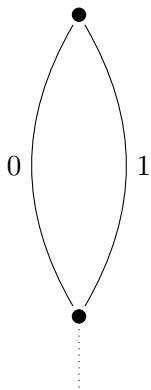
\begin{figure}[!htb]
\centering
\begin{tikzpicture}
\begin{scope}state without output/.append style={draw=none}%[every node/.style={circle,thick,draw}]
    \node (A) at (0,0) {};
    %\node (B) at (4,0) {};
    \node (C) at (0,-4) {};
    %\node (D) at (4,-5) {};

\end{scope}

%\begin{scope} %[>={Stealth[black]},
              every node/.style={fill=white,circle},
              every edge/.style={draw=black,very thick}]
              
    %\path [->] (A) edge [bend left=40] node {} (B);
    %\path [->] (B) edge [bend left=40] node {} (A);
    %\draw [->] (B) to  [out=320,in=40,looseness=9] (B);
    %\path [->] (D) edge node {$3$} (C);
    %\path [->] (A) edge node {$3$} (E);
    %\path [->] (D) edge node {$3$} (E);
    %\path [->] (D) edge node {$3$} (F);
    %\path [->] (C) edge node {$5$} (F);
    %\path [->] (E) edge node {$8$} (F); 
    %\path [->] (B) edge[bend right=60] node {$1$} (E); 

    %\draw[]  (A) node[draw=none][midway,above] {$b_1$} (B);

%\end{scope}

\begin{scope}
\path [-] (A) edge [bend right] node [midway,left] {$0$} (C);
\path [-] (A) edge [bend left] node [midway,right] {$1$} (C);
%\path [-] (D) edge node [near end,right]{$e_1$} (A);
 \draw (0,-4) node {} -- (0,-5) [dotted] node {};
  %\draw (4,-5) node {} -- (4,-6) [dotted] node {};

%\path [-] (B) edge [bend right=35] node  [midway,left]{$1$} (D) ;
%\path [-] (B) edge [bend right=20] node [midway,left]{$2$} (D);
%\path [-] (B) edge [bend right=5] node [midway,left]{$3$} (D);
%\path [-] (B) edge [bend left=5] node [midway,right]{$4$} (D);
%\path [-] (B) edge [bend left=20] node [midway,right]{$5$} (D);
%\path [-] (B) edge [bend left=35] node [midway,right]{$6$} (D);

\filldraw[black] (0,0) circle  (2.5pt) node[anchor=west] {};
%\filldraw[black] (4,0) circle (2.5pt) node[anchor=west] {2};
\filldraw[black] (0,-4) circle (2.5pt) node[anchor=west] {};
%\filldraw[black] (4,-5) circle (2.5pt) node[anchor=west] {2};
\end{scope}

\end{tikzpicture}
\caption{Bratteli Diagram $\mathsf{B}$ for binary rooted tree.}
\label{biroot}
\end{figure}
\end{exmp}

The following is an example of tile inflations that do not define an inverse semigroup action (with clopen domains). 

\begin{exmp}
    Let $\mathsf{X}=\{1,2\}$ and $(\mathsf{X}^{\omega},\sigma)$ be the full shift space. Let $(\mathcal{F},\sigma)$, where $\mathcal{F}\subset \mathsf{X}^{\omega}$ is defined by prohibiting the word $22$, be the \textit{Fibonacci subshift}. Figure \ref{fibosft} shows the Bratteli diagram $\mathsf{B}$ that models this subshift. The edge $e$ corresponds to the word $11$, $e_1$ corresponds to $12$, and $e_2$ corresponds to $21$. Consider the following tile inflation process. On each level of $\mathsf{B}$ there are two tiles, denoted $\mathcal{T}_{1,n}$ and $\mathcal{T}_{2,n}$, corresponding to $1,2$. Let $\mu=e^{\omega}$, $\sigma=(e_1e_2)^{\omega}$ and $\lambda=(e_2e_1)^{\omega}$. Let $\mu_n,\sigma_n,\lambda_n$ be their $n$-th truncations, respectively. There are 2 boundary points on $\mathcal{T}_{1,n}$, for $n\geq 1$, and 2 boundary points on $\mathcal{T}_{2,n}$, for $n\geq 2$. The boundary points on $\mathcal{T}_{1,n}$ are $\mu_n$ and $\lambda_n$ if $n$ is odd, and $\mu_n$ and $\sigma_n$ if $n$ is even. The boundary points on $\mathcal{T}_{2,n}$ are $\mu_{n-1}e_1$ and $\sigma_n$ if $n$ is odd, and $\mu_{n-1}e_1$ and $\lambda_n$ if $n$ is even. Let $\mathcal{S}=\{a,a^{-1}\}$. The following are tiles on levels $1$ and $2$, where only active edges are displayed, and each edge is a two-sided arrow labeled by $a$ from left to right and by $a^{-1}$ from right to left. We omit the labels.  
\begin{align*}
    &\mathcal{T}_{1,1}:\text{  }e\frac{ \ \quad  \quad \ }{}e_2,\text{       }\mathcal{T}_{1,2}:\text{  }e_1,\\
     \mathcal{T}_{2,1}:\text{  }ee&\frac{ \ \quad  \quad \ }{}e_2e\frac{ \ \quad  \quad \ }{}e_1e_2,\text{       }\mathcal{T}_{2,2}:\text{  }ee_1\frac{ \ \quad  \quad \ }{}e_2e_1.
\end{align*}
For $n\geq 3$, the tiles are defined by the following rules. The tile $\mathcal{T}_{1,n}$ is obtained by taking $1$ isomorphic copy $\mathcal{T}_{1,n-1}e$ of $\mathcal{T}_{1,n-1}$ and $1$ isomorphic copy $\mathcal{T}_{2,n-1}e_2$ of $\mathcal{T}_{2,n-1}$, and connecting $\sigma_{n-1}$ on $\mathcal{T}_{1,n-1}$ to $e^{n-2}e_1e_2$ on $\mathcal{T}_{2,n-1}$ by $a$ and $a^{-1}$, if $n$ is odd; connecting $\lambda_{n-1}$ on $\mathcal{T}_{1,n-1}$ to $e^{n-2}e_1e_2$ on $\mathcal{T}_{2,n-1}$ by $a$ and $a^{-1}$, if $n$ is even.  
%\begin{align*}
 %   \mathcal{T}_{1,n}=e\mathcal{T}_{1,n-1}\frac{ \ \quad x_n \quad \ }{} e_2\mathcal{T}_{2,n-1}. 
%\end{align*}
The tile $\mathcal{T}_{2,n}$ is an isomorphic copy $\mathcal{T}_{1,n-1}e_1$ of $\mathcal{T}_{1,n-1}$. Note that all the tiles are linear. 

By the construction, the infinite paths $\mu,\sigma,\lambda$ are boundary points on infinite tiles. We can extend $a$ continuously on $\sigma$ and $\lambda$ by letting $a(\sigma)=a(\lambda)=\mu$. But then $a$ is not injective, so we do not get an action of an inverse semigroup. Hence the above tile inflations do not define an inverse semigroup action on $\Omega(\mathsf{B})$. The transformation $a$ is called the \textit{adic transformation} (or the \textit{Vershik map}) on $\mathsf{B}$. See \cite[Subsection 1.3.4]{nek22}. 

    \begin{figure}[!htb]
\centering
\begin{tikzpicture}
\begin{scope}state without output/.append style={draw=none}%[every node/.style={circle,thick,draw}]
    \node (A) at (0,0) {};
    \node (B) at (4,0) {};
    \node (C) at (0,-5) {};
    \node (D) at (4,-5) {};

\end{scope}

%\begin{scope} %[>={Stealth[black]},
              every node/.style={fill=white,circle},
              every edge/.style={draw=black,very thick}]
              
    %\path [->] (A) edge [bend left=40] node {} (B);
    %\path [->] (B) edge [bend left=40] node {} (A);
    %\draw [->] (B) to  [out=320,in=40,looseness=9] (B);
    %\path [->] (D) edge node {$3$} (C);
    %\path [->] (A) edge node {$3$} (E);
    %\path [->] (D) edge node {$3$} (E);
    %\path [->] (D) edge node {$3$} (F);
    %\path [->] (C) edge node {$5$} (F);
    %\path [->] (E) edge node {$8$} (F); 
    %\path [->] (B) edge[bend right=60] node {$1$} (E); 

    %\draw[]  (A) node[draw=none][midway,above] {$b_1$} (B);

%\end{scope}

\begin{scope}
\path [-] (A) edge node [midway,left] {$e$} (C);
\path [-] (B) edge node [near end,left] {$e_2$} (C);
\path [-] (D) edge node [near end,right]{$e_1$} (A);
 \draw (0,-5) node {} -- (0,-6) [dotted] node {};
  \draw (4,-5) node {} -- (4,-6) [dotted] node {};

%\path [-] (B) edge [bend right=35] node  [midway,left]{$1$} (D) ;
%\path [-] (B) edge [bend right=20] node [midway,left]{$2$} (D);
%\path [-] (B) edge [bend right=5] node [midway,left]{$3$} (D);
%\path [-] (B) edge [bend left=5] node [midway,right]{$4$} (D);
%\path [-] (B) edge [bend left=20] node [midway,right]{$5$} (D);
%\path [-] (B) edge [bend left=35] node [midway,right]{$6$} (D);

\filldraw[black] (0,0) circle  (2.5pt) node[anchor=west] {1};
\filldraw[black] (4,0) circle (2.5pt) node[anchor=west] {2};
\filldraw[black] (0,-5) circle (2.5pt) node[anchor=west] {1};
\filldraw[black] (4,-5) circle (2.5pt) node[anchor=west] {2};
\end{scope}

\end{tikzpicture}
\caption{Bratteli Diagram $\mathsf{B}$ for $(\mathcal{F},\sigma)$.}
\label{fibosft}
\end{figure}
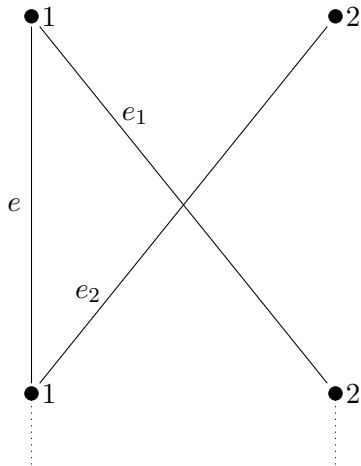
\end{exmp}

The following notion will be useful later. Let $\gamma_1=(e_1,\ldots ,e_n)$ and $\gamma_2=(f_1,\ldots ,f_n)$ be finite paths on $\mathsf{B}$ such that $\mathbf{r}_n(e_n)=\mathbf{r}_n(f_n)$. Define the partial homomorphism (called \textit{prefix switch}) $S_{\gamma_2,\gamma_1}$ by \[S_{\gamma_2,\gamma_1}(e_1,\ldots ,e_n,e_{n+1},e_{n+2},\ldots )=(f_1,\ldots ,f_n,e_{n+1},e_{n+2},\ldots ). \addtag\label{switch}\]
The following is straightforward.
\begin{prop}
    Let $(\gamma_1e_1,\gamma_2e_2,F)\in P_n$ be a connector. Then it defines a pair of prefix switches $S_{\gamma_2e_2,\gamma_1e_1}$ and $S_{\gamma_1e_1,\gamma_2e_2}$. 
\end{prop}

\subsection{Time-varying automata describing inverse semigroup actions}\label{timevary}
Let $G=\langle\mathcal{S}\rangle$ be an inverse semigroup constructed using a tile inflation, as above. It can be described by an $\omega$-deterministic time-varying automaton (Proposition \ref{Automaton}). Recall the definitions. 
\begin{defn}(\cite[Subsections 2.3.3-2.3.5]{nek22}) \label{nondetaut}
     Let $X_1,X_2\ldots $ and $X_1',X_2',\ldots $ be two sequences of \textit{alphabets} and $Q_1,Q_2\ldots $ be a sequence of \textit{sets of states}. A \textit{non-deterministic time-varying automaton} $\mathcal{A}$ consists of a sequence of transitions $T_1,T_2\ldots $ where $T_n\subset Q_n\times Q_{n+1}\times X_n\times X_n'$. The automaton is \textit{$\omega$-deterministic} if for each sequence $x_1x_2\ldots $ where $x_i\in X_i$, and each $q_1\in Q_1$, there exists at most one sequence of transitions of the form $(q_1,q_2,x_1,x_1'),(q_2,q_3,x_2,x_2')\ldots $ where $q_i\in Q_i$ and $x_i'=\lambda_i(q_i,x_i)$. Each initial state $q_1\in Q_1$ defines then a partial transformation $X_1\times X_2\times X_3\times\ldots \rightarrow X_1'\times X_2'\times X_3'\times\ldots $. A state $q\in Q_n$ is called \textit{nontrivial} if it defines a non-identity partial transformation $X_n\times X_{n+1}\times X_{n+2}\times\ldots \rightarrow X_n'\times X_{n+1}'\times X_{n+2}'\times\ldots $.
\end{defn}

\begin{defn}\label{seccttion}
Let $\mathcal{A}$ be a non-deterministic time-varying automaton. Its \textit{Moore diagram} consists of a set of vertices $Q=\bigcup_n Q_n$ and a set of edges $T=\bigcup_n T_n$, where $(q_i,q_{i+1},x,y)\in T_i$ is an arrow from $q_i\in Q_i$ to $q_{i+1}\in Q_{i+1}$ labeled by $x|y$. 
\end{defn}

\begin{defn}\label{finconpt} Let $\mathcal{T}_{v}$ be an $n$-th level tile. An \textit{$n$-th level boundary connection} is a triple $(F,\gamma_1,\gamma_2)$ where $F\in\mathcal{S}$, and $\gamma_1,\gamma_2$ are admissible paths (Definition \ref{admissible}) of length $n$ in the Bratteli diagram such that $F$ is an outgoing boundary edge at $\gamma_1$ and an incoming boundary edge at $\gamma_2$. If $(F,\gamma_1,\gamma_2)$ is a boundary connection, then $(F^{-1},\gamma_2,\gamma_1)$ is also a boundary connection.  % The path $\gamma_1$ is called a \textit{boundary connecting point} (abbr. \textit{connecting point}) of the finite tile $\mathcal{T}_v$. 
\end{defn}
\begin{prop} \label{Automaton}
 Consider the following $\omega$-deterministic time-varying automaton $\mathcal{A}$. The sequence of input-output alphabets is equal to the set of the edges $E_1,E_2,\ldots $ of $\mathsf{B}$. The set of states $Q_n$ is the union of the set of trivial states $1_v$ labeled by the vertices $v\in V_n$ and the set of $n$-th level boundary connections $(F,\gamma_1,\gamma_2)$.  

For every $e\in E_n$, define a transition from the state $1_{\mathbf{s}(e)}$ to $1_{\mathbf{r}(e)}$ labeled by $e|e$.  If $e_1$ and $e_2$ are edges such that $\gamma_1e_1,\gamma_2e_2$ are admissible with $F$ being an outgoing boundary edge at $\gamma_1e_1$ and incoming boundary edge at $\gamma_2e_2$, then for every boundary connection $(F,\gamma_1e_1,\gamma_2e_2)$, define a transition from $(F,\gamma_1,\gamma_2)$ to $(F,\gamma_1e_1,\gamma_2e_2)$ labeled by $e_1|e_2$. Otherwise, define a transition from $(F,\gamma_1,\gamma_2)$ to $1_{\mathbf{r}(e_1)}$ labeled by $e_1|e_2$. 

Then the set of initial states of the form $(F,v_1,v_2)$, for $v_1,v_2\in V_1$, defines the local homeomorphism $F$. Each non-initial non-trivial state has exactly $1$ incoming edge. 
\end{prop}
The proof is the same as that of \cite[Proposition 5.2.19]{nek22}. We omit it here.  

\begin{rmk}\label{identi}
Let $F\in\mathcal{S}$. By Proposition \ref{Automaton}, $F$ is identified with an initial state of $\mathcal{A}$. Let $\gamma\in\Omega_n(\mathsf{B})$. We also identify each state in $F|_{\gamma}$ with a transformation on $\sigma^{n-1}(\Omega(\mathsf{B}))$.
\end{rmk}

% Let us interpret Definitions \ref{seccttion} and \ref{fininf} in the context of the automaton $\mathcal{A}$ in Proposition \ref{Automaton}. 

\begin{defn}
 Let $e\in E$ be such that $\mathbf{s}(e)=v_1$. By Definition \ref{seccttion}, the \textit{section} of $F$ on $e$, denoted ${F|}_{e}$, is the set of states of $\mathcal{A}$ of the form $(F,ee_1,-)$ for some $e_1\in E$ with $\mathbf{s}(e_1)=\mathbf{r}(e)$. Similarly, let $\gamma\in\Omega^*(\mathsf{B})$ on which $F$ is defined. The \textit{section} of $F$ on $\gamma$ is the set of states of the form $(F,\gamma e_1,-)$ for some $e_1\in E$ with $\mathbf{s}(e_1)=\mathbf{r}(\gamma)$. If there exist state(s) $(F,\gamma,-)\in F|_{\eta_n}$ such that $\mathbf{s}(e)\neq\mathbf{r}(\gamma)$, then we define the states $(F,\gamma e,-)$ to be trivial. Let $\gamma\in\Omega(\mathsf{B})$ or $\Omega^*(\mathsf{B})$. The \textit{section of $F$ along $\gamma$} is the sequence of sets $\{F|_{\gamma_n}\}$.  
 \end{defn}
\begin{defn}\label{findir} 
     A state $(F,\gamma_1,\gamma_2)\in\mathcal{A}$ is said to be \textit{finitary} if $(F,\gamma_1\eta_1,\gamma_2\eta_2)$ is a trivial state for all $\eta_1\in\Omega_m(\mathsf{B})$ such that $F$ is defined on $\eta_1\gamma_1$, for some $m$. The state $(F,\gamma_1,\gamma_2)$ is said to be \textit{directed} if there is an $\eta\in\Omega(\mathsf{B})$ such that $(F,\gamma_1\eta_n,\gamma_2\eta_n')$ are all nontrivial, for all finite truncations $\eta_n$ of $\eta$.  Let $F\in\mathcal{S}$. The transformation $F$ is said to be \textit{finitary on some $\gamma\in\Omega_n(\mathsf{B})$} if $F|_{\gamma}$ consists of only finitary states. The transformation $F$ is said to be \textit{finitary} if there exists $n\in\mathbb{N}$ such that for all $\gamma\in\Omega_n(\mathsf{B})$ on which $F$ is defined, $F|_{\gamma}$ consists of only finitary or trivial states. The transformation $F$ is said to be \textit{directed along $\gamma$} if there is a nontrivial directed state in each ${F|}_{\gamma_n}$, for $\gamma\in\Omega(\mathsf{B})$ and all $n\in\mathbb{N}$. 
\end{defn}

We have the following equivalent definitions of boundary points in terms of the automaton $\mathcal{A}$. 

\begin{defn}\label{bdrypts} Let $\mathcal{T}$ be an infinite tile. Let $\gamma\in \mathcal{T}$ be a vertex. The point $\gamma$ is said to be a \textit{boundary point} of $\mathcal{T}$ if there is an $F\in\mathcal{S}$ whose sections along $\gamma$ are all nontrivial. 
\end{defn}

\begin{defn}\label{finbdrypt}
    Let $\mathcal{T}_v$ be a finite tile ending in $v\in V_{n+1}$. A point $\gamma\in\mathcal{T}_v$ is said to be a boundary point of $\mathcal{T}_v$ if there is an $F\in\mathcal{S}$ such that $(F,\gamma,\gamma')$ is a nontrivial state for some $\gamma'\in\Omega_n(\mathsf{B})$ (not necessarily ending in $v$). 
\end{defn}

\subsection{Minimal semigroup actions on $\Omega(\mathsf{B})$}\label{minacto}
Let $\mathcal{X}$ be a Cantor set and $G$ be a semigroup. The action $G\acts \mathcal{X}$ is said to be \textit{minimal} if all $G$-orbits are dense in $\mathcal{X}$. In this subsection, we give a characterization of minimal semigroup actions on Cantor sets (Proposition \ref{minact}) and give two proofs. In order for $\Omega(\mathsf{B})$ to be a Cantor set, i.e., for it to have no isolated points, we need the assumption that $\mathsf{B}$ is \textit{simple}.

\begin{defn} A Bratteli diagram $\mathsf{B}$ is said to be \textit{simple} if for every level $n$ there exists $m>n$ such that for every pair $v\in V_n$ and $u\in V_m$ there exists a path $(e_n,e_{n+1},\ldots ,e_{m-1})$ in $\mathsf{B}$ starting in $v$ and ending in $u$. 
\end{defn}

Now we discuss the characterization of minimal group actions on $\Omega(\mathsf{B})$ in terms of tiles. Let $G$ be an inverse semigroup defined by a tile inflation process.  For any $\alpha\in \Omega(\mathsf{B})$, the set of $G$-orbits of $\alpha$ is denoted $G\alpha$. 

 \begin{defn}\label{cofinal} Any pair of points $\zeta,\eta \in \Omega(\mathsf{B})$ are said to be \textit{cofinal} if they differ for finitely many edges. This defines an equivalence relation on $\Omega(\mathsf{B})$ called an \textit{tail-equivalence} relation. For any point $\eta\in\Omega(\mathsf{B})$, denote by $\text{Cof}(\eta)$ the set of points $\Omega(\mathsf{B})$ that are cofinal with $\eta$. The set $\text{Cof}(\eta)$ is called the \textit{cofinality class} of $\eta$. 
\end{defn}

\begin{prop}\label{minact}  Assume that all finite tiles are connected. Let $\mathsf{B}$ be a simple Bratteli diagram. Let $G\acts \Omega(\mathsf{B})$ by homeomorphisms. Then the action of $G$ is minimal if and only if for any pair of points (finite paths) $\alpha,\beta\in \Omega_n(\mathsf{B})$, there exists $m>n$ (depending on $\alpha,\beta$) such that all vertices on $\mathcal{T}_{\alpha}$ and $\mathcal{T}_{\beta}$ can be continued to $\Omega_m(\mathsf{B})$ and the continuations are in the same tile $\mathcal{T}_{\gamma}$ for some $\gamma\in \Omega_m(\mathsf{B})$. 
\end{prop}
The proof follows from the following results. 
\begin{lemma} \textup{(part of \cite[Theorem 3.11]{putnam})} \label{denseeq} A Bratteli diagram $\mathsf{B}$ is simple if and only if for any point $\eta\in \Omega(\mathsf{B})$, $\textup{Cof}(\eta)$ is dense in $\Omega(\mathsf{B})$. 
\end{lemma}
\begin{prop} \textup{(\cite[Proposition 2.3.6]{JC19})}\label{ctile} If the finite tiles are eventually connected, then $\textup{Cof}(\eta)\subset G\eta$ for all $\eta \in \Omega(\mathsf{B})$.
\end{prop}
\begin{proof} There exists $N\in\mathbb{N}$ such that the tiles $\mathcal{T}_{\gamma}$ are connected for all $\gamma\in\Omega_n(\mathsf{B})$ and $n\geq N$. Let $\eta\in\Omega(\mathsf{B})$ and $\zeta\in \textup{Cof}(\eta)$ be different from $\eta$. Write $\eta=(e_1,e_2,\ldots )$ and $\zeta=(f_1,f_2,\ldots )$. The there exists $m\geq N$ such that $e_k=f_k$ for all $k\geq m$. Let $\eta_m$ and $\zeta_m$ be respectively the $m$-th truncation of $\eta$ and $\zeta$. Then since $\mathcal{T}_{\eta_m}$ is connected, there exists $g\in G$ such that $g(\eta_m)=\zeta_m$ and $g$ acts as the prefix switch $S_{\eta_m',\eta_m}$ at $\eta_m\gamma$ for all allowed continuations $\gamma$. Hence $g(\eta)=\zeta$ and $\zeta\in G\eta$. 
\end{proof}

%\begin{prop}\label{mintile} Let $\mathsf{B}$ and $\eta$ be as above. The action $G\acts \Omega(\mathsf{B})$ is minimal if and only if the finite tiles are eventually connected to each other as isomorphic copies of subgraphs of a tile on deeper levels. 
%\end{prop}
\begin{proof}[Proof of Proposition \ref{minact}] ($\impliedby$)  Let $\eta\in \Omega(\mathsf{B})$ be any point and let $\eta_n$ be its $n$-th truncation, i.e., the first $n$ prefixes of $\eta$. For each $n$ there are finitely many tiles of the form $\mathcal{T}_{\gamma}$ on $\Omega_n(\mathsf{B})$. For each $\mathcal{T}_{\alpha}$, there exists $m>n$ such that there are continuations (depending on $\alpha$) of $\mathcal{T}_{\eta_n}$ and $\mathcal{T}_{\alpha}$ to $\Omega_m(\mathsf{B})$ that are subgraphs of the same $\mathcal{T}_{\gamma}$ in $\Omega_m(\mathsf{B})$. Since this is true for all $n$ and $\mathsf{B}$ is simple, the orbit $G\eta$ is eventually in the neighborhood of any point in $\Omega(\mathsf{B})$. Hence the action is minimal.

($\implies$) Let $\eta\in\Omega(\mathsf{B})$. Since the set of vertices of the tile $\mathcal{T}_{\eta}$ is equal to $\text{Cof}(\eta)\cap G\eta$, by Proposition \ref{ctile}, the vertices of $\mathcal{T}_{\eta}$ can be identified with $\text{Cof}(\eta)$. Since, by Lemma \ref{denseeq}, cofinality classes are dense, it follows that every $G$-orbit is dense. 
\end{proof}

\begin{rmk} For groups acting by homeomorphisms on regular rooted trees, the induced action on the boundary of the tree is minimal if and only if the action is level-transitive. See \cite[Proposition 6.5]{GNS}. If $G$ is a self-replicating self-similar group, any two sequences are on the same orbit if and only if they are cofinal. See the proof of \cite[Lemma 7.2]{vor12}. In our case, the minimal action might not be level-transitive. However, the elements in each orbit can be continued so they are ``connected" to the continuation of other orbits on deeper levels. 
\end{rmk}
\subsection{Graphs of actions and fragmentations}\label{GOA}

In this subsection, $G$ is assumed to be a group. Let $\mathcal{X}$ be a Cantor set, $G$ be a group generated by a finite set $S$ of homeomorphisms of $\mathcal{X}$, and $G\acts \mathcal{X}$ be a \textbf{faithful} minimal action. For $\zeta\in\mathcal{X}$, denote by $G_{\zeta}$ the stabilizer of $\zeta$. Denote by $G_{(\zeta)}$ the \textit{neighborhood stabilizer} of $\zeta$, which is defined to be the subgroup of all elements $g \in G$ such that $\zeta$ is an interior point of the set of fixed points of $g$. In other words, each $g\in G_{(\zeta)}$ fixes pointwise an open neighborhood of $\zeta$.  The \textit{orbital graph} of $\zeta$, denoted $\Gamma_{\zeta}$, has the orbit $G\zeta$ as its set of vertices, and for every $\eta\in G\zeta$ and every $s\in S$, there is an arrow from $\eta$ to $s(\eta)$ labelled by $s$. The graph $\Gamma_{\zeta}$ is naturally isomorphic to the \textit{Schreier graph} of $G$ modulo $G_{\zeta}$.

A \textit{germ} of an element $g\in G$ at $\zeta$ is an equivalence class of the pair $(g,\zeta)$, where two pairs $(g_1,\zeta)$ and $(g_2,\zeta)$ are equivalent if there exists an open neighborhood $U$ of $\zeta$ such that $g_1|_{U}=g_2|_{U}$. The \textit{graph of germs}, denoted $\widetilde{\Gamma}_{\zeta}$, is the Schreier graph of $G$ modulo $G_{(\zeta)}$. The vertices of $\widetilde{\Gamma}_{\zeta}$ are identified with germs of elements of $G$ at $\zeta$. Since $G_{(\zeta)}$ is a normal subgroup of $G_{\zeta}$, the map $hG_{(\zeta)}\mapsto hG_{\zeta}$ induces a Galois covering of graphs $\widetilde{\Gamma}_{\zeta}\rightarrow\Gamma_{\zeta}$ with the group of deck transformations $G_{\zeta}/G_{(\zeta)}$. The group $G_{\zeta}/G_{(\zeta)}$ is called the \textit{group of germs}. Recall some definitions. 

\defn\label{gregular} A point $\zeta\in\mathcal{X}$ is said to be \textit{$G$-regular} if its group of germs $G_{\zeta}/G_{(\zeta)}$ is trivial, i.e., if every $g\in G$ that fixes $\xi$ also pointwise fixes a neighborhood of $\zeta$. If $\zeta$ is not $G$-regular, then we say it is \textit{singular}. 
\defn\label{Germsin} (\cite[Definition 3.2.1]{JC19}) Let $G$ be a group generated by a finite set $S$. Let $\xi\in \mathcal{X}$ be a singular point. It is said to be \textit{germ-defining} if 
\begin{enumerate}
    \item $\xi$ is the only singular point \underline{of $S$} in the orbit $G\xi$,
    \item for any $\zeta\in G\xi$ and an element $g=s_ms_{m-1}\ldots s_{1}\in G_{\zeta},s_i\in S\cup S^{-1}$, such that $s_i\ldots s_1(\zeta)$ are distinct for $i=1,\ldots ,m-1$, $(g,\zeta)=(Id,\zeta)$. The element $g$ is called a \textit{nonloop cycle} at $\zeta$ on $\Gamma_{\xi}$.
\end{enumerate}
\defn\label{pure} Let $\xi\in\mathcal{X}$ be a germ-defining singular point. We say that $\xi$ is a \textit{Hausdorff singularity} if for every $g\in G_{\xi}\backslash G_{(\xi)}$, the interior of the set of fixed points of $g$ does not accumulate on $\xi$. Otherwise, $\xi$ is called a \textit{non-Hausdorff} singularity. If for every $g\in G_{\xi}$, the interior of the set of fixed points of $g$ accumulates on $\xi$, then $\xi$ is called a \textit{purely non-Hausdorff} singularity.  

\begin{rmk} Definition \ref{pure} is a slight modification of \cite[Definition 2.2]{nekrash18} by adding the assumption that $\xi$ is a germ-defining singular point.
\end{rmk}

\begin{defn}\label{frag} (\cite[Definition 4.1.1]{JC19}) Let $h$ be a finite order homeomorphism of a Cantor set $\mathcal{X}$ and $\mathcal{P}_h$ be a finite partition of $\text{Supp}(h)$ into $h$-invariant open sets. A \textit{fragmentation} of $h$ with respect to $\mathcal{P}_h$ is a finite group $F_h$ of homeomorphisms of $\mathcal{X}$ such that
\begin{enumerate}
    \item for all $g\in F_h$ and $P\in\mathcal{P}_h$, there exists $k\in\mathbb{Z}$ such that $g|_P=h^k|_P$,
    \item for all $g\in F_h$, $g|_{\text{Fix}(h)}=Id$,
    \item for all $P\in\mathcal{P}_h$, there exists $g\in F_h$ such that $g|_P=h^k|_P$. 
\end{enumerate}
Let $S$ be a finite generating set of the group $G$. The group generated by all the union of all $F_s$ (including trivial fragmentations which are just $\langle s \rangle$), for $s\in S$, denoted $F_G$, is called a \textit{fragmentation group}. 
\end{defn}
The elements of $\mathcal{P}_h$ are called the \textit{pieces} of the fragmentation. Each $F_h$ is always a finite abelian group. It is convenient to express fragmentations in terms of subdirect products. 
\begin{defn}\label{SUBD} Let $F_h$, $\mathcal{P}_h$ be as above. Each piece $P\in\mathcal{P}_h$ defines an epimorphism $\pi_P:F_h\rightarrow \langle h \rangle|_P$ by $\pi_P(g)=g|_P$. The product map $$(\pi_P)_{P\in\mathcal{P}_h}:F_h\rightarrow \prod_{P\in\mathcal{P}_h}\langle h \rangle|_P$$ is called a \textit{subdirect product} (of $\prod_{P\in\mathcal{P}_h}\langle h \rangle|_P$) or a \textit{subdirect embedding}. It is an embedding that is surjective on each factor. 
\end{defn}

Conversely, given a subdirect product, we can define a fragmentation as follows. Let $\alpha\in \prod_{P\in\mathcal{P}_h}\langle h \rangle|_P$ and $\zeta\in \mathcal{X}$. Define 
\begin{align*}
\alpha(\zeta)=
\begin{cases}
h^k(\zeta) \text{, if }\zeta\in P \text{ and } \alpha_P=h^k|_P,\\
\zeta \text{, otherwise. }
\end{cases}
\end{align*}
This gives a faithful action of $\alpha$ on $\mathcal{X}$. The set of all such $\alpha$ is the fragmentation $F_h$. 

The technique of fragmentations is used to construct groups with a purely non-Hausdorff singularity (as shown in the following proposition), which are the majority of examples in this paper. 

\begin{prop}
Let $(G,S)$ be a finitely-generated group acting on a Cantor set $\mathcal{X}$ as above. Let $\xi\in \mathcal{X}$ be a germ-defining singularity. Let $H$ be the group of germs of $\xi$. Then for any non-identity $h\in H$, there exists a finite partition $\mathcal{P}_h$ of $\mathcal{X}\backslash\{\xi\}$ into disjoint $h$-invariant open subsets such that each piece $P_i\in\mathcal{P}_h$ accumulates on $\xi$. For each partition $P_h$, there exists a fragmentation $F_h$ of $h$ such that $\xi$ is a purely non-Hausdorff singularity of the group $\langle S\backslash\{h\}\cup F_h \rangle$. 
\end{prop}

This follows directly from the results of \cite[Section 4.5]{JC19}.

\subsection{Various notions of repetitivity of orbital graphs}\label{repppttt}
Repetitivity of orbital graphs is a key ingredient to obtain subexponential growth. We collect here various notions of repetitivity that will be useful later. The following was proved in \cite[Proposition 2.5]{nekrash18}. 
\begin{prop}\label{repetitivity}
    Let $\mathcal{X}$ be a Cantor set. Suppose the action of a semigroup $G$ on $\mathcal{X}$ is minimal. Then for every $r>0$, there exists $R(r)>0$ such that for every $G$-regular point $\zeta \in \mathcal{X}$ and for every $\eta \in \mathcal{X}$, there exists a vertex $\eta'$ of $\Gamma_{\eta}$ such that $d(\eta,\eta')\leq R(r)$, and the rooted balls $B_{\zeta}(r)\subset \Gamma_{\zeta}$ and $B_{\eta'}(r)\subset \Gamma_{\eta}$ are isomorphic.
\end{prop}

\begin{defn}\label{repp}
    The action of $G$ on $\mathcal{X}$ is said to be \textit{repetitive} if it satisfies the conclusion of Proposition \ref{repetitivity}. The action is said to be \textit{linearly repetitive} if there exists $K>1$ such that the function $R(r)$ from Proposition \ref{repetitivity} satisfies $R(r)<Kr$ for all $r\geq 1$. 
\end{defn}

We also have an equivalent definition of linear repetitivity. 

\begin{defnis}{repp}
    The orbital $\Gamma$ is \textit{linearly repetitive} if there exists a constant $L>1$ such that for every pair of vertices $v,v'\in\Gamma$ and every $R>0$ there exists an isomorphic embedding of $B_v(R)$ into $B_{v'}(LR)$. 
\end{defnis}

Let $\Delta\subset \Gamma$ be a finite connected subgraph. Recall that the diameter of $\Delta$, denoted Diam$(\Delta)$, is defined to be $\max_{u,v\in V(\Delta)}d(u,v)$, where $d$ is the metric on $\Gamma$ counting the smallest number of paths (the length of a geodesic) between any pair of vertices $u,v$ modulo multiple edges and loops. The distance between $2$ disjoint connected subgraphs $\Delta$ and $\Delta'$ is the shortest distance between the vertices in  $\Delta$ and $\Delta'$.
\begin{prop}
 Let $\Gamma$ be an infinite orbital graph. Then $\Gamma$ is linearly repetitive if and only if for any finite connected subgraph $\Delta\subset \Gamma$, there exists a subgraph $\Delta'$ isomorphic to $\Delta$ such that the distance between $\Delta$ and $\Delta'$ is bounded above by $C\cdot\textup{Diam}(\Delta)$ for a fixed constant $C>0$. 
\end{prop}
\begin{proof}
    $(\implies)$ Suppose $\Gamma$ is linearly repetitive. Let $B_{v}(R)$ be a ball of radius $R$. Then there exists  $v'\in\Gamma$ with $d(v,v')<CR$ such that the ball $B_{v'}(R)$ is isomorphic to $B_v(R)$. Note that $R=\frac{1}{2}\textup{Diam}(B_v(R))$. Let $\Delta$ be any finite subgraph of $\Gamma$. Then $\Delta$ is contained in some $B_v(R)$ for $R=\frac{1}{2}\textup{Diam}(\Delta)$. Hence the necessity is true. 
    
    $(\impliedby)$ Conversely, we can replace $\Delta$ by any ball of radius $R$. Hence the sufficiency is true.  
\end{proof}
The following is a stronger version of repetitivity.
\begin{defn}\label{substrrep}(\cite[Definition 3.3.7]{JC19}) Let $\Delta$ be a finite connected graph and $\Gamma$ be a locally finite connected infinite graph. The graph $\Delta$ is said to be \textit{strongly repetitive in $\Gamma$} if there exists infinitely many isomorphic copies $\{\Delta_i\}_{\i\in I}$ of $\Delta$ in $\Gamma$ such that $\Gamma\backslash E(\bigcup_{i\in I}\Delta_i)$ consists of finite connected components of bounded size. The set $\{\Delta_i\}_{\i\in I}$ is called a \textit{$\Delta$-sieve}.
\end{defn}

\begin{defn}\label{strrep} Let $\Gamma$ be as above. The graph $\Gamma$ is called \textit{strongly repetitive} if for every finite connected subgraph $\Delta$ of $\Gamma$, there exists a $\Delta$-sieve. 
\end{defn}
Notice that in general the orbital graphs are not strongly repetitive. However, it will be useful to assume some finite subgraphs $\Delta$ are strongly repetitive (i.e. they form $\Delta$-sieves) in some special cases. 

\subsection{Subexponential growth estimates via orbital graphs }\label{subestivia}
We use right action notations here. Let $G$ be a group. Let $S$ be a finite symmetric generating set of $G$. Every element $g\in G$ can be written as a product $g=s_1\ldots s_k$ with $s_i\in S$, for $i=1,\ldots ,k$. Actions of $g$ on some elements will be denoted $\zeta\cdot g$ or $\zeta\cdot s_1\ldots s_k$. Let $l\in \mathbb{N}$ be such that $g=s_1,\ldots ,s_l$ is the shortest representation of $g$. Then we can define a length function $|\cdot|:G\rightarrow \mathbb{R}$ by $|g|=l$. The growth function for $G$, denoted $\gamma_G(R)$, is defined to be the number of $g\in G$ with $|g|\leq R$. Let $f_1$ and $f_2$ be two growth functions. The function $f_1$ is said to be \textit{dominated by} $f_2$, denoted $f_1\preccurlyeq f_2$, if there exists $C\geq 1$ such that $f_1(R)\leq f_2(CR)$ for all $R\geq 1$. The function $f_1$ is said to be \textit{equivalent} to $f_2$, denoted $f_1 \sim f_2$, if $f_1\preccurlyeq f_2$ and $f_2\preccurlyeq f_1$. 

\begin{defn}\label{stnpor}{(\cite[Definition 2.2]{BNZ})} Let $w=s_1\ldots s_l$ be a word in $S^*$, and let $v$ be a vertex in $\Gamma$. Denote by $R_w(v)$ the maximum over $0\leq i\leq n$ of the distances $d(v,v\cdot s_1\ldots s_i)$. The \textit{standard portrait} $\mathcal{P}_w$ of $w$ is defined to be the collection of the (isomorphism classes of the) bi-rooted graphs $(B_v(R_w(v)),v,v\cdot w)$. Let $\gamma_{\Gamma}(R)=\max_{v\in\Gamma}\# B_v(R)$ be the maximal size of a ball of radius $R$ in $\Gamma$, and let $\delta_{\Gamma}(R)$ be the number of isomorphism classes of rooted graphs $(B_v(R),v)$. 
\end{defn}

Let $g\in G$. Denote by $L(g)$ the minimal size of standard portraits $\mathcal{P}_w$ among all words $w\in S^*$ representing $g$, and set $L(R)=\max_{|g|\leq R}L(g)$. We have the following estimate for the growth of $G$.
\begin{prop}\label{estiviaorb} \textup{(\cite[Proposition 2.3]{BNZ})}
 The growth function $\gamma_G(R)$ of $G$ satisfies $$\gamma_G(R)\leq (\delta_{\Gamma}(R)\gamma_{\Gamma}(R))^{L(R)}.$$   
\end{prop}

In particular, if $\delta_{\Gamma}(R)$ and $\gamma_{\Gamma}(R)$ are bounded above by polynomials, and $L(R)$ is bounded by $R^{\alpha}$ for some $\alpha \in (0,1)$, then the growth of $G$ is bounded from above
by $\exp(\log R \cdot R^{\alpha})$.

Given a standard portrait $\mathcal{P}_{w}$ for $w\in S^*$, define $$N_p(w)=\left(\sum\limits_{(B,x,y)\in\mathcal{P}_{w}}(\# B)^p \right)^{1/p}. $$ Let $N_p(g)$ be the infimum of $N_p(w)$ over all $w$ representing $g$. Note $N_1(w)=\#\mathcal{P}_w$, so $N_1(g) = L(g)$. The following proposition was proved in \cite[Corollary 2.6]{BNZ}. 
\begin{prop} \label{upperbound}
    Suppose the graph of the action $\Gamma$ has polynomial growth and is linearly repetitive. Denote by $N_p(R)$ the maximum of $N_p(g)$ for all elements $g\in G$ of length at most R. Then for every $p \geq 1$ there exists $K > 0$ such that the growth of $G$ satisfies
$$\gamma_{G}(R) \leq R^{KN_p(R)}.$$
\end{prop}
By the above proposition, we will develop methods to find upper bounds for $N_p(R)$ in order to bound the growth function $\gamma_G(R)$. 

\section{Bounded type}\label{inct}
Throughout the rest of this paper, we use right action notations. Let $G$ be a group or an inverse semigroup. Let $g\in G$, $e\in E$ and $\gamma\in\Omega^*(\mathsf{B})$. 
%Taking sections of $g$ will be denoted ${}_{e}{|g}$ and ${}_{\gamma}{|g}$. 
For $\zeta\in\Omega(\mathsf{B})$, the orbital graph $\Gamma_{\zeta}$ is (isomorphic to) the Schreier graph of the right coset space $G_{\zeta}\backslash G$, and the graph of germs $\widetilde{\Gamma}_{\zeta}$ is the Schreier graph of the right coset space $G_{(\zeta)}\backslash G$. Let $F\in\mathcal{S}$. The germ of $F$ at $\gamma$ is denoted $\gamma F$. The action of $F$ on $\gamma$ will be denoted $\gamma\cdot F$. Actions of (semi)group elements on an infinite path on the Bratteli diagram will be denoted $x\cdot g$ or $\gamma\cdot g$. Paths on the Bratteli diagram are read from the right to the left. An infinite periodic path in $\Omega(\mathsf{B})$ will be denoted $(e_n\ldots e_1)^{\omega}$ or $\gamma^{\omega}$ for some finite path $\gamma$. Given $\gamma=\ldots e_2e_1\in\Omega(\mathsf{B})$, the notation $\gamma_n=e_n\ldots e_2e_1$ means the $n$-th truncation of $\gamma$. All tiles are assumed to be connected.

\subsection{Definitions of inverse semigroups and groups of bounded type}
Let $\mathsf{B}$ and $\mathcal{S}$ be as in Subsection \ref{tileinf}. In this subsection, we define inverse semigroups and groups of bounded type and describe them using automata. In this paper, we will focus on inverse semigroups and groups defined by \textit{stationary} automata. 
\begin{defn}\label{bddtpp}
    A tile inflation process on $\mathsf{B}$ is of \textit{bounded type} if it satisfies the following conditions.
 \begin{enumerate}   
  \item \label{cond1}  All finite tiles have uniformly bounded cardinalities of the set of boundary points and there exist only finitely many boundary points of infinite tiles. 
  \item \label{cond2} All finite tiles are linearly repetitive in all infinite tiles.
\end{enumerate}
    
\end{defn}

\begin{defn}\label{semibd}
    The inverse semigroup generated by $\mathcal{S}$, denoted $G=\langle \mathcal{S}\rangle$, is of \textit{bounded type} if it is defined by a tile inflation process of bounded type.
\end{defn}

In certain nice cases, we have groups determined by tile inflations of bounded type. For example, the tile inflation process directly defines a group in Example \ref{addingmach}, or, after a relabeling of orbital graphs, the inverse semigroup can be turned into a group. See Subsection \ref{2tiles} for an example of the relabeling process. Hence we have the following.

\begin{defn}\label{grpbd}
  Let $S$ be a finite labeling set of a tile inflation process of bounded type. The inverse semigroup $G$ generated by $S$, denoted $G=\langle S\rangle$, is a \textit{group of bounded type} if, after extending the action of some labels at boundary points of infinite tiles, the tile inflation process determines perfectly labeled orbital graphs labeled by elements in $S$. 
\end{defn}
\begin{rmk}
Since we do not always have groups defined by tile inflations, we will focus on the properties of inverse semigroups in this section. 
\end{rmk}

\subsection{Stationary case}
 \begin{defn}\label{stabra} A Bratteli diagram $\mathsf{B}=((V_n)_{n=1}^{\infty},(E_n)_{n=1}^{\infty}),\mathbf{s}_n,\mathbf{r}_n)$ is said to be \textit{stationary} if all the sequences $V_n, E_n, \mathbf{s}_n,\mathbf{r}_n$ are constant. Denote $E_n=E$, $V_n=V$, $\mathbf{s}_n=\mathbf{s}$ and $\mathbf{r}_n=\mathbf{r}$.
\end{defn}

Let $\mathsf{B}$ be stationary. Let $\sigma:\Omega(\mathsf{B})\rightarrow\Omega(\mathsf{B})$ be the shift map, i.e., $\sigma(\ldots e_2e_1)=\ldots e_3e_2$. We also consider $\sigma:\Omega^*(\mathsf{B})\rightarrow\Omega^*(\mathsf{B})$ defined by deleting the first edge on the right of each finite path.  

\begin{defn}\label{staaut}
 Let $\mathcal{A}$ be an automaton described in Proposition~\ref{Automaton}. It is said to be \textit{stationary} if the following hold. 
\begin{enumerate}
 \item The automaton $\mathcal{A}$ is defined on a stationary Bratteli diagram $\mathsf{B}$, meaning that the sequence of alphabets is the constant set $E$. 

 \item Let $B_n,C_n\subset\Omega^*(\mathsf{B})$ be respectively the set of boundary points and the set of connecting points on level $n$ tiles. The map $\sigma$ restricted on each $B_n$ and $C_n$ are well-defined maps between these sets, i.e., $\sigma:B_{n+1}\rightarrow B_n$ and $\sigma:C_{n+1}\rightarrow C_n$, for all $n\in\mathbb{N}$.

 \item There exist finite sets $B,P$ such that there are bijections $\varphi_n:B\rightarrow B_n$ and $\psi_n: P\rightarrow C_n$ satisfying $\sigma\circ\varphi_{n+1}=\varphi_n$ and $\sigma\circ\psi_{n+1}=\psi_n$, for all $n\in\mathbb{N}$.
%  \item The sequence of sets of states is constant, meaning that for all $n\in\mathbb{N}$, there is a exactly $1$ state on level $n$ of $\mathsf{B}$ defining the same transformation on $\sigma^{n-1}(\Omega(\mathsf{B}))=\Omega(\mathsf{B})$. 
%  \item The sequence of the sets of transitions is constant. 
\end{enumerate}
 \end{defn}
 \begin{rmk}
Since the nontrivial states in $\mathcal{A}$ are boundary connections (Definition \ref{finconpt}) $(F,\gamma_1,\gamma_2)$, where $\gamma_1,\gamma_2$ are either admissible boundary points or connecting points, the definition above implies the following.
\begin{enumerate}
\item The sequence of sets of states of $\mathcal{A}$ is constant, meaning that for all $n\in\mathbb{N}$, there is exactly $1$ state on level $n$ of $\mathsf{B}$ defining the same transformation on $\sigma^{n-1}(\Omega(\mathsf{B}))=\Omega(\mathsf{B})$. 

\item The sequence of the sets of transitions of $\mathcal{A}$ is constant. 
\end{enumerate}
The finite sets $B,P$ will be described at the end of this subsection. 
 \end{rmk}

% \begin{defn}\label{stagrp}
%  The tile inflation process and the inverse semigroup or group $G$ are called \textit{stationary} if they are defined by a stationary automaton $\mathcal{A}$.  
% \end{defn}
The elements of the generating set $\mathcal{S}$ are identified with the initial states of $\mathcal{A}$. See Remark \ref{identi}. To describe the inverse semigroup $G=\langle\mathcal{S}\rangle$ in terms of an automaton, we need the notion of the \textit{full automaton} of $\mathcal{A}$, denoted $\mathfrak{F}(\mathcal{A})$, whose set of initial states is $G$ and the sequence of input-output alphabets is equal to that of $\mathcal{A}$. The states and transition functions of $\mathfrak{F}(\mathcal{A})$ are described as follows. 

Since $\mathcal{A}$ is $\omega$-deterministic (Definition \ref{nondetaut}), we have the following notion of taking sections in $\mathcal{A}$ and $\mathfrak{F}(\mathcal{A})$. 

\begin{defn}
Given $\gamma\in\Omega(\mathsf{B})$ accepted by an initial state $F\in\mathcal
{A}$ (note we have identified transformations in $\mathcal{S}$ with states in $\mathcal{A}$), there is a unique choice of state in each set ${}_{\gamma_n}{|F}$ (see Definition \ref{seccttion}) along $\gamma$, denoted ${}_{(\gamma,n)}{|F}$. The \textit{sections of $F$ along $\gamma$} is the sequence $\{{}_{(\gamma,n)}{|F}\}$. Similarly, let $g=F_1...F_m\in G$ for composable $F_i$. There is a unique state in $\mathfrak{F}(\mathcal{A})$, denoted ${}_{(\gamma,n)}{|g}$, satisfying 
 \[{}_{(\gamma,n)}{|g}= {}_{(\gamma,n)}{|F_1}\cdot{}_{(\gamma\cdot F_1,n)}{|F_2}\cdot\ldots \cdot{}_{(\gamma\cdot F_1\ldots F_{m-1},n)}{|F_m}.\addtag\label{seccc}\] 
 The \textit{sections of $g$ along $\gamma$} is the sequence $\{{}_{(\gamma,n)}{|g}\}$.
\end{defn}

The product on the right-hand side of (\ref{seccc}) is well-defined. Indeed, since $\mathcal{A}$ is $\omega$-deterministic, $F_1$ maps $\gamma$ to a unique path $\gamma\cdot F_1$. Hence the sections of $F_2$ along $\gamma\cdot F_1$ have the same properties. In other words, there is a unique choice of element ${}_{(\gamma\cdot F_1,n)}{|F_2}\in{}_{(\gamma\cdot F_1)_n}{|F_2}$ and $F_2$ maps $\gamma\cdot F_1$ to a unique path $\gamma\cdot F_1F_2$. These properties hold for each $F_i$ for $i\in\{2,\ldots ,m\}$. Hence the set of (non-initial) states of $\mathfrak{F}(\mathcal{A})$ consists of elements of the form (\ref{seccc}).

The transitions of $\mathfrak{F}(\mathcal{A})$ have the following form. Write $\gamma=\ldots e_{n+1}e_n\ldots e_2e_1$ and $\gamma'=\gamma\cdot g=\ldots e_{n+1}'e_n'\ldots e_2'e_1'$. Then there is an arrow in the Moore diagram of $\mathfrak{F}(\mathcal{A})$ from ${}_{(\gamma,n)}{|g}$ to ${}_{(\gamma,n+1)}{|g}$ labeled by $e_{n+1}|e_{n+1}'$. 

The following is straightforward.

\begin{prop}
  If $\mathcal{A}$ is stationary, then the set of states of $\mathfrak{F}(\mathcal{A})$ is $G$. 
\end{prop}

In other words, if $\mathcal{A}$ is stationary, we say that the inverse semigroup $G$ is \textit{generated by} the automaton $\mathcal{A}$. 

\begin{defn}\label{finstaaut} Let $\mathcal{A}$ be a time-varying automaton over a simple and stationary Bratteli diagram $\mathsf{B}$. The \textit{minimization} of $\mathcal{A}$, denoted $\mathcal{A}'$, is defined by identifying the states in $\mathcal{A}$ that define the same actions. The automaton $\mathcal{A}$ is said to be \textit{finite-state} if the set of states on the automaton $\mathcal{A}'$ is finite.
\end{defn}
The following follows directly from the definition above. 

\begin{prop}\label{minicoin}
  Let $\mathcal{A}$ be an automaton described in Proposition \ref{Automaton}. Then $\mathfrak{F}(\mathcal{A}')=\mathfrak{F}(\mathcal{A})'$, whose set of initial states is $G$.  
\end{prop}

\begin{prop}\label{fiinstat}
 Let $\mathcal{A}$ be a stationary automaton. Then $\mathcal{A}$ is finite-state.
\end{prop}

\begin{proof}
    Let $B,P$ be finite sets with bijections $\varphi_n:B\rightarrow B_n$ and $\psi_n: P\rightarrow C_n$ satisfying $\sigma\circ\varphi_{n+1}=\varphi_n$ and $\sigma\circ\psi_{n+1}=\psi_n$, for all $n\in\mathbb{N}$. This implies that $\sigma$ is a bijection when restricted on $B_n$ and $C_n$. We first note that any boundary point $\gamma$ of an infinite tile must be periodic. Otherwise, since all $\sigma^{i}(\gamma)$, $i\in\mathbb{N}$, are boundary points of infinite tile(s), there must be infinitely many boundary points of infinite tiles, which is a contradiction. 

Now let $(F,\eta,\eta')$ be a finitary state. Then $\eta$ must have the form $e_k\ldots e_l\sigma^i(\gamma_{dm})$ for some boundary point $\gamma$ of an infinite tile with period $m$, and $d$ can be $0$. By assumption, all $\sigma^j(\eta)$ have the same preimage in $P$. Since $P$ is finite and each $\varphi$ is a bijection, there must be finitely many different possibilities $e_k\ldots e_l$ following different $\gamma_{dm}$. The same argument holds for $(\eta,\eta',F)$ being a connector. Hence the automaton $\mathcal{A}$ is finite-state.
\end{proof}

Recall that a \textit{path} in the Moore diagram of an automaton $\mathcal{A}$ is a sequence of edges that join a sequence of states (vertices of the Moore diagram). A \textit{nontrivial infinite path} of $\mathcal{A}$ is an infinite sequence of edges on which all states are nontrivial. If the minimization automaton $\mathcal{A}'$ is finite, then we call the corresponding finite path in the Moore diagram of $\mathcal{A}'$ obtained from identifying identical states in a nontrivial infinite path of $\mathcal{A}$ a \textit{directed cycle}. 

 \begin{prop}
\label{period} Let $\mathsf{B}$ be a simple, stationary Bratteli diagram whose sets of vertices and edges are finite on each level. Suppose the number of edges is at least $2$ on each level. Let $G$ be an inverse semigroup defined by a stationary automaton $\mathcal{A}$. Suppose $G$ has finitely many boundary points on infinite tiles. Then each boundary point of an infinite tile is a periodic sequence of edges on $\mathsf{B}$. Any two directed cycles in the minimization automaton $\mathcal{A}'$ are disjoint, meaning that they consist of different states and the states in different cycles are not connected by paths.   
\end{prop}

\begin{proof}
 By Proposition \ref{fiinstat}, $\mathcal{A}$ is finite-state. Let $(F,e_1,e_1')\in\mathcal{A}$ be an initial state. Since $G$ is stationary is generated by a stationary automaton, the transformation determined by $(F,e_2e_1,e_2'e_1')$ is also in $G$. It follows that $(H,e_2,e_2')$ is also an initial state of $\mathcal{A}$. The statement is true in general, i.e., if $(F,e_n\ldots e_1,e_n'\ldots e_1')$ is an element of the section of  $(F,e_1,e_1')$, then each $(H_i,e_i,e_i')$ is an initial state. 
 
 Let $F\in\mathcal{S}$. Suppose $F$ is directed along $\gamma\in\Omega(\mathsf{B})$. Then by definition, $\gamma$ is a boundary point of an infinite tile. Since the automaton $\mathcal{A}$ is finite-state, there exist smallest $N,M$, $\gamma_1\in\Omega_{N}(\mathsf{B})$, $\gamma_2\in\Omega_{N+M}(\mathsf{B})$ and a state $(F,\gamma_1,\gamma_1')$ on the directed path of $(F,v_1,v_2)$ in the Moore diagram of $\mathcal{A}$, such that $F$ acting on Dom($F$)$\gamma_1$ is the same as $F$ acting on Dom($F$)$\gamma_2$. Note that Dom($F$)$\gamma_1$ is homeomorphic to Dom($F$)$\gamma_2$ since $\mathsf{B}$ is simple and stationary. This is true for every directed initial state of $\mathcal{A}$. Hence the Moore diagram of the minimization $\mathcal{A}'$ is a finite union of finite cycles and finitely many finitary paths possibly with edges connecting them. It follows that all boundary points of infinite tiles are periodic. 

Now suppose there are directed cycles in $\mathcal{A}'$ that are not disjoint. We distinguish three cases.
\begin{enumerate}
    \item Two directed cycles are connected by a path. Then the number of boundary points of infinite tiles is infinite since we can trace arbitrary long words on the first cycle and then go to the second cycle, contradicting the fact that there are finitely many boundary points. Hence this case cannot happen. 
    \item There exist $\gamma_1\neq\gamma_2$ starting and ending at the same vertex of $\mathsf{B}$ such that there exists a state $(F,\eta_1,\eta_2)$ that is in the section of some $F\in\mathcal{S}$ such that $(F,\gamma_1\eta_1,\gamma_1'\eta_2)$, $(F,\gamma_2\eta_1,\gamma_2'\eta_2)$ and $(F,\eta_1,\eta_2)$ define the same transformation on $\Omega(\mathsf{B})$. Then we can define infinite paths $\gamma\in\Omega(\mathsf{B})$ by concatenating $\gamma_1,\gamma_2$ in any arbitrary order for any arbitrary number of appearances. Then there are infinitely many such $\gamma$ and each of them is an boundary point of an infinite tile since the sections of $F$ along them are all nontrivial. This is a contradiction. Hence this case cannot happen. 
    \item There exist two cycles that share at least one state in the middle (a non-initial state). Then by a similar reasoning as Case 2, there are infinitely many boundary points of infinite tiles. Hence this case cannot happen either.  
\end{enumerate}
The above are all the cases, so the proof is done. 
\end{proof}
Based on the proposition above and the fact that $\mathcal{A}$ is $\omega$-deterministic, let us make the following technical assumptions.
\begin{Assumptions}[leftmargin=5\parindent]
  %  \item  Each $F\in\mathcal{S}$ is only defined at a single vertex $v_1\in V$, so that it is identified with the set $\{(F,v_1,-)\}$;
    \item\label{A2} For each nontrivial finitary state in $\mathcal{A}$ there is a \textbf{unique} finite path in $\Omega^*(\mathsf{B})$ along which the sections of this state is nontrivial; 
    \item\label{A3} Two different nontrivial finitary states are nontrivial along different finite paths. These two paths may have common beginnings but at least the last edges are different, or they have different lengths. 
\end{Assumptions} 
These are valid since if there is a finitary state that has two outgoing edges to different nontrivial states, then we can replace this state with two different states that are deterministic along two different finite paths. By the above assumptions and Proposition \ref{period}, each finitary transformation is identified with a singleton consisting of a finitary state, while the sections of each directed transformation are sets consisting of a unique directed state and several finitary states. 

Let us describe the sets $B,P$. 

By Proposition \ref{fiinstat}, $\mathcal{A}$ is finite-state. Let $(F,e_1,e_1')\in\mathcal{A}$ be a directed state. By Proposition \ref{period}, there exists a unique periodic path $\gamma\in\Omega(\mathsf{B})$ and a smallest $m\in\mathbb{N}$ such that $(F,\gamma_{m+1},\gamma_{m+1}')$ and $(F,e_1,e_1')$ define the same transformation on $\Omega(\mathsf{B})$, where $\gamma_m'=(\gamma\cdot F)_m$. This implies that $\gamma$ has period $m$. Write $\gamma=(e_m\ldots e_1)^{\omega}$. It follows that $\gamma'$ is also a periodic path of period $m$ since $(F^{-1},e_1',e_1)$ has the same properties along $\gamma'$ as $(F,e_1,e_1')$ along $\gamma$. By the first paragraph in the proof of Proposition \ref{period}, we also have that each $(H_i,e_i,e_i')$ and $(H_i,e_{i}\ldots e_m\ldots e_i,e_{i}'\ldots e_m'\ldots e_i')$ define the same transformations, where $H_1=F$. It follows that a subset of $B_n$ has the form  \[B_{\gamma,n}=\{\gamma_n,(\sigma(\gamma))_n,\ldots ,(\sigma^{m-1}(\gamma))_n \}.\] The elements in the set above will repeat if $n>m$. We identify the repeating elements. It follows that given an boundary point of an infinite tile $\gamma$ with period $m$, there are $m-1$ other boundary points of infinite tile(s) of the form $\sigma^i(\gamma)$, for $i=1,..,m-1$, while each $B_{\gamma,n}$ has cardinality $m$, for $n\in\mathbb{N}$.  By a similar reasoning, all finitary states are initial states, and the boundary points corresponding to finitary states on level $n$ have the form $e_n\ldots e_{l+1}\gamma_l$, where $\gamma$ is a boundary point of an infinite tile and $l$ can be $0$. If $l\geq 1$, then $\sigma(e_n\ldots e_{l+1}\gamma_l)$ is also a boundary point on level $n-1$. Hence the cardinalities of each $B_n$ are the same. 

\textbf{Now we define the set $B$}. Denote by $\Omega^{-}(\mathsf{B})$ the set of all right-infinite paths of $\mathsf{B}$. In other words, $\Omega^{-}(\mathsf{B}):=\{e_1e_2\ldots e_i\ldots: \mathbf{s}(e_i)=\mathbf{r}(e_{i+1})\text{ for }i\in\mathbb{N}\}$. Let $\mathsf{s}:\Omega^{-}(\mathsf{B})\rightarrow \Omega^{-}(\mathsf{B})$ be the shift map, i.e., $\mathsf{s}(e_2e_1e_m\ldots )=e_1e_m\ldots $.  Give $\gamma=(e_m\ldots e_1)^{\omega}\in\Omega(\mathsf{B})$, let $\gamma^{-}:=(e_m\ldots e_1)^{-\omega}$ the right-infinite path determined by $\gamma$.

Let $B$ be the set of all right infinite paths in $\mathcal{A}'$ that end in nontrivial states in $\mathcal{A}'$, i.e., they are of the form $e_k\ldots e_l\mathsf{s}^j(\gamma^{-})$ for $\gamma$ being any boundary point of an infinite tile, where $ j\in\{0,\ldots ,m-1\}$ and $e_k\ldots e_j$ can possibly be empty. Since $\mathcal{A}'$ is finite, $B$ is finite. For a finite path $\eta$, denote by $L(\eta)$ its length. Given $\gamma\in\Omega(\mathsf{B})$ a boundary point of an infinite tile, suppose its period is $m$. Define map $\varphi_n:B\rightarrow B_n$ by the following rules:  
\[\varphi_n(\gamma^{-})=\sigma^{|n-dm|}(\gamma_{dm}), \addtag\label{dirdir}\]
where $d$ is such that $L(\sigma^{|n-dm|}(\gamma_{km}))=n$; and
 \begin{align}\label{finfin}
 \varphi_n(e_k\ldots e_l\mathsf{s}^j(\gamma^{-}))=
    \begin{cases}
    e_k\ldots e_{l+c}, \text{ if $n<L(e_k\ldots e_l)$ and $L(e_k\ldots e_{l+c})=n$,}  \\
    e_k\ldots e_l\sigma^i(\mathsf{s}^j(\gamma)_{dm}), \textup{if $n\geq L(e_k\ldots e_l)$},\\
    \end{cases}
\end{align}
where $i,d\geq 0$ are such that $L(e_k\ldots e_l\sigma^i(\mathsf{s}^j(\gamma)_{dm}))=n$. 

Then the map $\varphi_n$ is a bijection. Indeed, for any $\eta\in B_n$, we can append to its right an (eventually) periodic sequence of edges that are on a directed cycle of $\mathcal{A}'$, and thus $\varphi_n$ is surjective. Let $x_1,x_2\in B$. If they are mapped to the same element on the right-hand side of (\ref{dirdir}), (\ref{finfin}), then there is only a unique way to continue the path to the right to obtain an element of $B$ (here we used Assumptions \hyperref[A2]{1} and \hyperref[A3]{2} in the last subsection), which implies $x_1=x_2$. Hence $\varphi_n$ is injective. 

\textbf{Now we define the set $P$}. For $C_n$, we note that given a connector $(\gamma_1,\gamma_2,F)\in P_n$, $F$ acts as a prefix switch $S_{\gamma_2,\gamma_1}$ (see (\ref{switch})) on all paths starting with $\gamma_1$. Hence $(F,e\gamma_1,e\gamma_2)$ is a trivial state for all allowed continuations $e$. Since $G$ is stationary, it also follows that $(\sigma^i(\gamma_1),\sigma^i(\gamma_2),H_i)$ are connectors on levels $n-i$ for $i=0,\ldots ,n-1$, and each $(H_i,e_i\sigma^i(\gamma_1),e_i\sigma^i(\gamma_2))$ is a trivial state for all allowed continuations $e_i$. It follows that we can define $P$ to be the set of right infinite paths on $\mathcal{A}'$ ending in the trivial state. The bijection $\psi_n:P\rightarrow C_n$ is defined the same way as (\ref{finfin}). The relations $\sigma\circ\varphi_{n+1}=\varphi_n$ and $\sigma\circ\psi_{n+1}=\psi_n$ are easily verified. 

The construction above is summarized in the following proposition. 

\begin{prop} Let $B$ be the set of all right infinite paths in $\mathcal{A}'$ ending in nontrivial states, and $P$ be the set of right infinite paths ending in the trivial state of $\mathcal{A}'$. Then $B,P$ are finite. Let $\varphi_n:B\rightarrow B_n$ be the map defined by the rules (\ref{dirdir}), (\ref{finfin}), and let $\psi_n:P\rightarrow C_n$ be the map defined by by the rule (\ref{finfin}). Then for all $n\in\mathbb{N}$, $\sigma\circ\varphi_{n+1}=\varphi_n$ and $\sigma\circ\psi_{n+1}=\psi_n$. 
\end{prop}

\subsection{Contractions}
\begin{defn}\label{contra}
    Let $G$ be a inverse semigroup of bounded type defined by a stationary finite-state automaton $\mathcal{A}$. The inverse semigroup $G$ is called \textit{contracting} if there exists a finite subset $\mathcal{N}\subset G$ such that for each $g\in G$, there exists $k(g)\in \mathbb{N}$ such that ${}_{(\gamma,n)}{|g}\in \mathcal{N}$ for all paths $\gamma\in\Omega(\mathsf{B})$ on which $g$ is defined and for all $n\geq k(g)$. The minimal set $\mathcal{N}$ with this property is called the \textit{nucleus} of $G$. 
\end{defn}
The nucleus $\mathcal{N}$ is described in the following proposition. 

\begin{prop}\label{bdtimpliescon}
   If $G$ is an inverse semigroup of bounded type defined by a stationary finite-state automaton $\mathcal{A}$, then $G$ is contracting. The nucleus $\mathcal{N}$ of $G$ is a finite-state automaton consisting of elements $g\in G$ that are elements in sections of states of $\mathcal{N}$ that are in the directed cycles of the Moore diagram of its minimization $\mathcal{N}'$. 
\end{prop}
\begin{proof}
 Since the automaton $\mathcal{A}$ is finite-state, we can assume without loss of generality that all directed initial states of $\mathcal{A}$ belong to directed cycles of $\mathcal{A}'$. In other words, for every directed $F\in\mathcal{S}$, there exists $\gamma\in\Omega(\mathsf{B})$ such that ${}_{(\gamma,n)}{|F}=F$, for some $n\in\mathbb{N}$. Let $n_1$ be a common multiple of the lengths of all non-trivial cycles in $\mathcal{A}'$, $\mathcal{S}_1$ be the set of directed elements in $\mathcal{S}$, and $\mathcal{S}_0$ be the set of finitary elements in $\mathcal{S}$. We can make $n_1$ sufficiently big (bigger than the lengths of all finitary paths in $\mathcal{A}'$) so that all allowed sections of elements in $\mathcal{S}$ are either finitary or belong to $\mathcal{S}_1$. Let $\mathcal{N}_1'$ be the set of all $g$ that are non-identity products of composable elements in $\mathcal{S}$ such that there exists a unique $\gamma(g)\in\Omega(\mathsf{B})$ such that ${}_{(\gamma(g),n_1)}{|g}=g$ and for all other words $\eta\in\Omega_{n_1}(\mathsf{B})$ on which $g$ are defined the sections ${}_{\eta}{|g}$ consist of only finitary states. Notice that since $G$ is a group of bounded type, $\gamma(g)$ must be unique. The reason is similar to Case 2 in the proof Proposition \ref{period}. Indeed, suppose otherwise. Let $\gamma_1\neq \gamma_2\in\Omega(\mathsf{B})$ whose $n_1$-th truncations start and end at the same vertex such that ${}_{(\gamma_1,n_1)}{|g}={}_{(\gamma_2,n_1)}{|g}=g$. Define infinite paths $\gamma$ by concatenating $\gamma_1,\gamma_2$ in any arbitrary order for any arbitrary number of appearances. Then all sections of $g$ along $\gamma$ will be nontrivial. Thus all these $\gamma$ are boundary points, which is a contradiction to the assumption that there are finitely many boundary points. \textbf{Claim:} the set $\mathcal{N}_1'$ is finite. This is because each of its element $g$ is uniquely defined by how elements in $\Omega_{n_1}(\mathsf{B})$ (a finite set) can be possibly moved by $g$, and by its sections on elements in $\Omega_{n_1}(\mathsf{B})$. By the assumption, the elements in sections of $g$ are either itself or finitary on $\Omega_{n_1}(\mathsf{B})$. Since the number of finitary elements is finite, the set $\mathcal{N}_1'$ must be finite. 

Let $g=F_1\ldots F_k$ for some composable $F_1,\ldots ,F_n\in\mathcal{S}$. Denote by $l_1(g)$ the number of $F_i\in\mathcal{S}_1$ for $F_i$ in the product of $g$. \textbf{Claim:} for each such $g$, there exists a number $k$ such that for every $\gamma\in\Omega_{n_1k}{(\mathsf{B}})$ on which $g$ is defined, the elements in the sections of $g$ on $\gamma$ are either in $\mathcal{N}_1'$ or finitary. This is done by induction. 

If $l_1(g)=1$, then $g=H_1FH_2$ where $H_1,H_2$ are products of elements in $\mathcal{S}_0$ and $F\in\mathcal{S}_1$. Hence $H_1,H_2$ are finitary. Then there exists a $k$ such that the sections ${}_{(\gamma,n_1k)}{|H_i}$ are trivial for all $\gamma\in\Omega(\mathsf{B})$ on which $H_i$ are defined. Then ${}_{(\gamma,n_1k)}{|H_1FH_2}={}_{(\gamma\cdot H_1,n_1k)}{|F}$. Since $F\in\mathcal{S}$, the above statement is true. Hence the statement is true for $l_1(g)=1$. Now suppose that the claim is true for all $g$ such that $l_1(g)<m$. For every $\eta\in\Omega(\mathsf{B})$ on which $F_i$ is defined, the element ${}_{(\eta,n_1k_i)}{|F_i}$ is either equal to $F_i$ or is finitary. Consequently, either ${}_{(\gamma,n_1k)}{|g}=g$ for a unique $\gamma$ and ${}_{(\eta,n_1k)}{|g}$ is finitary for $\eta\neq\gamma\in\Omega(\mathsf{B})$, where $k$ is a common multiple of each $k_i$, or $l_1({}_{(\eta,n_1k)}{|g})<l_1(g)$ for every $\eta\in\Omega_{n_1k}$. Here $l_1({}_{(\eta,n_1k)}{|g})$ is counting the directed elements on the right-hand side of Formula (\ref{seccc}). In the first case, we have $g\in\mathcal{N}_1'$, and in the second case, the claim follows from the induction hypothesis.
 \end{proof}

The following shows that the boundary points of the infinite tiles are the same as the boundary points determined by the automaton $\mathcal{N}$. 
 \begin{prop}\label{agree}
   Let $G$ be a stationary inverse semigroup of bounded type. Then the set of boundary points determined by $\mathcal{A}$ is equal to the set of boundary points determined by the automaton $\mathcal{N}$.   
 \end{prop}
 \begin{proof}
   Assume, without loss of generality, that all initial directed states of $\mathcal{A}$ are in directed cycles of $\mathcal{A}'$. By Proposition \ref{bdtimpliescon}, $\mathcal{N}$ is finite-state whose initial states are in the directed cycles of $\mathcal{N}'$, while the directed cycles are disjoint. 

   We show that the set of infinite paths on which all directed cycles of $\mathcal{A}'$ are directed along is the same as that of $\mathcal{N}'$. Denote these two sets by $\mathsf{D}$ and $\mathsf{D}_{\mathcal{N}}$ respectively. Note that $\mathsf{D}$ is exactly the set of boundary points of infinite tiles by definition. Let $\gamma$ be a boundary point of an infinite tile. Then there exists $F\in\mathcal{S}$ such that $F$ is directed along $\gamma\in\mathsf{D}$. Hence $F\in\mathcal{N}$ and $\gamma\in\mathsf{D}_{\mathcal{N}}$. It follows that $\mathsf{D} \subset\mathsf{D}_{\mathcal{N}}$. On the other hand, let $\gamma\in\mathsf{D}_{\mathcal{N}}$. Write $\gamma=(e_n\ldots e_1)^{\omega}$ for each $e_i\in E$. Let $g$ be an initial state on the directed cycle along $\gamma$, i.e., ${}_{(\gamma,n)}{|g}=g$. Write $g=F_1\ldots F_m$, for $F_1,\ldots ,F_m\in\mathcal{S}$. Let $H'_1,H'_2,\ldots ,H'_n$ be the beginning of the edge labeled by $e_1,\ldots ,e_n$ in the Moore diagram of $\mathcal{N}$. Note that $H_1'=g$ and all $H_j'$ for $j\neq 1$ are obtained by taking sections of $g$ along $\gamma$. Write $H_j'=H_{j1}\ldots H_{j,m_j}$ where $H_{jk}$ are sections of some $F\in\mathcal{S}$. Let $\sigma:\Omega(\mathsf{B})\rightarrow\Omega(\mathsf{B})$ be the shift map. We have $H_2'={}_{(\gamma,1)}{|H'_1}={}_{(\gamma,1)}{|g}= {}_{(\gamma,1)}{|F_1}\cdot{}_{(\gamma\cdot F_1,1)}{|F_2}\cdot\ldots \cdot{}_{(\gamma\cdot F_1\ldots F_{m-1},1)}{|F_m}$, and $H'_{j+1}={}_{(\sigma^{j-1}(\gamma),1)}{|H'_j}={}_{e_j}{|H_{j1}}\cdot{}_{(\sigma^{j-1}(\gamma)\cdot H_{j1},1)}{|H_{j2}}\cdot\ldots \cdot{}_{(\sigma^{j-1}(\gamma)\cdot H_{j1}\ldots H_{j,m_j-1},1)}{|H_{j,m_j}}$, for $j=2,\ldots ,n-1$, and ${}_{(\sigma^{n-1}(\gamma),1)}{|H'_n}=H_1'=g$. Since each $H_j'\in\mathcal{A}_{\mathcal{N}}$, all $m_j$ must be the same; otherwise, there will be some $H_{j,k}$ whose sections are trivial, and this $H_{j,k}$ will not be in the product of $H'_j$. Hence we produced a directed cycle $H_{j+1,1}={}_{(\sigma^{j-1}(\gamma),1)}{|H_{j1}}$ whose initial state is $H_{11}=F_1\in\mathcal{S}$. It follows that this cycle is in $\mathcal{A}'$ and thus $\gamma\in\mathsf{D}\implies \mathsf{D}_{\mathcal{N}}\subset\mathsf{D}$.  
\end{proof}
 
\subsection{Polynomial growth of orbital graphs}
Recall that given a finite graph $\Gamma$, the \textit{diameter} of $\Gamma$, denoted Diam$(\Gamma)$, is defined to be $\max_{u,v\in V(\Gamma)}d(u,v)$. In other words, the diameter of $\Gamma$ is the length of the longest geodesic of $\Gamma$. Let $v\in \Gamma$ be a vertex. Denote by $B_v(R)$ the ball of radius $R$ centered at $v$, and by $\#B_v(R)$ the number of vertices of $\Gamma$ enclosed by $B_v(R)$. The graph $\Gamma$ is said to have \textit{polynomial growth} if the growth of balls of radius $R$ is polynomial, i.e., there exist constants $C,d\geq 1$ such that \[C^{-1}R^d\leq \#B_v(R)\leq CR^d.\addtag\label{polygrow}\] In this subsection, we show that the growth of orbital graphs is bounded above by polynomials. All tiles are assumed to be connected in this subsection.

\begin{prop}
    Let $G$ be a stationary finite-state inverse semigroup of bounded type acting on $\Omega(\mathsf{B})$. Let $\gamma\in\Omega(\mathsf{B})$ be any point. Then the growth of the orbital graph $\Gamma_{\gamma}$ is bounded above by a polynomial of the form $CR^d$. 
\end{prop}
\begin{proof}
    It is enough to prove the statement for the infinite tile $\mathcal{T}_{\gamma}$ since $\Gamma_{\gamma}$ either coincides with $\mathcal{T}_{\gamma}$ or is obtained by connecting a bounded number of infinite tiles at their boundary points. For $g\in G$, denote by $l(g)$ the word length of $g$, i.e., the minimal length of the representation as a product of elements in $\mathcal{S}$. We need the following lemma.
    \begin{lemma}
        Suppose $G$ is contracting. Let  $g\in G$ and $\eta\in\Omega_n(\mathsf{B})$ on which $g$ is defined, write 
         \[{}_{(\zeta,n)}{|g}= {}_{(\zeta,n)}{|F_1}\cdot{}_{(\zeta\cdot F_1,n)}{|F_2}\cdot\ldots \cdot{}_{(\zeta\cdot F_1\ldots F_{m-1},n)}{|F_m}\] as in Formula (\ref{seccc}), for all $\zeta\in\Omega(\mathsf{B})$ having $\eta$ as a finite truncation. Define $l({}_{\eta}{|g})$ to be the number of non-identity transformations on the right-hand side in the formula above. 
        
        Then there exists a number $M>0$ and a positive integer $n>0$ such that the inequality
        \[l({}_{(\zeta,n)}{|g})\leq\frac{l(g)}{2}+M\] holds.
    \end{lemma}
     We note that for a different $\zeta'$ having $\eta$ as a finite truncation, $l({}_{(\zeta',n)}{|g})$ might have a different value. 
    \begin{proof}

       % is because the vertices of $\mathcal{T}_{\gamma}$ are points in $\Omega(\mathsf{B})$ while $\eta$ is a finite truncation of some $\zeta\in\mathcal{T}_{\gamma}$, and the automaton $\mathcal{A}$ defining $G$ is eventually deterministic.  

        Let $M$ be the maximal length of the elements in the nucleus $\mathcal{N}$. Then there exists a number $n\in\mathbb{N}$ such that for every $g\in G$ with length not more than $2M$ and every path $\eta\in\Omega_n(\mathsf{B})$ on which $g$ is defined, we have ${}_{\eta}{|g}\subset\mathcal{N}$.  

        Let $g\in G$ be any element. It can be written in the form $g=g_1\ldots g_kg_{k+1}$, where $k=\left\lfloor\frac{l(g)}{2M}\right\rfloor$, $l(g_i)=2M$ for all $1\leq i\leq k$ and $l(g_{k+1})<2M$. Then for every $\eta\in\Omega_n(\mathsf{B})$ on which $g$ is defined, and for every $\zeta\in\Omega(\mathsf{B})$ having $\eta$ as a finite truncation, the element ${}_{(\zeta,n)}{|g}$ can be written in the form $h_1\ldots h_kh_{k+1}$, where $h_i\in\mathcal{N}$. It follows that 
        \[l({}_{\eta}{|g})\leq(k+1)M=\left(\left\lfloor \frac{l(g)}{2M}\right\rfloor +1 \right)M\leq\left( \frac{l(g)}{2M}+1\right)M=\frac{l(g)}{2}+M.\]
    \end{proof}
    
    Let $M,n$ be as in the lemma above. Let $\sigma:\Omega(\mathsf{B})\rightarrow\Omega(\mathsf{B})$ be the shift map. Let $v\in\mathcal{T}_{\gamma}$ be a generic point. The cardinality of the ball $B(v,R)=\{v\cdot g:l(g)\leq R\}$, denoted $\#B(v,R)$, is not greater than
    \[|E|^n\cdot\#B(\sigma^n(v),\frac{R}{2}+M),\]
    since the $n$-th iteration of the shift map $\sigma^n$ maps $B(v,r)$ to $B(\sigma^n(v),\frac{R}{2}+M)$, and every point in $\Omega(\mathsf{B})$ has at most $|E|^n$ preimages under $\sigma^n$. 
    
    Let $k=\left\lfloor \frac{\log R}{\log 2}\right\rfloor+1$. Then $(\frac{1}{2})^kR<1$ and 
    \[\#B(v,R)\leq |E|^{nk}\#B(\sigma^n(v),R')\] where 
    \[R'=(\frac{1}{2})^kR+(\frac{1}{2})^{k-1}M+(\frac{1}{2})^{k-2}M+\ldots +\frac{1}{2}M+M<1+\dfrac{M}{1-1/2}=1+2M.\]

    Let $w\in\Omega(\mathsf{B})$ be any point. Then $\#B(w,R')<|\mathcal{S}|^{R'}$. Hence
    \[\#B(v,R)<|\mathcal{S}|^{R'}\cdot|E|^{nk}<|\mathcal{S}|^{R'}\cdot |E|^{n\left(\frac{\log R}{\log 2}+1\right)}\]
    \[=K_1\cdot \exp\left(n\dfrac{\log R\log|E|}{\log 2}+n\log|E|  \right)=K_2\cdot R^{\frac{n\log|E|}{\log 2}},\]
where $K_1=|\mathcal{S}|^{R'}$ and $K_2=K_1\cdot |E|^n$. Let $C=K_2$ and $d=\frac{n\log|E|}{\log 2}$. The proof is complete. 
\end{proof}

\subsection{Repetitivity of finite tiles}\label{reppft}
The results in this subsection will be used in Section \ref{exmmmmp}. 

\begin{prop}\label{expdiam}
    Let $\{\mathcal{T}_{i,n}\}$ be the set of $n$-th level tiles, for $n\in\mathbb{N}$, whose vertices end in the same $v_i\in V$. Then \textup{Diam}$(\mathcal{T}_{i,n})$ is asymptotically equivalent to $d'(\gamma_{1,n},\gamma_{2,n})$ where $\gamma_{1,n},\gamma_{2,n}$ are boundary points on $\mathcal{T}_{i,n}$ with $\{d'(\gamma_{1,n},\gamma_{2,n})\}$ growing exponentially. Consequently, the sequence of diameters of a tile inflation has exponential growth.    
\end{prop}
\begin{proof}
    The proof is following the idea of that of \cite[Theorem V.3]{Bon} but in a simplified fashion. First note that the diameter is not smaller than $d'(u,v)$ for any $u,v\in\mathcal{T}_{i,n}$. Hence $\textup{Diam}(\mathcal{T}_{i,n})$ grows at least exponentially. We need to show the other inequality. 

    Let $l$ be a geodesic path that represents the diameter of $\mathcal{T}_{i,n}$. Let $v$ and $u$ be the end vertices of $l$. The vertices $v,u$ may or may not be boundary points of $\mathcal{T}_{i,n}$ and may or may not be boundary points of any isomorphic tiles $\mathcal{T}_{j,n-1}$ contained in $\mathcal{T}_{i,n}$. The path $l$ can be represented in the form $[l_1]s_1l_2s_2\ldots s_m[l_{m+1}]$, where each $l_k$ lies in some isomorphic copy of $\mathcal{T}_{j,n-1}$, for various $j\in\{1,2,\ldots ,|V|\}$, and $s_k\in S$ is the label of a boundary edge of $\mathcal{T}_{j,n-1}$. The brackets mean that the beginning and ending paths may just be the boundary edges. Assume without loss of generality that $l_1$ and $l_{m+1}$ has length greater than $1$, and that $v$ and $u$ are not boundary points of $\mathcal{T}_{i,n}$. For each $l_k$, denote by $p_k$ its beginning vertex and $q_k$ its ending vertex. Note $p_1=v$ and $q_{m+1}=v$. Then 
    \[d'(v,u)=d'(v,q_1)+\sum\limits_{k=2}^{m}d'(p_k,q_k)+d'(p_{m+1},u)+m.\] The $m$ at the end counts the number of boundary edges $s_k$. Since the distances between boundary points grow exponentially, we can choose $a>1$ such that $d'(p_k,q_k)\leq a^{n-1}$. Let $C$ be the largest number of $(n-1)$-th tiles that are contained in an $n$-th level tile. Then we have $m\leq C$ and \[d'(v,u)\leq d'(v,q_1)+d'(p_{m+1},u)+C(a^{n-1}+1).\addtag\label{neq1}\]
    We can do similar estimates for $d'(v,q_1)$ and $d'(p_{m+1},u)$. The geodesic connecting $v$ and $q_1$ induces geodesics connecting various isomorphic copies of tiles on level $n-2$. This path is represented by $l'=l'_1s'_1l'_2s'_2\ldots s'_{m'}l'_{m'+1}$. The beginning point of $l_1'$ is $v$, the ending point of $l'_{m'+1}$ is $q_1$, and all other boundary points are denoted $p'_{k'},q'_{k'}$ for $k'\in\{2,\ldots ,m\}$. It follows that \[d'(v,q_1)\leq d'(v,q_1')+C(a^{n-2}+1).\] By a similar reasoning, we have \[d'(p_{m+1},u)\leq d'(p'_{m+1},u)+C(a^{n-2}+1).\] Continuing this process, we have \[d'(v,q_1)\leq d'(v,q_1^{(n+1)})+C\sum\limits_{\nu=1}^{n-2}(a^{\nu-2}+1)\leq C_1a^n, \addtag\label{neq2}\] and \[d'(p_{m+1},u)\leq d'(p_{m^{(n+1)}},u)+C\sum\limits_{\nu=1}^{n-2}(a^{\nu-2}+1)\leq C_1a^n,\addtag\label{neq3}\] where the upper script $(n+1)$ is the number of primes labeling the boundary points. Combining inequalities (\ref{neq1}) (\ref{neq2}) (\ref{neq3}), we have that \textup{Diam}$(\mathcal{T}_{i,n})$ is bounded above by an exponential function. This completes the proof. 
\end{proof}
\begin{prop}\label{linlin}
    Let $G$ be an inverse semigroup of bounded type acting minimally on $\Omega(\mathsf{B})$. Suppose the distances of all boundary points on finite tiles grow exponentially. Then all finite tiles are linearly repetitive in all infinite tiles. 
\end{prop}
\begin{proof}
    Since $\mathsf{B}$ is simple, for any $u,v\in V$, for $n\in \mathbb{N}$, there exists $m>n$ such that $v\in V_n=V$ is connected to $u\in V_m=V$. Since the action is minimal, for the same $n\in \mathbb{N}$, and any $i,j\in\mathbb{N}$, there exists $m'>n$ such that the tile $\mathcal{T}_{i,n}$ is contained as an isomorphic copy in $\mathcal{T}_{j,m}$. For each $n$, let $m$ be the smallest number such that every tile on level $m$ contains every tile on level $n$. In other words, all tiles on level $m$ are covered by all tiles on level $n$ that are connected to each other at their boundary points. Since $\mathsf{B}$ is stationary, $m-n$ is a constant. Denote $L=m-n$. Then all tiles on level $m+L$ will contain all tiles on level $m$. 

    Fix $\mathcal{T}_{j,m+L}$. It contains at most $|E|^L$ copies of tiles on level $m$ and at most $|E|^{2L}$ copies of tiles on level $n$. Since the distances of all boundary points grow exponentially, by Proposition \ref{expdiam}, the diameters of all tiles on level $n$ are asymptotically equivalent to $a^n$ for some $a>1$. Hence the gap between consecutive appearance of $\mathcal{T}_{i,n}$ in $\mathcal{T}_{j,m+L}$ has diameter not more than $(2|E|^L-2)a^n\sim (2|E|^L-2)\cdot\textup{Diam}(\mathcal{T}_{i,n})$. This is because each $\mathcal{T}_{j,m}$ contains at least one copy of $\mathcal{T}_{j,n}$, and if $\mathcal{T}_{j_1,m}$, $\mathcal{T}_{j_1,m}$ are adjacent to each other on $\mathcal{T}_{j,m+L}$, then the gap between $\mathcal{T}_{i,n}$ in $\mathcal{T}_{j_1,m}$ and $\mathcal{T}_{i,n}$ in $\mathcal{T}_{j_2,m}$ has diameter not more than the sum of the diameters of all tiles on level $n$ contained in them. Since infinite tiles are obtained from finite tile inflations, the statement is true for any $\mathcal{T}_{i,n}$ in any infinite tiles. 
\end{proof}

For groups acting by tree automorphisms (e.g. contracting self-similar groups), we have the following. 
\begin{prop}\label{strrepp}
Let $G$ be a group of bounded type acting minimally by tree automorphisms on the boundary $\partial T_d$ (i.e. $\Omega(\mathsf{B})$ where the vertex set $V$ of $\mathsf{B}$ is a singleton) of a regular rooted tree $T_d$ of degree $d$. Then each finite tile is strongly repetitive in all the infinite tiles on $\partial T_d$. 
\end{prop}
\begin{proof}
There is only one tile on each level of $T_d$. Denote by $\mathcal{T}_n$ the $n$-th level tile. Since $G$ is a group of bounded type, the number of boundary points on finite tiles is bounded. Suppose, without loss of generality, that all finite tiles have the same number of boundary points, denoted $\gamma_{1,n},..,\gamma_{k,n}$ on $\mathcal{T}_n$. Then $\mathcal{T}_{n+1}$ is obtained by taking $|E|$ copies of $\mathcal{T}_n$, adding a letter $e\in E$ to each vertex on $\mathcal{T}_n$, denoted $e\mathcal{T}_n$, and connect all $e\mathcal{T}_n$ at the points $e\gamma_{i,n}$ by the corresponding generators $s\in S$ such that ${}_{e\gamma_{i,n}}{|s}=Id$ and $e\gamma_{i,n}\cdot s=e_1\gamma_{j,n}$ for $e_1\in E$. Since the action is minimal, it is level-transitive. The tile $\mathcal{T}_{n+1}$ will have all points on level $(n+1)$ of $T_d$ as its vertices. Hence the gap between two consecutive isomorphic copies of $\mathcal{T}_{n}$ in $\mathcal{T}_{n+1}$ has a bounded size (not more than $k|S|$). Now let $\mathcal{T}_{\eta}$ be an infinite tile. Then it is obtained by taking tile inflations of $\mathcal{T}_{\eta_n}=\mathcal{T}_n$ where $\eta_n$ is the $n$-th truncation of $\eta$. It follows that the gap between consecutive isomorphic copies of $\mathcal{T}_n$ in $\mathcal{T}_{\eta}$ has size not more than $k|S|$. Hence $\mathcal{T}_n$ is strongly repetitive in any infinite tiles for all $n\in\mathbb{N}$. 
\end{proof}

\begin{rmk}
    Here we showed a weaker version of repetitivity. Namely, we do not require that any finite connected subgraphs are (linearly/strongly) repetitive as in the original definitions in Subsection \ref{repppttt}. We only require that this holds for all finite tiles. 
\end{rmk}

\subsection{Incompressible elements}
Let $\mathsf{B}=(V,E,\mathbf{s}_n,\mathbf{r}_n)$ be a simple Bratteli diagram. Let $G$ be an inverse semigroup of bounded type generated by $\mathcal{S}$ as before.
%Let $S$ be the generating set of $G$. 

\begin{defn}
     \label{trajectory} Let $v$ be a vertex in an orbital/tile graph. Let $g=F_1\ldots F_m$ for composable $F_i\in\mathcal{S}$. The \textit{trajectory} or \textit{walk} of $g$ starting at $v$ is the path \[v\cdot F_1,v\cdot F_1F_2,\ldots ,v\cdot F_1F_2\ldots F_m.\] The points $v$ and $v\cdot F_1F_2\ldots F_m$ are called \textit{initial} and \textit{final} vertices of the trajectory of $g$. 
\end{defn}

\begin{defn}\label{return} Let $\xi\in\Omega(\mathsf{B})$ be a boundary point of an infinite tile of $G$ on $\mathcal{T}_{\xi}$ (as a subgraph of the orbital graph $\Gamma_{\xi}$). Let $g=F_1\ldots F_m$, for composable $F_1,\ldots ,F_m\in \mathcal{S}$, be an element in $G$. We say that $g$ is a \textit{return word} if the final vertex of the trajectory of $g$ starting at $\xi\cdot F_1^{-1}$ is $\xi\cdot F_m$, $F_1$ and $F_m$ are labels of boundary edges at $\xi$ on the tile $\mathcal{T}_{\xi}$, while all edges in between labeled by $F_2,\ldots ,F_{m-1}$ on the trajectory of $g$ belong to $\mathcal{T}_{\xi}$ and are not boundary edges of $\mathcal{T}_{\xi}$. 
\end{defn}
\begin{defn}\label{incom} An element $g=F_1\ldots F_n$ is said to be \textit{incompressible} if it does not contain any return subword.
\end{defn}

\section{Main theorem}\label{mmain}
\begin{thm}\label{main}  Let $\mathsf{B}=(V,E,\mathbf{s},\mathbf{r})$ be a simple and stationary Bratteli diagram. Let $G$ be a group of bounded type of homeomorphisms of $\Omega(\mathsf{B})$ defined by a stationary automaton $\mathcal{A}$. Suppose $G$ acts minimally on $\Omega(\mathsf{B})$. If the set of incompressible elements of $G$ is finite, then the growth function $\gamma_G(R)\preccurlyeq \exp(R^{\alpha})$  for some $\alpha\in (0,1)$.
\end{thm}
\begin{rmk}\label{allcon1} By virtue of Subsection \ref{reppft}, we assume that $\mathsf{B}$ satisfies that for any $v\in V$ there is at least $1$ edge $e\in E$ such that $v\in V_n$ is connected to $v\in V_{n+1}$ by $e$. In other words, each tile $\mathcal{T}_{i,n+1}$ contains at least $1$ copy of $\mathcal{T}_{i,n}$. Otherwise, we may pass to the least common multiple $M$ with respect to all $i\in\{1,\ldots ,|V|\}$ such that $\mathcal{T}_{i,n+M}$ contains at least $1$ copy of $\mathcal{T}_{i,n}$, for all $n\in\mathbb{N}$. 
\end{rmk}

\begin{rmk}\label{allcon2}
Since the automaton $\mathcal{A}$ is stationary, which implies all boundary points on infinite tiles are periodic sequences of paths in $\Omega(\mathsf{B})$, it follows that given $\mathcal{T}_{i,n}$ and a boundary point $\xi_n$ which is a truncation of boundary point $\xi$ on an infinite tile there exists $M$ such that the tile $\mathcal{T}_{i,n+M}$ has a boundary point of the form $\xi_{n+M}$. In other words, there is an isomorphic copy of $\mathcal{T}_{i,n}$ on $\mathcal{T}_{i,n+M}$ that shares a boundary point with $\mathcal{T}_{i,n+M}$. Hence, by combining Remark \ref{allcon1}, we may assume that all tiles $\mathcal{T}_{i,n+1}$ on level $n+1$ contain a copy of $\mathcal{T}_{i,n}$ that shares a boundary point with it.
\end{rmk}

Denote by $S$ the generating set of $G$. Since $\mathsf{B}$ is stationary, the number of different tiles on each $\Omega_n(\mathsf{B})$ are the same. Each tile $\mathcal{T}_{i,n}$ on $\Omega_n(\mathsf{B})$, for $i=1,\ldots ,|V_n|$, corresponds to a vertex in $V_n=V$.

\begin{defn}\label{traverse} Let $\mathcal{T}_{i,n}$ be a tile. Let $w=s_1s_2\ldots s_m\in S^*$ be a word. A \textit{traverse} of $\mathcal{T}_{i,n}$ is the trajectory (Definition \ref{trajectory}) of a subword of $w$, denoted $w'=s_i\ldots s_j$ for $1\leq k<j \leq l$, such that the initial and final vertices of the trajectory of $w'$ are boundary points of $\mathcal{T}_{i,n}$ and all other vertices of the trajectory of $w'$ are inside $\mathcal{T}_{i,n}$ different from the boundary points.  
\end{defn}

We also need the following technical definition. 
\begin{defn}\label{reduced}
    Let $w=s_1\ldots s_m\in S^*$ be a word. We say that $w$ is \textit{reduced} if all $s_i\neq Id$ and either all consecutive $s_i,s_{i+1}$ are labels of boundary edges at different boundary points of infinite tile(s) or at least one of them is finitary. Otherwise, the word $w$ is said to contain \textit{non-reduced subwords}.  
\end{defn}

Fix a word $w=s_1\ldots s_m\in S^*$. For each tile $\mathcal{T}_{i,n}$, denote by $\Theta_{i,n}$ the set of all traverses of $w$ on the tile $\mathcal{T}_{i,n}$, and let $\Theta_n=\bigcup_{i}\Theta_{i,n}$. Define a generating function 
$$F_{w,i}(t)=\sum\limits_{n=0}^{\infty}T_{i,n}t^n,$$ where $T_{i,n}:=\#\Theta_{i,n}$ denotes the number of traverses of $w$ on $\mathcal{T}_{i,n}$. Note that this is a finite sum. Also Let $$F_{w}(t)=\sum\limits_{n=0}^{\infty}T_nt^n,$$ where $T_n=\sum_{i}T_{i,n}$. In other words, $F_{w}(t)=\sum_{i}F_{w,i}$. 

Following the idea of \cite{BNZ}, we have the following. 
\begin{prop} \label{ultesti}
 Suppose that there are $t_0>1$ and a function $\Psi(t)$ such that $F_{w}(t)\leq |w|\cdot\Psi(t)$ for all $w\in S^*$, all $t<t_0$. Set $\beta=\limsup_{n\rightarrow \infty}(\sum_i|\mathcal{T}_{i,n}|)^{1/n}$. Then the growth function $\gamma_G(R)$ of $G$ satisfies $$\gamma_G(R)\preccurlyeq \exp(R^{\alpha})$$ for all $\alpha>\frac{\log \beta}{\log\beta+\log t_0}$. 
\end{prop}

\begin{proof}
Each $\mathcal{T}_{i,n+1}$ is obtained by taking copies of $\mathcal{T}_{i_1,n},\ldots ,\mathcal{T}_{i_k,n}$, for $i_1,\ldots ,i_k\in\{1,\ldots ,|V_n|\}$, adding a new letter labeled by an edge starting at a point in $V_n$ to each vertex, and connecting them at one of their boundary points. Especially, since $\mathsf{B}$ is simple and the action is minimal, we assume that each $\mathcal{T}_{i,n+1}$ contains at least one isomorphic copy of $\mathcal{T}_{i,j}$ for all $j<n+1$ (see the remark after Theorem \ref{main}). Hence some traverses of $\mathcal{T}_{i,n+1}$ will induce traverses of $\mathcal{T}_{i,j}$, by taking restricting of the corresponding word on an isomorphic copy $\mathcal{T}_{i,j}$ on $\mathcal{T}_{i,n+1}$. Let $(B,x,y)\in \mathcal{P}_w$ be as in Definition \ref{stnpor}. Then there exists $m$ such that $\mathcal{T}_{i,m}$ contains an isomorphic copy of $B$ %\sout{Since the orbital graphs have polynomial growth}
and $\sum_i|\mathcal{T}_{i,m}|/\#B$ is uniformly bounded. By linear repetitivity, there also exists $n$ such that we can choose isomorphic copies $\mathcal{T}_{i,n}$ contained in $B$ inside $\mathcal{T}_{i,m}$ with gaps not more than $C\cdot \textup{Diam}(\mathcal{T}_{i,n})$ between consecutive appearances, and $\#B/|\mathcal{T}_{i,n}|$ is bounded. This gives a traverse of $\mathcal{T}_{i,n}$. Hence we defined a map $\chi_i:\mathcal{P}_w\rightarrow \Theta_{i,n}$, mapping $(B,x,y)$ to $\chi_i(B,x,y)$ which is a traverse of $\mathcal{T}_{i,n}$. The bi-rooted graph $(B,x,y)$ is uniquely determined by the traverses and the corresponding isomorphic copies of $\mathcal{T}_{i,n}$ inside $\mathcal{T}_{i,m}$. Since the numbers of boundary points of $\mathcal{T}_{i,n}$ are finite and invariant, and $|\mathcal{T}_{i,m}|/|\mathcal{T}_{i,n}|$ is uniformly bounded, and the number of isomorphic copies of $\mathcal{T}_{i,n}$ in $\mathcal{T}_{i,m}$ is uniformly bounded, it follows that there is a constant $C_i$ such that the map $\chi$ is at most $C_i$-to-one. Hence the summand $(\#B)^{p-1}$ in the definition of $N_p(w)^p$ corresponds the coefficient of $t^n$ in the definition of $F_{w,i}(t)$. The correspondence is given by the map $\chi_i$. 

Let $\beta=\lim\sup\limits_{n\rightarrow \infty}(\sum_i|\mathcal{T}_{i,n}|)^{1/n}$. 
%\sout{Since the orbital graphs have polynomial growth} 
There exists $K_i>0$ such that $|\mathcal{T}_{i,n}|\leq K_i\beta_1^n$ for any $\beta_1>\beta$ and all $n\geq 1$. Consider $F_{w,i}(\beta_1^{p-1})$. Then the map $\chi_i$ associates a summand $(\# B)^{p-1}$ of $N_p(w)^p$ to a summand of the coefficient of $\beta_1^{n(p-1)}$ of $F_{w,i}(\beta_1^{p-1})$. And the ratio 
$$\dfrac{(\#B)^{p-1}}{\beta_1^{n(p-1)}}\leq \Bigl(K_0\dfrac{\#B}{\sum_i|\mathcal{T}_{i,m}|}\Bigr)^{p-1}.$$
Hence it is bounded by a constant not depending on $n$, $w$. Since the map $\chi_i$ is at most $C_i$-to-one, there exists $C_{1,i}(p)$ such that $$N_p(w)^p\leq C_{1,i}(p)F_{w,i}(\beta_1^{p-1}).$$
By the assumption, $$F_{w,i}(\beta_1^{p-1})\leq |w|\Psi_{i}(\beta_1).$$

The above process is done for each $i$. Let $C_{2}(p)=(\max_i C_{1,i}(p))\cdot(\max_i \Psi_{i}(\beta_1))$. Then $$N_{p}(w)^p\leq C_2(p)|w|$$ for all $p$ such that $\beta_1^{p-1}<t_0$, i.e., for all $p<1+\frac{\log t_0}{\log \beta_1}$. 

Since the orbital graphs have polynomial growth, it follows from Proposition \ref{upperbound} that the growth function $\gamma_{G}(R)\preccurlyeq \exp(R^\alpha)$ for all $\alpha>(1+\frac{\log t_0}{\log \beta_1})^{-1}$. 
\end{proof}

It remains to find a function $\Psi(t)=\max_i \Psi_{i}(t)$ as in the proof of the above proposition. This is done in the following proof. 

\begin{proof}[Proof of Theorem \ref{main}.]
Since the set of incompressible elements is finite, there exists $M>0$ such that every reduced word $w\in S^*$ that has length greater than $M$ contains a return.  By minimality, we may assume, without loss of generality, that each $\mathcal{T}_{i,n+1}$ contains all isomorphic copies of all tiles on level $n$. Fix a word $w=s_1\ldots s_m\in S^*$. Let $\tau$ be a traverse of $\mathcal{T}_{i,n+1}$ (viewed as a subword of $w$). Then $\tau$ induces a finite sequence of traverses of tiles on the $n$-th level. Denote this sequence by $\Phi(\tau)=(\tau_1,\tau_2,\ldots ,\tau_k)$. Let $\mathbf{w}=x_1x_2\ldots x_k$ be the sequence of labels of the incoming edges of each induced traverse at each boundary point. Each $x_i$ is a product of elements in $S$, and $x_1$ can be any edge going into the tile $\mathcal{T}_{i,n+1}$. 
\begin{defn}\label{lastmoment}
    Define a map $\phi:\Theta_{i,n+1}\rightarrow\Theta_{i,n}$ by letting $\phi(\tau)$ be the \textbf{last} traverse induced by $\tau$ of an isomorphic copy of $\mathcal{T}_{i,n}$ on which there is at least one boundary point that is a truncation of a boundary point of $\mathcal{T}_{i,n+1}$ (see Remark \ref{allcon2}). Call the map $\phi$ the \textit{last-moment map}, and $\phi(\tau)$ is called the \textit{last induced traverse}.
\end{defn}
We distinguish $2$ cases of $\mathbf{w}$: 1. $\mathbf{w}$ is a reduced word; 2. $\mathbf{w}$ contains non-reduced subwords. 

Suppose first the word $\mathbf{w}$ is reduced. Since each tile $\mathcal{T}_{i,n+1}$ is partly obtained by embedding from  $\mathcal{T}_{i,n}$, there always exists a copy of $\mathcal{T}_{i,n}$ in $\mathcal{T}_{i,n+1}$ such that at least one of its boundary points is embedded in $\mathcal{T}_{i,n+1}$ to be a boundary point of $\mathcal{T}_{i,n+1}$. It follows that we can combine all last-moment maps on each $\Theta_{i,n+1}$ to be a single map $\phi:\Theta_{n+1}\rightarrow\Theta_n$. This map is injective since the trajectory of the induced traverse on $\mathcal{T}_{i,n}$ uniquely determines the traverse of $\mathcal{T}_{i,n+1}$ that induces it. We would like to show $\phi$ is non-surjective. 

Now let $\tau$ be a traverse of any tile on the $(n+1)$-st level. Write \[\Phi(\tau)=(\tau_1,\ldots ,\tau_j,\tau_{j+1},\ldots ,\tau_k)\] where $\tau_j=\phi(\tau)$. Consider the sub-label $\mathbf{w}'=x_jx_{j+1}\ldots x_l$ of $\mathbf{w}$, for $j< l \leq k$. If $\mathbf{w}'$ itself is a return word, then $\tau_j$ cannot be in the image of $\phi$ since $\tau_{l}$ is the last one. Now suppose $\mathbf{w}$ is long enough (with length greater than $M$). Then by the assumption, $\mathbf{w}$ contains a return subword. Let $x_{l_1}\ldots x_{l_2}$ represent this subword, and denote the corresponding sequence of traverse by $(\tau_{l_1},\ldots ,\tau_{l_2})$. Note that the trajectory of the induced traverses might not return to the same vertex, but if we move this subword to the corresponding tile that defines this return, its trajectory will return to the same vertex on that tile. Then the traverse $\tau_{l_1}$ cannot be in the image of $\phi$ for the same reason. It follows that for each $M$ induced traverses, there is at least $1$ that is not in the image of $\phi$. This implies that $$T_{n+1}\leq \dfrac{M-1}{M}T_n.$$ Then the generating function $F_w(t)=\sum\limits_{n=1}^{\infty}T_nt^n$ satisfies
\[F(t)-T_0\leq \frac{M-1}{M}tF(t),\] \[\implies F(t)\leq\dfrac{T_0}{1-\frac{M-1}{M}t}.\addtag\label{redd}\] When $t< \frac{M}{M-1}$, the right hand side is positive. Note that $T_0=|w|$. Hence we can take $\Psi(t)=\frac{1}{1-\frac{M-1}{M}t}$. 

% By Proposition \ref{ultesti}, the growth function \[\gamma_G(R)\preccurlyeq \exp(R^{\alpha})\] for every $\alpha>\dfrac{\log(\beta)}{\log(\beta)+\log(\frac{M}{M-1})}$, where $\beta=\lim\sup\limits_{n\rightarrow \infty}(\sum_i|\mathcal{T}_{i,n}|)^{1/n}$.

Now suppose $\mathbf{w}$ contains non-reduced subwords. Let $n\in\mathbb{N}$ be big enough. Let $r<n$ be an integer satisfying the following properties. For any traverse $\tau$ of $\mathcal{T}_{i,n+1}$ as above, denote by $\Phi(\tau)=(\tau_1,\ldots ,\tau_k)$ the sequence of traverse of tiles on level $r$ of $\mathsf{B}$. Let $\mathbf{w}_r=x_1\ldots x_k$ be the sequence of labels of incoming edges at each boundary point, where each $x_i$ is the section of some $s\in S$ on a boundary point or connecting point $\gamma\in\Omega_{r-1}(\mathsf{B})$, and $\mathbf{w}_r$ is a reduced word. In other words, the word $\mathbf{w}_r$ is an element in the group $G_r=G$, where $G_r$ is defined by tile inflations starting on level $r$ of $\mathsf{B}$. We can choose $r$ such that no consecutive $x_i,x_{i+1}$ are labels of boundary edges (of $G_r$). For example, if $x_i$ is the label of a boundary edge, then $x_{i+1}$ is the label of an edge at a connecting point, which is finitary. Define the last moment map $\phi:\Theta_{i,n+1}\rightarrow\Theta_{i,r}$ in the same way as Definition \ref{lastmoment}. By the assumption, every subword of $\mathbf{w}_r$ of length $M$. By a similar argument as above, for each $M$ induced traverses, there is at least $1$ that is not in the image of $\phi$. This implies that 
\[T_{n+1}\leq\dfrac{M-1}{M}T_r.\]
The generating function $F_w(t)$ satisfies
\[Mt^{-(n+1-r)}(F_w(t)-(T_0+T_1t+\ldots +T_{n-r}t^{n-r}))\leq (M-1)F_w(t)\] \[\implies F_w(t)\leq |w|\dfrac{g(t)}{t^{-(n+1-r)}-\frac{M-1}{M}},\addtag\label{nonredd}\] where $g(t)$ is positive. Notice that $n-r$ is constant. Taking $\Psi(t)=\dfrac{g(t)}{t^{-(n+1-r)}-\frac{M-1}{M}}$, then $F_w(t)\leq|w|\cdot\Psi(t)$ whenever $t<\sqrt[n+1-r]{\frac{M}{M-1}}$. 

Combining (\ref{redd}) and (\ref{nonredd}), and by Proposition \ref{ultesti}, we conclude that the growth function
 \[\gamma_G(R)\preccurlyeq \exp(R^{\alpha})\] for every \[\alpha>\max\left\{\dfrac{\log(\beta)}{\log(\beta)+\log(\frac{M}{M-1})},\dfrac{\log(\beta)}{\log(\beta)+\log\left(\sqrt[n+1-r]{\frac{M}{M-1}}\right)}\right\}\] 
 \[=\dfrac{\log(\beta)}{\log(\beta)+\log\left(\sqrt[n+1-r]{\frac{M}{M-1}}\right)},\]
 where $\beta=\lim\sup\limits_{n\rightarrow \infty}(\sum_i|\mathcal{T}_{i,n}|)^{1/n}$.
The proof is complete.
\end{proof}
\begin{rmk}
    The choice of $\Psi(t)$ is not unique. There might be better ways of choosing such a function. See the example in Subsection \ref{k110} and the examples in \cite{BNZ}. In the case when $\mathbf{w}=x_1\ldots x_k$ contains non-reduced subwords, we cannot directly apply the method in the first case to show non-surjectivity of $\phi$. For example, in the extreme case when all $x_1,..,x_k$ are labels of boundary edges at the same boundary point, it is impossible to find return subwords of $\mathbf{w}=x_1\ldots x_k$.
\end{rmk}

\section{Examples}\label{exmmmmp}
All the examples discussed here use the technique of \textit{fragmentations}. See Subsection \ref{GOA} for the definition. The examples in Subsections \ref{k110}, \ref{fggrp} are \textit{self-similar groups}. Let us recall the definitions from \cite{nekrash05} (in right action notations).  
\begin{defn}
    Let $\mathsf{X}$ be a finite alphabet and $\mathsf{X}^*$ be the set of all finite words over $\mathsf{X}$. The set $\mathsf{X}^*$ can be identified with a regular rooted tree of degree $d=|\mathsf{X}|$. Denote by $\textup{Aut}(\mathsf{X}^*)$ the group of automorphisms of $\mathsf{X}^*$.  A \textit{(faithful) self-similar group} ($G,\mathsf{X}$) is a group $G\leq \textup{Aut}(\mathsf{X}^*)$ together with a faithful action on $\mathsf{X}^*$ such that for every $g\in G$ and $x\in X$  there exist $h\in G$ and
$y\in \mathsf{X}$ such that
\begin{equation*}
vx\cdot g=v\cdot h y
\end{equation*}
for all $v\in \mathsf{X}^*$. 
\end{defn}
Self-similar groups can be presented by wreath recursions.
\begin{defn}\label{wreath}
     Write $\mathsf{X}=\{x_1,\ldots ,x_d\}$. The map 
\begin{equation*}
\begin{split}
&\psi: G\rightarrow S_d\wr G :=S_d\ltimes G^d\\
&\psi(g)=\sigma({}_{x_1}{|g},\ldots ,{}_{x_d}{|g})
\end{split}
\end{equation*}
where $S_d$ is the symmetric group on $\mathsf{X}$ and $\sigma\in S_d$, is called the \textit{wreath recursion} of $G$. We often write $g=\sigma(g_1,\ldots ,g_d)$ where $g_i={}_{x_i}{|g}$.
\end{defn}

Before presenting the examples, let us prove Theorem \ref{estimate} which applies to Subsections \ref{k110}---\ref{fggrp}. 

\subsection{Subexponential growth estimate for groups with a purely non-Hausdorff singularity} \label{53}
\begin{defn}\label{bridge} Let $G$ be a group of bounded type with finite cycles with a unique germ-defining singularity $\xi\in \Omega(\mathsf{B})$. Let $\xi_n$ be the $n$-th truncation of $\xi$. Let $H$ be the group of germs at $\xi$. Let $e_n$ be a non-identity subset of $H$, for $n\in \mathbb{N}$. A \textit{bridge of rank $n$}  at $\xi$, denoted $\mathcal{B}_{\xi,e_n}$ (abbr. $\mathcal{B}_{\xi_n}$ or $\mathcal{B}_n$), is a finite graph consisting of $|e_n|$ isomorphic copies of $\mathcal{T}_{\xi_{n}}$ and connecting them at $\xi_n$ by all the elements in $e_n$, while all elements in  $H\backslash e_n$ are loops at $\xi_n$. The elements in $e_n$ are called \textit{connectors}.
\end{defn}

In other words, a bridge is a union of several isomorphic copies of the same tile connected at the same boundary points. 

\begin{defn}\label{typeh} Let $w\in S^*$ be a word. For $h\in H$, a \textit{traverse of type $h$} of $\mathcal{T}_{\xi_m}$ is a subword of $w$ whose trajectory is a traverse of $\mathcal{T}_{\xi_m}$ starting at a boundary point of (an isomorphic copy in $\mathcal{T}_{\xi}$) $\mathcal{T}_{\xi_m}$ different from $\xi_m$, and ending at a boundary point of $\mathcal{T}_{\xi_m}$ different from $\xi_m$, and this subwalk is lifted to $\widetilde{\Gamma}_{\xi}$ to a walk starting in a branch of the graph of germs $\widetilde{\Gamma}_{\xi}$ corresponding to $Id\in H$, ending in a branch corresponding to $h\in H$. The set of all traverses of type $h$ of $\mathcal{T}_{\xi_m}$ for $h\in H$ is called \textit{traverses of type $H$}.
\end{defn}

\begin{thm}\label{estimate}  Let $G$ be a group of bounded type. Suppose the following conditions hold:
%generated by a finite set $S$ of finite order homeomorphisms acting faithfully and minimally on a Cantor set $\mathcal{X}$. 
\begin{enumerate}
      \item the group $G$ has a unique purely non-Hausdorff singularity $\xi\in\Omega(\mathsf{B})$;
      \item the group of germs $H=G_{\xi}/G_{(\xi)}$ is finite abelian with invariant domains;
      \item\label{partition} let $\mathcal{P}$ be a finite partition of $\Omega(\mathsf{B})\backslash\{\xi\}$. For every $h\in H$ there exists a piece $P\in\mathcal{P}$ accumulating on $\xi$ such that ${}_{P}{|h}={}_{P}{|Id}$, and given any piece $P\in\mathcal{P}$, there exists $h\in H$ such that ${}_{P}{|h}={}_{P}{|Id}$ (restriction on a subset, not to confuse with taking sections);      
    \item the orbital graph $\Gamma_{\xi}$ is thin and $1$-ended;
    \item The tiles of truncations of $\xi$ are strongly repetitive.
\end{enumerate}
Then $G$ has finitely many incompressible elements. 
\end{thm}
\rmk This theorem, together with Proposition \ref{genebd} below, is a partial generalization of \cite[Theorem 3.1]{BNZ}. The group $G$ is also periodic by \cite[Theorem 3.4.1]{JC19}. 
\begin{proof}
It is enough just to consider the graphs of traverses for each $\mathcal{T}_{\xi_n}$, since a traverse of other Write $\mathcal{P}=\{P_1,\ldots ,P_d\}$. 
%Recall from Subsection \ref{GOA} that
For each piece $P_i\in\mathcal{P}$, let $\{\zeta_n\}\subset P_i$ be a sequence of regular points in converging to $\xi$. Then the sequence of rooted graphs $(\Gamma_{\zeta_n},\zeta_n)$ converges to a limit graph $\Lambda_{P_i}$ (abbr. $\Lambda_{i}$). By assumption 3, for $n$ big enough, there exist some $h\in H$ such that the germs $(h,\zeta_n)$ are trivial, and for other $h'\neq Id \in H$, the germs $(h',\zeta_n)$ are non-trivial (by the definition of a singular point). Denote by  $H_i$ the group generated by $h\in H$ with trivial germs $(h,\zeta_n)=(Id,\zeta_n)$. Then $\Lambda_i$ is obtained by taking $|H/H_i|$ copies of $\Gamma_{\xi}$, each of which is identified with $\{h\}\times \Gamma_{\xi}$ for $h\notin H_i$, connecting them at $\xi$ by the Cayley graph of $H/H_i$, while elements in $H_i$ become loops at each $(h,\xi)$. Recall from Subsection \ref{GOA} that we then have a Galois covering $\lambda_i:\widetilde{\Gamma}_{\xi}\rightarrow \Lambda_{i}$ given by the rule 
\begin{align*}
   \lambda_i(h,v)= 
   \begin{cases}
       (h,v) \text{ if } h\notin H_i,\\
       (Id,v) \text{ if }h\in H_i,\\
   \end{cases}
\end{align*}
and also the Galois covering $\Lambda_i\rightarrow\Gamma_{\xi}$ that collapses all boundary edges to the loops at $\xi$. Since $\Lambda_i$ is the limit of orbital graphs of regular points, it follows from strong repetitivity that any finite central part $\Delta$ of $\Lambda_i$ forms a $\Delta$-sieve (Definitions \ref{substrrep}, \ref{strrep}) in any orbital graph. These finite subgraphs $\Delta$ contain isomorphic copies of bridges of rank $n$ with connectors $e_i$ for some $n\in\mathbb{N}$. Denote these bridges by $\mathcal{B}_{n,i}$ (not to confuse with the lower indices of tiles since here we only consider one type of tiles). Also by assumption 3, for each $h\in H$, there is a bridge $\mathcal{B}_{n,i}$ such that $h\notin e_i$. Fix $n$ big enough and let $l>0$ be the smallest integer such that the tile $\mathcal{T}_{\xi_n}$ contains at least one copy of all bridges $\mathcal{B}_{n-l,i}$, for $i=1,\ldots ,d$. Fixing $w\in G$. Hence a traverse induced by $w$ of type $h$ of $\mathcal{T}_{\xi_m}$ will induce traverses of type $h'$ of $\mathcal{T}_{\xi_{m-l}}$ by crossing a bridge $\mathcal{B}_{n-l,j}$. But this subwalk is lifted to a central part of $\widetilde{\Gamma}_{\xi}$ and mapped by $\lambda_{j'}$ to $\Lambda_{j'}$ that starts and ends in the same branch $\{Id\}\times \Gamma_{\xi}$. Hence this subwalk cannot cross the bridge $\mathcal{B}_{n-l,j'}$.  It follows from strong repetitivity that there is an isomorphic copy of a subwalk of $w$ in $\widetilde{\Gamma}_{\xi}$ that starts and ends in an isomorphic copy $(\mathcal{T}_{\xi_n},Id)$ and crosses the boundary edges of $\widetilde{\Gamma}_{\xi}$, for some $n<M$. Hence $w$ contains a return subword. Let $M$ be fixed. It follows that every word that has length bigger than $M$ will contain a return subword. The proof is complete. 
\end{proof}

To obtain an explicit subexponential growth estimate of $G$, we need the following. 

\begin{prop}\label{genebd} Let $G$ be the group in Theorem \ref{estimate}, and $\mathcal{P}=\{P_1,\ldots ,P_d\}$ be the partition in Condition \ref{partition} of Theorem \ref{estimate} . Then, for $\alpha>\frac{\log(\beta)}{\log(\beta)+\log(\frac{d}{d-1})}$, where $\beta=\lim\sup\limits_{n\rightarrow \infty}(\sum_i|\mathcal{T}_{\xi_n}|)^{1/n}$, the growth function $\gamma_G(R)\preccurlyeq \exp(R^{\alpha})$. 
 \end{prop}
 \begin{proof} 
 By the assumption, for each $h\in H$, there exists $e_{n-l,j}$ such that $h\in e_{n-l,j}$ and $e_{n-l,j'}$ such that $h\in H\backslash e_{n-l,j'}$, for $j\neq j'$ and $1\leq j,j'\leq d$. 
 %Consider a subgraph of the graph of traverses for $\mathcal{T}_{\xi_n}$, whose initial vertices is $\Theta_{\xi_n}$ and other vertices are $\Theta_{\xi_{n-l,j}}$, and $j=1,\ldots ,d$. 
 Let $w\in S^*$ be a word. Denote by $T_{n,h}$ the number of traverses of type $h$ of $\mathcal{T}_{\xi_n}$ for the word $w$. Also let $T_n=\sum\limits_{h\in H\backslash\{Id\}}T_{n,h}$. For $h\in H$, let $j_h$ be such that $h\notin e_{n-l,j_h}$. We have $$ T_{n} \leq \sum\limits_{h\in e_{n-l,j_h}}T_{n-l,h}.$$ Summing over $j_h$ and noticing that $j_h$ exhausts $1,\ldots ,d$, we have $$dT_n\leq \sum\limits_{j_h}\sum\limits_{h\in e_{n-l,j_h}}T_{n-l,h}\leq (d-1)T_{n-l}.$$ This inequality holds if we replace $n$ by $n-sl$ for every $s$ such that $n-(s+1)l\geq 1$. Hence $$dt^{-l}(F(t)-(T_0+T_1t+\ldots +T_{l-1}t^{l-1}))\leq (d-1)F(t)$$ $$\implies F(t)\leq |w|\dfrac{g(t)}{t^{-l}-\frac{d-1}{d}},$$ where $g(t)$ is positive. The last inequality holds since $T_n\leq |w|$ for all $n$. Taking $\Phi_{\xi}(t)=\dfrac{g(t)}{t^{-l}-\frac{d-1}{d}}$, then $F(t)\leq|w|\cdot\Phi_{\xi}(t)$ whenever $t<\sqrt[l]{\frac{d}{d-1}}$. By Proposition \ref{ultesti}, the conclusion follows. 
 \end{proof}

\subsection{A fragmentation of $\mathfrak{K}(11,0)$} \label{k110}

This example was communicated to the author by V.~Nekrashevych. Let $\mathsf{X}=\{0,1\}$. Consider the group $\mathfrak{K}(11,0)$ generated by 
\begin{align*}
    \begin{cases}
     a=\sigma,\\
     b=(1,a),\\
     c=(c,b).
    \end{cases}
\end{align*}

The action of $\mathfrak{K}(11,0)$ on $\mathsf{X}^{\omega}$ (the space of left infinite words over $\mathsf{X}$) is minimal since it is level-transitive.  
We can fragment the generator $c$ such that the point $\xi=0^{\omega}$ is a purely non-Hausdorff singularity. Define $W_n=\{0,1\}^{\omega}1\underbrace{00\ldots 0}_{n\text{-times}}$. Define pieces $$P_0=\bigcup_{k=0}^{\infty}W_{3k}\text{, }P_1=\bigcup_{k=0}^{\infty}W_{3k+1}\text{, }P_2=\bigcup_{k=0}^{\infty}W_{3k+2}.$$
Then the set of pieces $\mathcal{P}:=\{P_0,P_1,P_2\}$ forms a partition of $\{0,1\}^{\omega}\backslash\{\xi\}$ and they accumulate on $\xi$. Let $c_0$ act as $c$ on $P_0\bigcup P_1$, $c_1$ act as $c$ on $P_0\bigcup P_2$, and $c_2$ act as $c$ on $P_1\bigcup P_2$. In terms of wreath recursions, we have 
\begin{align*}
    \begin{cases}
     c_0=(c_1,b),\\
     c_1=(c_2,b),\\
     c_2=(c_0,1).
    \end{cases}
\end{align*}
Let $G:=\langle a,b,c_0,c_1,c_2 \rangle$. It follows that the point $\xi$ is the purely non-Hausdorff singularity of $G$. The group of germs $H=\langle c_0,c_1,c_2 \rangle\cong (\mathbb{Z}/2\mathbb{Z})^2$.

\begin{figure}[!htb]
    \centering
    \includegraphics[scale=.60]{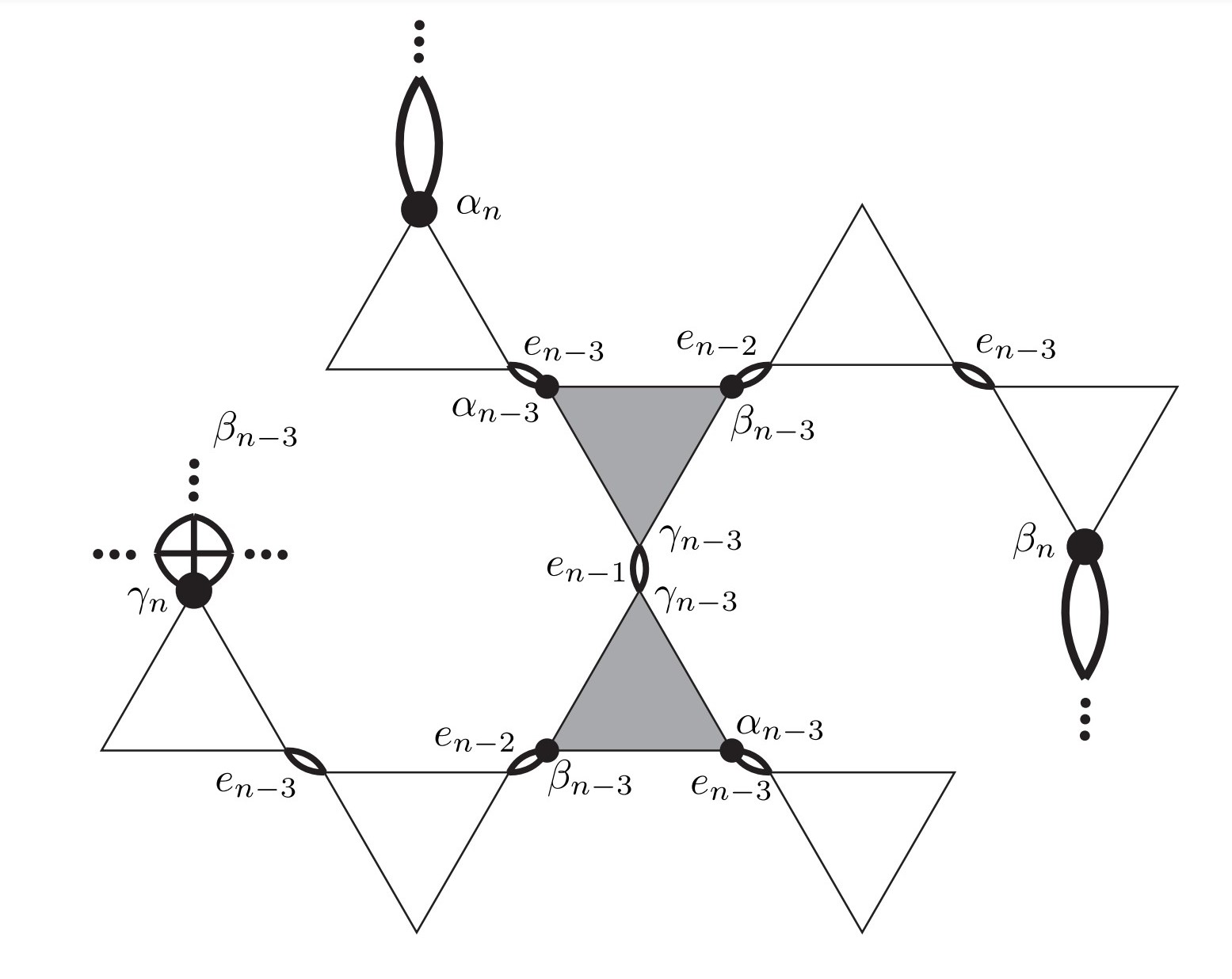}
    \caption{Constructing $\mathcal{T}_{n}$ from $\mathcal{T}_{n-3}$.}
    \label{ntile}
\end{figure}

Each finite tile $\mathcal{T}_n$ has three boundary points: $\gamma_n=0^{n}$, $\beta_n=10^{n-1}$, $\alpha_n=110^{n-2}$. Then $\gamma_{n+1}=0\gamma_{n}$, $\beta_{n+1}=1\gamma_{n}$, and $\alpha_{n+1}=1\beta_{n}$. The $(n+1)$-th level tile $\mathcal{T}_{n+1}$ is constructed as follows. Take two copies $0\mathcal{T}_n $ and $1\mathcal{T}_n $ of $\mathcal{T}_n$ and connect $0\alpha_n =0110^{n-2}$ with $1\alpha_n =1110^{n-2}$ by the an edge $e_n$ with labels
\begin{equation*} \label{connector}
\begin{split}
e_0 & = \frac{ \quad \quad a \quad \quad }{}, \\
e_1 & = \frac{\quad \quad b \quad \quad }{}, \\
e_{3k-1} & = \frac{ \ \quad c_0,c_1 \quad \ }{},\\
e_{3k} & = \frac{ \ \quad c_0,c_2 \quad \ }{},\\
e_{3k+1} & = \frac{ \ \quad c_1,c_2 \quad \ }{},
\end{split}
\end{equation*}
for $k\geq 1$. Figure \ref{ntile} shows the construction of $\mathcal{T}_{n}$ from $\mathcal{T}_{n-3}$. The relative positions of the boundary points $\gamma_n$, $\beta_n$, $\alpha_n$ were obtained using model graphs of the hull graph $[\alpha,\beta,\gamma]_n$. See \cite[Example 5.1.2]{JC19}. 
The orbital graph of $\xi=0^{\omega}$, denoted $\Gamma_{\xi}$, is the inductive limit of $\mathcal{T}_{n}\mapsto 0\mathcal{T}_{n}\subset \mathcal{T}_{n+1}$. By Proposition \ref{strrepp}, each finite tile is strongly repetitive. 

The graph of germs $\widetilde{\Gamma}_{\xi}$ is obtained by taking $|H|$ ($=4$) copies of $\Gamma_{\xi}$ connected by the Cayley graph of $H$. Let $h\in H$. A \textit{rank $n$ traverse of type $h$} is a walk in $\widetilde{\Gamma}_{\xi}$ starting from a point in the set $\{(Id, 0^{\omega}\alpha_n),(Id,0^{\omega}\beta_n)\}$ ending in a point in the set $\{(h,0^{\omega}\alpha_n ),(h,0^{\omega}\beta_n)\}$ and touching these points only once at the beginning and the end of the walk. Let $g\in G$ be a word such that it is a rank $n$ traverse. At the last time that $g$ crosses the Cayley graph of $H$, the remainder of the walk will be a traverse of an isomorphic copy of the tile $\mathcal{T}_n$ in $\widetilde{\Gamma}_{\xi}$ starting in $\gamma_n$ and ending in one of the points in the set $\{(h,0^{\omega}\alpha_n ),(h,0^{\omega}\beta_n)\}$ without touching any of them in between. This moment is surrounded by the last moment the walk entered a copy of $\mathcal{T}_{n-3}$ and the first moment it exited the other copy of $\mathcal{T}_{n-3}$ adjacent to the connector $e_{n-1}$. These two copies are highlighted on Figure \ref{ntile}. Hence we have a traverse of rank $n-3$ (when lifting these two copies of $\mathcal{T}_{n-3}$ to the central part of $\widetilde{\Gamma}_{\xi}$). 

Denote by $T_{n,h}$ the number of rank $n$ traverses of type $h$, and by $T_n$ the total number of rank $n$ traverses, i.e., $T_n=\sum\limits_{h\in\{c_0,c_1,c_2\}}T_{n,h}$. By the argument above, we have $$T_n\leq \sum\limits_{h\in e_{n-1}} T_{n-3,h}.$$ By similar arguments, we have $$T_{n+1}\leq \sum\limits_{h\in e_{n}} T_{n-2,h}$$ and $$T_{n+2}\leq \sum\limits_{h\in e_{n+1}} T_{n-1,h}.$$

Consider the generating functions $$F_h(t)=\sum\limits_{n=0}^{\infty}T_{n,h}t^n,\text{   } F(t)=\sum\limits_{h\in\{c_0,c_1,c_2\}}F_h(t)\text{ for } t>0.$$ Since $T_n\leq |g|$, we have
$$t^{-i}(F(t)-C_i(t)|g|)\leq \sum\limits_{h\in e_{n-4+i}}F_h(t) \text{ for }i=3,4,5,$$ for some polynomials $C_i$ independent of $g$ and $t$. Notice that each $h\in \{c_0,c_1,c_2\}$ appears exact twice in the connectors $e_{n-1},e_n,e_{n+1}$. Add up the inequalities and we have
$$(t^{-3}+t^{-4}+t^{-5})(F(t))-|g|C(t)\leq 2F(t).$$ Hence
\begin{equation}\label{ineq1}
F(t)\leq |g|\dfrac{C(t)}{t^{-5}+t^{-4}+t^{-3}-2}    
\end{equation}
for some function $C(t)>0$ independent of $g$ and $t>0$ such that the denominator on the right-hand side is positive. Solving the equation $x^5+x^4+x^3-2=0$, we get a real positive root $\eta\approx 0.9028$. It follows that for all $t<\eta^{-1}$ the right hand side of (\ref{ineq1}) is positive. By Theorem \ref{estimate} and \cite[Proposition 2.7]{BNZ}, the growth function of this group is dominated by $\exp(R^{\alpha})$ for every $\alpha>\frac{\log 2}{\log 2-\log \eta}\approx 0.8715$. 

\begin{figure}
    \centering
    \includegraphics[scale=.50]{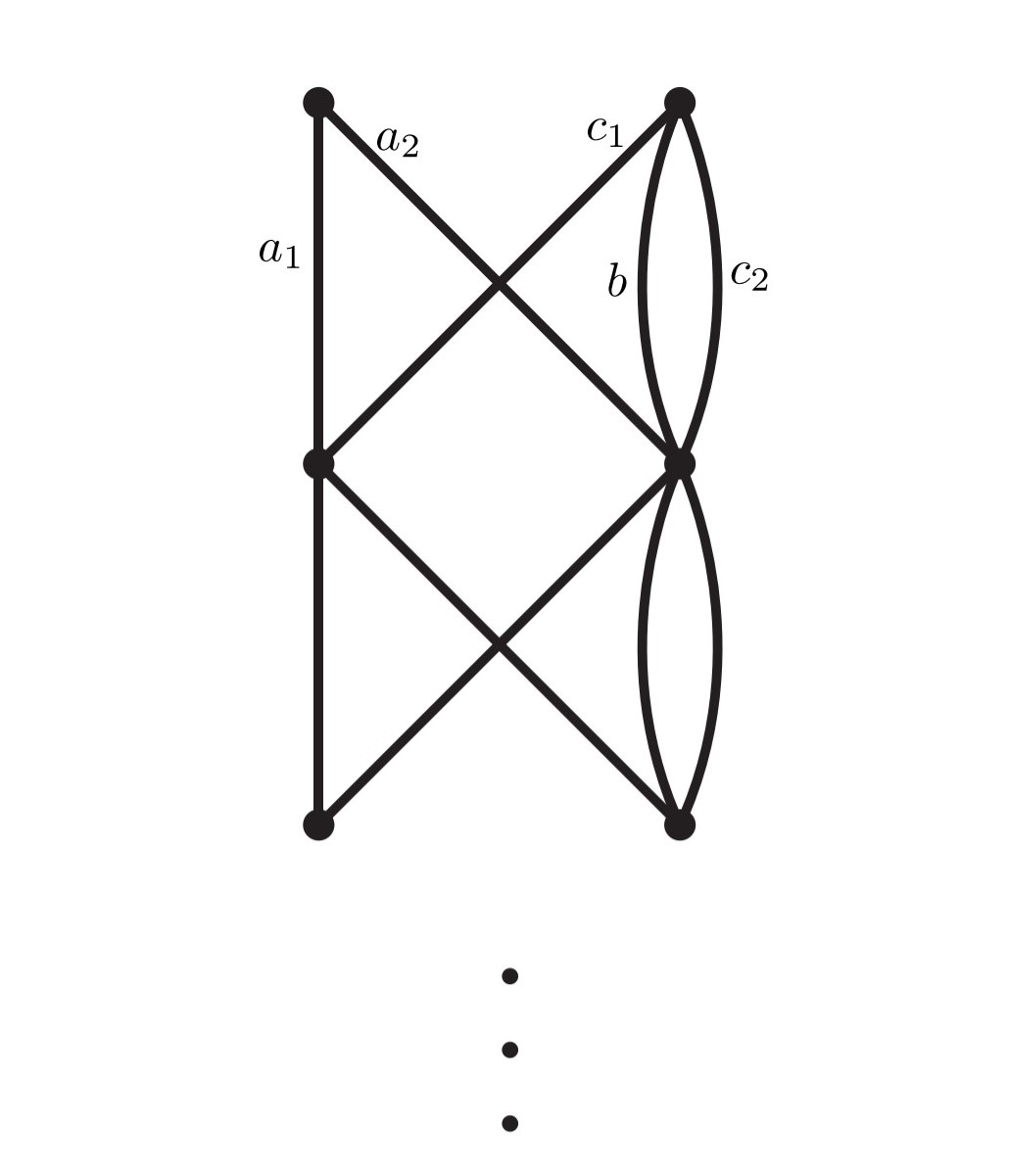}
    \caption{Bratteli Diagram $\mathsf{B}$ related to Penrose tiling.}
    \label{Brat}
\end{figure}

\subsection{A fragmentation of a group related to Penrose tiling}\label{2tiles}
Let $\mathsf{B}$ be the Bratteli diagram shown in Figure \ref{Brat}. 
%See \cite[Section 3]{BGN} for a description of Penrose tiling and a semigroup acting on it. 
Let $V_n=V=\{v_1,v_2\}$ be the set of $n$-th level vertices of $\mathsf{B}$. Let the set of edges $E_n=E=\{a_1,a_2,c_1,b,c_2\}$.
%Let $\Omega(\mathsf{B})$ and $\Omega_n(\mathsf{B})$ be as in Definition \ref{brat}. 
Let $G_0$ be an inverse semigroup acting on $\Omega(\mathsf{B})$ defined by the following rules of tile inflation with the singular point $c_2^{\omega}$.

\begin{figure}[!htb]
\centering
\begin{minipage}{.5\textwidth}
  \centering
  \includegraphics[width=\linewidth]{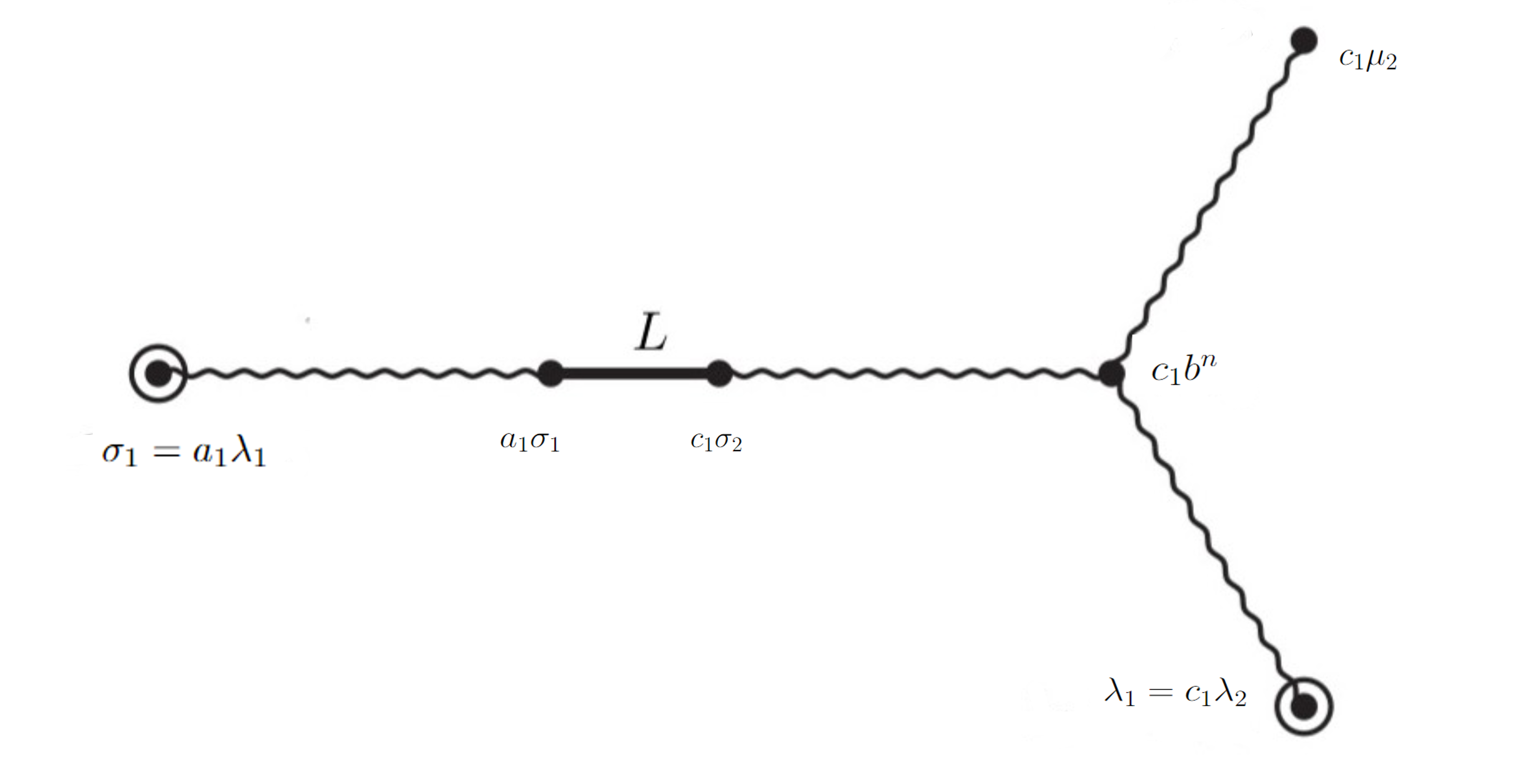}
  \subcaption{$\mathcal{T}_{1,n+1}$}
  \label{figT1}
\end{minipage}%
\begin{minipage}{.5\textwidth}
  \centering
  \includegraphics[width=\linewidth]{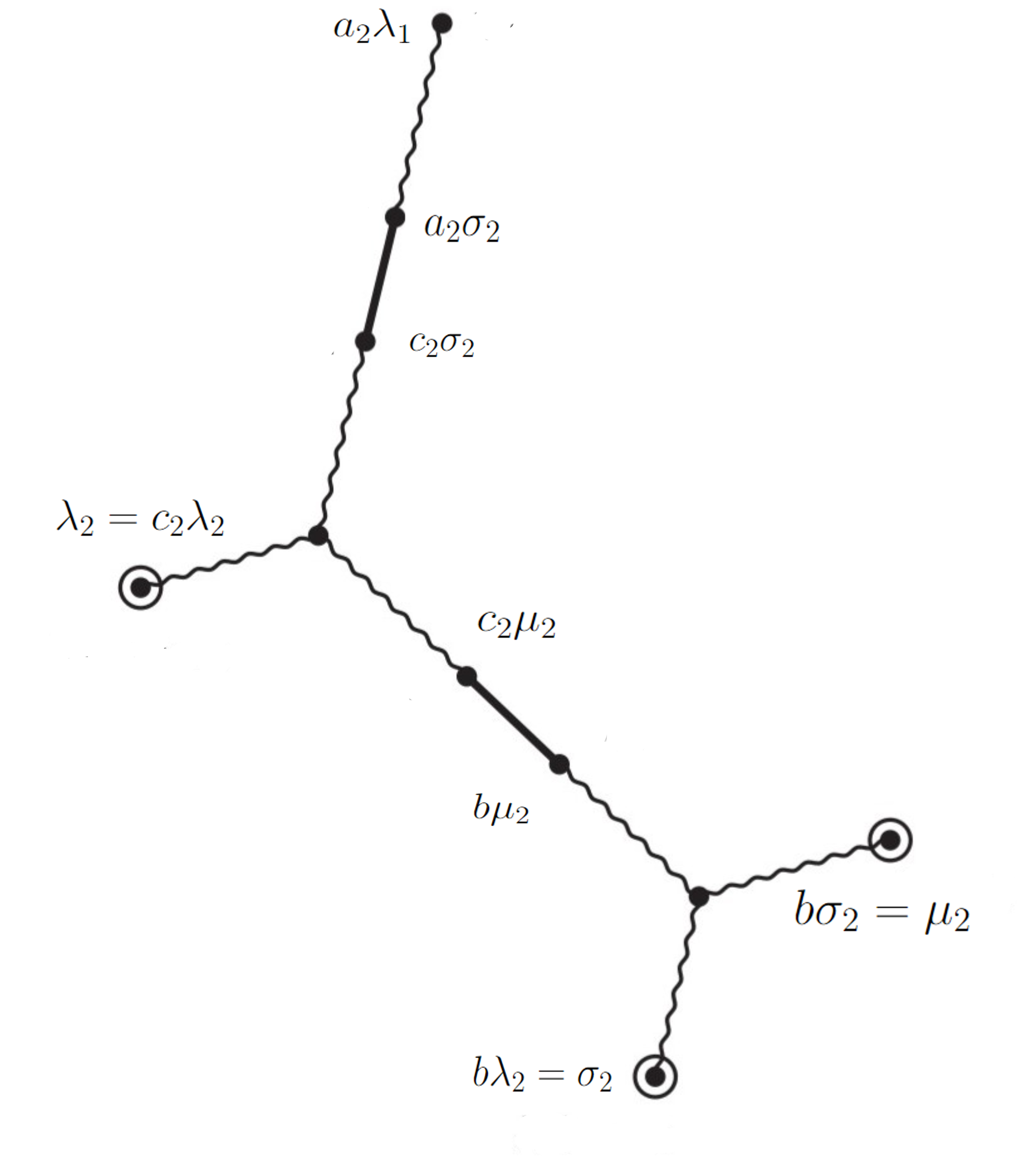}
  \subcaption{$\mathcal{T}_{2,n+1}$}
  \label{figT2}
\end{minipage}
\caption{Constructing the $(n+1)$-th level tiles from the $n$-th level tiles.}
\label{T1T2}
\end{figure}

Let $\mathcal{T}_{1,n}$ and $\mathcal{T}_{2,n}$ be the two types of $n$-th level tiles whose edges are two-sided arrows. Assume that they are connected for $n\leq 2$. For $n\geq 3$, their boundary points are $\partial\mathcal{T}_{1,n}=\{\lambda_{1,n},\sigma_{1,n}\}$ and $\partial\mathcal{T}_{2,n}=\{\lambda_{2,n},\sigma_{2,n},\mu_{2,n}\}$, where $\lambda_{1,n}=c_1c_2^{n-1}$, $\sigma_{1,n}=c_1a_1c_2^{n-2}$, $\lambda_{2,n}=c_2^{n}$, $\sigma_{2,n}=bc_2^{n-1}$ and $\mu_{2,n}=bbc_2^{n-2}$. In particular, we have $\lambda_{1,n}=c_1\lambda_{2,n-1}$, $\sigma_{1,n}=a_1\lambda_{1,n-1}$, $\lambda_{2,n}=c_2\lambda_{2,n-1}$, $\sigma_{2,n}=b\lambda_{2,n-1}$, and $\mu_{2,n}=b\sigma_{2,n-1}$. The tile $\mathcal{T}_{1,n+1}$ is obtained by taking a copy $a_1\mathcal{T}_{1,n}$ of $\mathcal{T}_{1,n}$ and a copy $c_1\mathcal{T}_{2,n}$ and connect the points $a_1\sigma_1$ with $c_1\sigma_2$ by an edge labelled by $L$. The tile $\mathcal{T}_{2,n+1}$ is obtained by taking two copies $b\mathcal{T}_{2,n}$ and $c_2\mathcal{T}_{2,n}$ of $\mathcal{T}_{2,n}$ and connect the point $b\mu_2$ to $c_2\mu_2$ by an edge also labelled by $L$. Also connect the point $c_2\sigma_2$ on $c_2\mathcal{T}_{2,n}$ with the point $a_2\sigma_1$ on $a_2\mathcal{T}_{2,n}$ labelled by $L$. See Figure \ref{bdorb} for a description of the continuations of boundary points and locations of connecting points (the lower indices $n$ are dropped). See Figure \ref{T1T2} for the construction of the $(n+1)$-th level tiles from the $n$-th level tiles, where the points highlighted by circles are boundary points and all connections are highlighted by thick lines. Since all tiles are connected to each other as isomorphic copies on deeper levels, by Proposition \ref{minact}, the action of $G_0$ on $\Omega(\mathsf{B})$ is minimal.

\begin{prop}
 The largest distance between boundary points on both $\mathcal{T}_{1,n}$ and $\mathcal{T}_{2,n}$ grow exponentially with respect to $n$. 
\end{prop}
\begin{proof}
Let us prove the statement for $\mathcal{T}_{2,n}$. The statement for $\mathcal{T}_{1,n}$ will follow from a similar argument. Let us estimate the distance the distance from $\lambda_{2,n}$ to $\sigma_{2,n}$, denoted $d(\lambda_{2,n},\sigma_{2,n})$. It follows from Figure \ref{figT2} that
\[d(\lambda_{2,n},\sigma_{2,n})=d(\lambda_{2,n},c_2\mu_{2,n-1})+d(b\mu_{2,n-1},\sigma_{2,n})+1 \addtag\label{eq1}\] where 
\[d(\lambda_{2,n},c_2\mu_{2,n-1})=d(\lambda_{2,n-1},\mu_{2,n-1})\addtag\label{eq2}\] and 
\[d(b\mu_{2,n-1},\sigma_{2,n})=d(\mu_{2,n-1},\lambda_{2,n-1})\addtag\label{eq3}.\]

The right-hand side of (\ref{eq2}) and (\ref{eq3}) are the distance between $\mu_2$ and $\lambda_2$ in $\mathcal{T}_{2,n-1}$, which is equal to 
\[d(\lambda_{2,n-2},\mu_{2,n-2})+d(\mu_{2,n-2},\sigma_{2,n-2})+1,\addtag\label{eq4}\]
while 
\[d(\mu_{2,n-2},\sigma_{2,n-2})=d(\sigma_{2,n-3},\lambda_{2,n-3}).\addtag\label{eq5}\]
It follows that the distance $d(\lambda_{2,n},\sigma_{2,n})$ satisfies the recursive formulas determined in (\ref{eq1})---(\ref{eq5}) for all $n\geq 4$, and for $1\leq n\leq 3$ the distances are at least $1$. Hence the sequence $\{d(\lambda_{2,n},\sigma_{2,n})\}_{n\in\mathbb{N}}$ grows exponentially. Since $d(\lambda_{2,n},\sigma_{2,n})$ is bounded above by the largest distance between points in $\mathcal{T}_{2,n}$, which is bounded above by $|E|^n=5^n$, this completes the proof. 
\end{proof}

By Propositions \ref{expdiam}, \ref{linlin}, we have the following.

\begin{cor}
The action of $G_0$ on $\Omega(\mathsf{B})$ is linearly repetitive and thus $G_0$ is an inverse semigroup of bounded type. 
\end{cor}

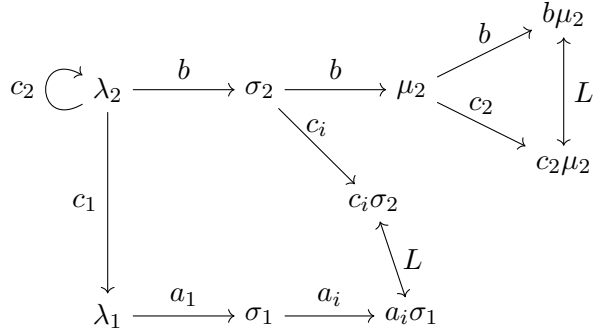
\begin{figure}[!htb]
\centering
\begin{tikzpicture}
\begin{scope}state without output/.append style={draw=none}%[every node/.style={circle,thick,draw}]
    \node (A) at (0,0) {$\lambda_2$};
    \node (B) at (2,0) {$\sigma_2$};
    \node (C) at (4,0) {$\mu_2$};
    \node (E) at (6,-1) {$c_2\mu_2$};
    \node (F) at (6,1) {$b\mu_2$};
    \node (G) at (0,-3) {$\lambda_1$};
    \node (H) at (2,-3) {$\sigma_1$};
    \node (I) at (4,-3) {$a_i\sigma_1$};
    \node (J) at (3.5,-1.5) {$c_i\sigma_2$};
\end{scope}

\begin{scope} %[>={Stealth[black]},
              every node/.style={fill=white,circle},
              every edge/.style={draw=black,very thick}]
              
    %\path [->] (A) edge [bend left=40] node {} (B);
    %\path [->] (B) edge [bend left=40] node {} (A);
    %\draw [->] (B) to  [out=320,in=40,looseness=9] (B);
    %\path [->] (D) edge node {$3$} (C);
    %\path [->] (A) edge node {$3$} (E);
    %\path [->] (D) edge node {$3$} (E);
    %\path [->] (D) edge node {$3$} (F);
    %\path [->] (C) edge node {$5$} (F);
    %\path [->] (E) edge node {$8$} (F); 
    %\path [->] (B) edge[bend right=60] node {$1$} (E); 

    %\draw[]  (A) node[draw=none][midway,above] {$b_1$} (B);
    \path [->] (A) edge node [midway,above] {$b$} (B);
    \path [->] (B) edge node [midway,above] {$b$} (C);
    \path [->] (C) edge node [midway,above] {$b$} (F);
    \path [->] (C) edge node [midway,above] {$c_2$} (E);
    \path [<->] (E) edge node [midway,right] {$L$} (F);
    \path [->] (B) edge node [midway,above] {$c_i$} (J);
    \path [->] (A) edge node [midway,left] {$c_1$} (G);
    \path [->] (G) edge node [midway,above] {$a_1$} (H);
    \path [->] (H) edge node [midway,above] {$a_i$} (I);
    \path [<->] (I) edge node [midway,right] {$L$} (J);
    %\draw [->] (A) to  [out=90,in=140,looseness=8] edge node [midway,left] {$c_2$} (A);
    \path [->] (A) edge [out=210, in=150, looseness=5, "$c_2$"] (A);
    
\end{scope}

\end{tikzpicture}
\caption{Continuations of boundary points and locations of connecting points, where $i=1,2$.}
\label{bdorb}
\end{figure}

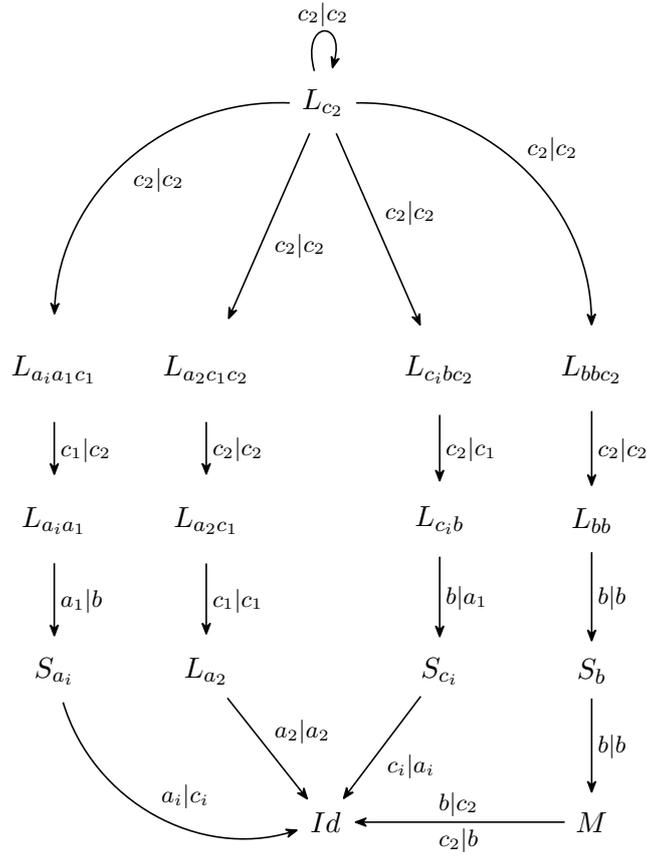
\begin{figure}
\centering
\begin{tikzpicture}[->,>={Stealth[round]},shorten >=1pt,%
                    auto,node distance=2cm,on grid,semithick,
                    inner sep=2pt,bend angle=45, state/.style={circle,inner sep=2pt}]
  \node[state]         (A)                    {$L_{c_2}$};
  \node[state]         (B) [below left=5cm  of A]  {$L_{a_ia_1c_1}$};
  \node[state]         (C) [right=  of B]       {$L_{a_2c_1c_2}$};
  \node[state]         (E) [below right=5cm of A]       {$L_{bbc_2}$};
  \node[state]         (D) [left=of E]       {$L_{c_ibc_2}$};
 
  \node[state]         (F) [below = of B]   {$L_{a_ia_1}$};
  \node[state]         (G) [below = of F]   {$S_{a_i}$};
  
  \node[state]         (H) [below = of C]   {$L_{a_2c_1}$};
  \node[state]         (I) [below = of H]   {$L_{a_2}$};

  \node[state]         (J) [below = of D]   {$L_{c_ib}$};
  \node[state]         (K) [below = of J]   {$S_{c_i}$};
  
  \node[state]         (L) [below = of E]   {$L_{bb}$};
  \node[state]         (M) [below = of L]   {$S_b$};
  \node[state]         (N) [below = of M]   {$M$}; 

  \node[state]         (O) [left = 3.5cm of N] {$Id$};

  \path [every node/.style={font=\footnotesize}]
        (A) edge [loop above] node {$c_2|c_2$} (A)
            edge [bend right] node {$c_2|c_2$} (B)
            edge              node {$c_2|c_2$} (C)
            edge              node {$c_2|c_2$} (D)
            edge [bend left]  node {$c_2|c_2$} (E)
            
        (B) edge              node {$c_1|c_2$} (F)
        (F) edge              node {$a_1|b$} (G)
        (G) edge [bend right] node {$a_i|c_i$} (O)

        (C) edge              node {$c_2|c_2$} (H)
        (H) edge              node {$c_1|c_1$} (I)
        (I) edge              node {$a_2|a_2$} (O)

        (D) edge              node {$c_2|c_1$} (J)
        (J) edge              node {$b|a_1$}   (K)
        (K) edge              node {$c_i|a_i$} (O)

        (E) edge              node {$c_2|c_2$} (L)
        (L) edge              node {$b|b$}     (M)
        (M) edge              node {$b|b$}     (N)
        (N) edge              node [above] {$b|c_2$} node [below] {$c_2|b$} (O);

        %(E) edge [bend left]  node {1,0,R} (A);

\end{tikzpicture}
\caption{Partial minimization automaton $\mathcal{A}'$ generating $G_0$ as an inverse semigroup, for $i=1,2$. }
\label{automag0}
\end{figure}
Figure \ref{automag0} shows the minimization automaton $\mathcal{A}'$ generating the inverse semigroup $G_0$. Notice that the subscripts of the states are read from right to left. Each state on the figure is also an initial state. The state $L_{c_2}$ is the unique directed state, and its section \[{}_{c_2}{|L_{c_2}}=\{L_{c_2},L_{a_1a_1c_1},L_{a_2a_1c_1},L_{a_2a_1c_2},L_{c_1bc_2},L_{c_2bc_2},L_{bbc_2}\},\] while all other states are finitary, and their sections are all identified with singletons. Table \ref{cyclab} is the list of states to be relabeled under the same labeling. Different sets of states are separated by double horizontal lines. All the states in the same set have disjoint domains. The state $M$ is only defined at $v_2$ and takes values at $v_2$. Hence we can extend $M$ identically to the rest of $\Omega(\mathsf{B})$. For the label $S$, we extend it identically on $\Omega(\mathsf{B})b$. For the label $L$, we extend it identically to the rest of $\Omega(\mathsf{B})$ and relabel $L_{c_2}$ to be $L$ since it has order $2$ and has disjoint domain and range from those of all other $L$-transformations. In this way, the domain of $L$ forms a partition of $\Omega(\mathsf{B})$. It follows that the group $G_0$ is generated by the set $\{L,M,S\}$. This group was studied in \cite[Section 4.5]{JNS16}.  
\begin{table}[htb!]
\centering
\begin{tabular}{|c | c | c | c|} 

 \hline
 States & Domain 1 & Domain 2 & Relabel \\ [0.5ex] 
 \hline
 $S_{a_i},S_{c_i}$ & $\Omega(\mathsf{B})a_i$ & $\Omega(\mathsf{B})c_i$ & S \\ 
 \hline\hline
 $L_{a_ia_1},L_{c_ib}$ & $\Omega(\mathsf{B})a_ia_1$ & $\Omega(\mathsf{B})c_ib$ & L \\

 $L_{a_ia_1c_1},L_{c_ibc_2}$ & $\Omega(\mathsf{B})a_ia_1c_1$ & $\Omega(\mathsf{B})c_ibc_2$ &  \\

$L_{c_2bb},L_{bbb}$ & $\Omega(\mathsf{B})c_2bb$ & $\Omega(\mathsf{B})bbb$ & \\

$L_{c_2bbc_2},L_{bbbc_2}$ & $\Omega(\mathsf{B})c_2bbc_2$ & $\Omega(\mathsf{B})bbbc_2$ & \\
\hline
\end{tabular}
\caption{List of states with disjoint domians and relabels. }\label{cyclab}
\end{table}

Now we can fragment the generator $L$ such that the point $\xi=c_2^{\omega}$ is the purely non-Hausdorff singularity. Consider the set $\{b,c_2\}^{\omega}\backslash \{\xi\}\subset\Omega(\mathsf{B})\backslash \{\xi\}$. Let $\mathcal{P}_{\xi}$ be the following partition of $\{b,c_2\}^{\omega}\backslash \{\xi\}$ analogous to the partition of $\{0,1\}^{\omega}\backslash \{0^{\omega}\}$ in Example \ref{k110}. Define $W_n=\{b,c_2\}^{\omega}b \underbrace{c_2c_2\ldots c_2}_{n\text{-times}}$. Define pieces  $$P_0=\bigcup_{k=0}^{\infty}W_{3k}\text{, }P_1=\bigcup_{k=0}^{\infty}W_{3k+1}\text{, }P_2=\bigcup_{k=0}^{\infty}W_{3k+2}.$$ 
Then the set of pieces $\mathcal{P}_{\xi}:=\{P_0,P_1,P_2\}$ forms a partition of $\{b,c_2\}^{\omega}\backslash\{\xi\}$ and they accumulate on $\xi$. Let $L_0$ act as $L$ on $P_0\bigcup P_1$, $L_1$ act as $L$ on $P_0\bigcup P_2$, and $L_2$ act as $L$ on $P_1\bigcup P_2$. Let $\mathcal{P}$ be any partition in three pieces of $\Omega(\mathsf{B})\backslash \{\xi\}$ such that each piece in $\mathcal{P}$ contains a piece in $\mathcal{P}_{\xi}$. Notice that although the fragmentation of $L$ is partially defined by the action on a subset, it can be extended to the whole $\Omega(\mathsf{B})\backslash \{\xi\}$ by extending the actions of $L_i$, $i=0,1,2$, on larger pieces. This corresponds to introducing multiple edges to the connections on each tile. Let $G$ be the group generated by $L_1,L_2,L_3$ and other generators of $G_0$. It then follows that the point $\xi$ is the purely non-Hausdorff singularity of $G$. The group of germs $H=\langle L_0,L_1,L_2 \rangle\cong (\mathbb{Z}/2\mathbb{Z})^2$. We can apply Theorem \ref{estimate} to obtain an estimate of upper bound for the growth function of $G$. 

Denote by $\mathcal{T'}_{i,n}$ the tiles of $G$ for $i=1,2$ and $n\geq 3$. The orbital graph of $\xi$, denoted $\Gamma_{\xi}$, is the inductive limit of $\mathcal{T'}_{2,n}\mapsto c_2\mathcal{T'}_{2,n}\subset\mathcal{T'}_{2,n+1}$. Notice that there are two places of connections on the unfragmented $\mathcal{T}_{2,n+1}$:
$$c_2\sigma_{2,n}\frac{ \quad \quad L \quad \quad }{}a_2\sigma_{2,n}$$ and $$c_2\mu_{2,n}\frac{ \quad \quad L \quad \quad }{}b\mu_{2,n}.$$

We can just consider the traverses of $\mathcal{T}'_{2,n}$, since any traverse of $\mathcal{T}'_{1,n}$ induced by $\mathcal{T}'_{2,n+1}$ necessarily contains a return. 
The construction of $\mathcal{T'}_{2,n}$ from $\mathcal{T'}_{2,n-3}$ looks exactly the same as Figure \ref{ntile}, up to a relabelling of boundary points and connections. The connectors for $k\geq 1$ are
\begin{equation*} \label{connector}
\begin{split}
e_{3k-1} & = \frac{ \ \quad L_0,L_1 \quad \ }{},\\
e_{3k} & = \frac{ \ \quad L_0,L_2 \quad \ }{},\\
e_{3k+1} & = \frac{ \ \quad L_1,L_2 \quad \ }{}. 
\end{split}
\end{equation*}
Hence by a similar argument as in Subsection \ref{k110}, the growth function for the group $G$ is dominated by $\exp(R^{\alpha})$ for every $\alpha>\frac{\log \beta}{\log\beta-\log\eta}$, where $\beta=\limsup_{n\rightarrow\infty}|\mathcal{T'}_{2,n}|^{1/n}$ and $\eta\approx 0.9028$ is the real positive root of the equation $x^5+x^4+x^3-2=0$. Let us calculate $\beta$. From Figure \ref{T1T2}, we see that $|\mathcal{T'}_{2,n}|=2|\mathcal{T'}_{2,n-1}|+|\mathcal{T'}_{1,n-1}|$ and $|\mathcal{T'}_{1,n}|=|\mathcal{T'}_{1,n-1}|+|\mathcal{T'}_{2,n-1}|$. Then 
\begin{equation*} \label{sizeT2n}
\begin{split}
|\mathcal{T'}_{2,n}| & = |\mathcal{T'}_{1,n-1}|+2|\mathcal{T'}_{2,n-1}|\\
 & = 3|\mathcal{T'}_{1,n-2}|+5|\mathcal{T'}_{2,n-2}|\\ 
 & = 8|\mathcal{T'}_{1,n-3}|+13|\mathcal{T'}_{2,n-3}|\\
 & =\\
 &\ldots \\
 & = F_{2n} |\mathcal{T'}_{1,1}|+F_{2n+1}|\mathcal{T'}_{2,1}|,\\
\end{split}
\end{equation*}
where $F_k$ is the $k$-th term of the Fibonacci sequence. Then 
$$\beta=\lim_{k\rightarrow\infty}(2F_{2k})^{1/k}=\frac{3+\sqrt{5}}{2}.$$
Hence for every $\alpha>\frac{\log \beta}{\log\beta-\log\eta}\approx 0.9040$, the growth function is dominated by $\exp(R^{\alpha})$.

\subsection{A fragmentation of the Fabrykowski-Gupta group} \label{fggrp} 
Let $\mathsf{X}=\{0,1,2\}$. Denote by $G$ the Fabrykowski-Gupta generated by the following wreath recursions
\begin{align*}
    \begin{cases}
    a=(012),\\
    b=(a,Id,b).\\
    \end{cases}
\end{align*}
Both generators have order $3$, but the element $ba$ has infinite order. By \cite[Theorem 1]{BS01}, the growth function $\gamma_G(R)$ satisfies  $$\exp({R^{\frac{\log 3}{\log 6}}})\preccurlyeq\gamma_G(R)\preccurlyeq \exp({\frac{R(\log\log R)^2}{\log R}}),$$ and thus the group $G$ has intermediate growth but does not have a bounded power in the exponent, since   
\begin{align*}
\lim_{R\rightarrow \infty}\dfrac{R^{\alpha}}{\frac{R(\log\log R)^2}{\log R}} =
    \begin{cases}
   0, \text{ if $\alpha\in(0,1)$,}\\
   \infty, \text{ if $\alpha=1.$}\\
    \end{cases}
\end{align*}

\begin{figure}[!htb]
    \centering
    \includegraphics[scale=.50]{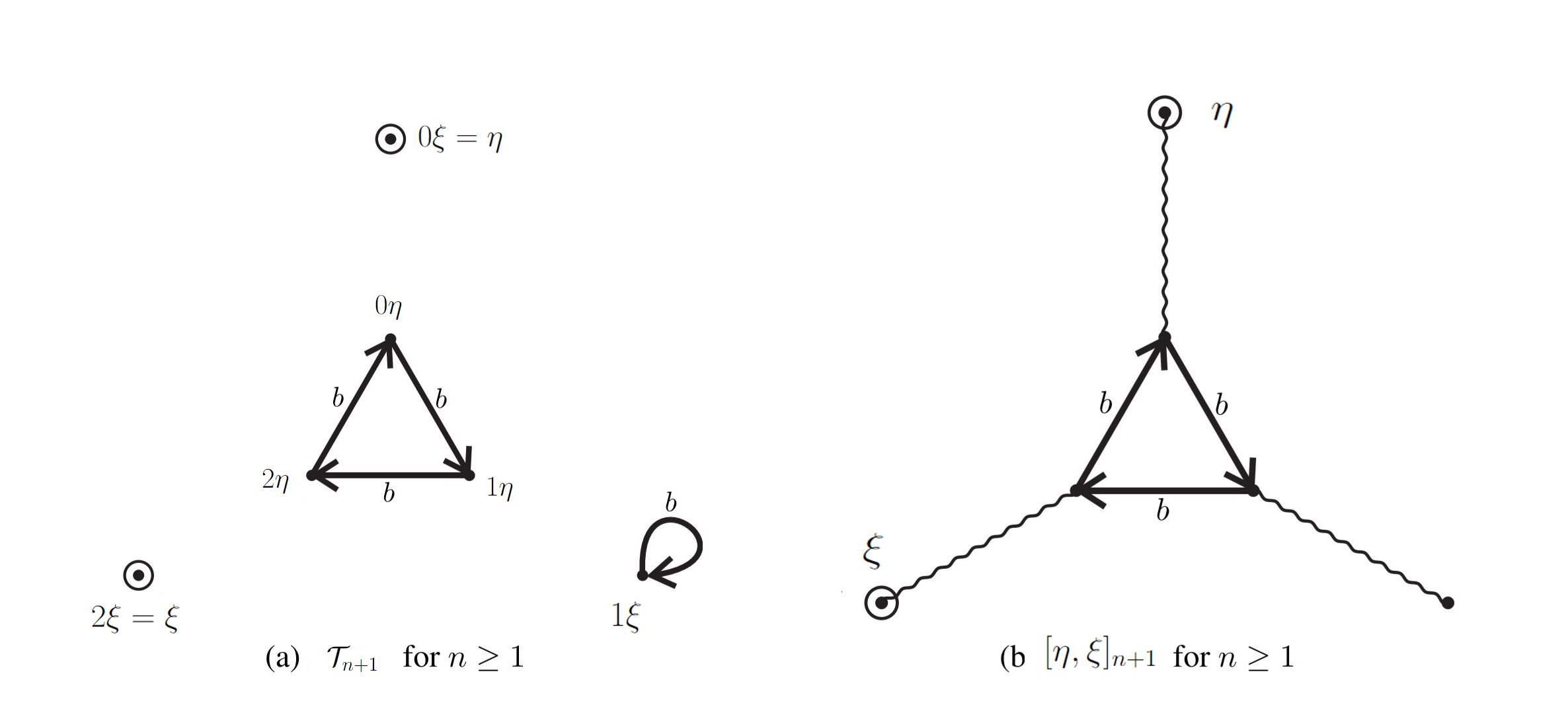}
    \caption{Tile inflations for the Fabrykovski-Gupta group.}
    \label{fgtile}
\end{figure}

The action of $G$ on $\mathsf{X}^{\omega}$ is minimal since it is level-transitive. Tile inflations of this group is shown in Figure \ref{fgtile}. By Proposition \ref{strrepp}, each finite tile is strongly repetitive. 
\\
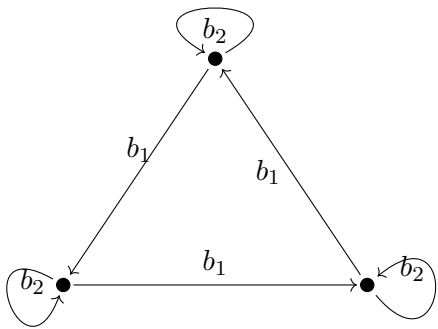
\begin{figure}
\centering
\begin{minipage}{.5\textwidth}
\begin{tikzpicture}
\begin{scope}state without output/.append style={draw=none}%[every node/.style={circle,thick,draw}]
    \node (A) at (3,0) {};
    \node (B) at (1,-3) {};

    \node (C) at (5,-3) {};
   
\end{scope}

\begin{scope} %[>={Stealth[black]},
              every node/.style={fill=white,circle},
              every edge/.style={draw=black,very thick}]
              
    %\path [->] (A) edge [bend left=40] node {} (B);
    %\path [->] (B) edge [bend left=40] node {} (A);
    %\draw [->] (B) to  [out=320,in=40,looseness=9] (B);
    %\path [->] (D) edge node {$3$} (C);
    %\path [->] (A) edge node {$3$} (E);
    %\path [->] (D) edge node {$3$} (E);
    %\path [->] (D) edge node {$3$} (F);
    %\path [->] (C) edge node {$5$} (F);
    %\path [->] (E) edge node {$8$} (F); 
    %\path [->] (B) edge[bend right=60] node {$1$} (E); 

    %\draw[]  (A) node[draw=none][midway,above] {$b_1$} (B);
    \path [->] (A) edge node [midway,above] {$b_2$} (B);
   
    \path [->] (B) edge node [midway,above] {$b_2$} (C);
    \path [->] (C) edge node [midway,left] {$b_2$} (A);
   
    %\draw [->] (A) to  [out=90,in=140,looseness=8] edge node [midway,left] {$c_2$} (A);
    \path [->] (A) edge [out=30, in=150, looseness=15, "$b_1$"] (A);
     \path [->] (B) edge [out=150, in=250, looseness=15, "$b_1$"] (B);
      \path [->] (C) edge [out=310, in=40, looseness=15, "$b_1$"] (C);
      
\filldraw[black] (3,0) circle (2.5pt) node[anchor=west]{};
\filldraw[black] (1,-3) circle (2.5pt) node[anchor=west]{};
\filldraw[black] (5,-3) circle (2.5pt) node[anchor=west]{};
\end{scope}
\end{tikzpicture}
\subcaption*{Connector (1): ker$(\pi_i)=\langle b_1 \rangle$}
\end{minipage}
\begin{minipage}{.5\textwidth}

\begin{tikzpicture}

\begin{scope}state without output/.append style={draw=none}%[every node/.style={circle,thick,draw}]
    \node (A) at (3,0) { };
    \node (B) at (0,-3) { };

    \node (C) at (6,-3) { };
   
\end{scope}

\begin{scope} %[>={Stealth[black]},
              every node/.style={fill=white,circle},
              every edge/.style={draw=black,very thick}]
              
    %\path [->] (A) edge [bend left=40] node {} (B);
    %\path [->] (B) edge [bend left=40] node {} (A);
    %\draw [->] (B) to  [out=320,in=40,looseness=9] (B);
    %\path [->] (D) edge node {$3$} (C);
    %\path [->] (A) edge node {$3$} (E);
    %\path [->] (D) edge node {$3$} (E);
    %\path [->] (D) edge node {$3$} (F);
    %\path [->] (C) edge node {$5$} (F);
    %\path [->] (E) edge node {$8$} (F); 
    %\path [->] (B) edge[bend right=60] node {$1$} (E); 

    %\draw[]  (A) node[draw=none][midway,above] {$b_1$} (B);
    \path [->] (A) edge [bend right=10] node [midway,above] {$b_2$} (B);
   
    \path [->] (B) edge [bend right=10] node [midway,above] {$b_2$} (C);
    \path [->] (C) edge [bend right=10] node [midway,right] {$b_2$} (A);
   
    %\draw [->] (A) to  [out=90,in=140,looseness=8] edge node [midway,left] {$c_2$} (A);
    \path [->] (A) edge [bend left=10] node [midway,below] {$b_1$} (B);
   
    \path [->] (B) edge [bend left=10] node [midway,above] {$b_1$} (C);
    \path [->] (C) edge [bend left=10] node [midway,left] {$b_1$} (A);
    
    \filldraw[black] (3,0) circle (2.5pt) node{}[anchor=west];
\filldraw[black] (0,-3) circle (2.5pt) node{}[anchor=west];
\filldraw[black] (6,-3) circle (2.5pt) node{}[anchor=west];
\end{scope}

\end{tikzpicture}
\subcaption*{Connector (2): ker$(\pi_i)=\langle b_1^2b_2,b_1b_2^2\rangle$}
\end{minipage}

\begin{minipage}{.5\textwidth}
\begin{tikzpicture}
\begin{scope}state without output/.append style={draw=none}%[every node/.style={circle,thick,draw}]
    \node (A) at (3,0) { };
    \node (B) at (0,-3) { };

    \node (C) at (6,-3) { };
   
\end{scope}

\begin{scope} %[>={Stealth[black]},
              every node/.style={fill=white,circle},
              every edge/.style={draw=black,very thick}]
              
    %\path [->] (A) edge [bend left=40] node {} (B);
    %\path [->] (B) edge [bend left=40] node {} (A);
    %\draw [->] (B) to  [out=320,in=40,looseness=9] (B);
    %\path [->] (D) edge node {$3$} (C);
    %\path [->] (A) edge node {$3$} (E);
    %\path [->] (D) edge node {$3$} (E);
    %\path [->] (D) edge node {$3$} (F);
    %\path [->] (C) edge node {$5$} (F);
    %\path [->] (E) edge node {$8$} (F); 
    %\path [->] (B) edge[bend right=60] node {$1$} (E); 

    %\draw[]  (A) node[draw=none][midway,above] {$b_1$} (B);
    \path [->] (A) edge [bend right=10] node [midway,above] {$b_2$} (B);
   
    \path [->] (B) edge [bend right=10] node [midway,above] {$b_2$} (C);
    \path [->] (C) edge [bend right=10] node [midway,right] {$b_2$} (A);
   
    %\draw [->] (A) to  [out=90,in=140,looseness=8] edge node [midway,left] {$c_2$} (A);
    \path [<-] (A) edge [bend left=10] node [midway,below] {$b_1$} (B);
   
    \path [<-] (B) edge [bend left=10] node [midway,above] {$b_1$} (C);
    \path [<-] (C) edge [bend left=10] node [midway,left] {$b_1$} (A);
    
\filldraw[black] (3,0) circle (2.5pt) node[anchor=west]{};
\filldraw[black] (0,-3) circle (2.5pt) node[anchor=west]{};
\filldraw[black] (6,-3) circle (2.5pt) node[anchor=west]{};
\end{scope}

\end{tikzpicture}
\subcaption*{Connector (3): ker$(\pi_i)=\langle b_1b_2 \rangle$}
\end{minipage}

\begin{minipage}{.5\textwidth}

\begin{tikzpicture}
\begin{scope}state without output/.append style={draw=none}%[every node/.style={circle,thick,draw}]
    \node (A) at (3,0) { };
    \node (B) at (1,-3) { };

    \node (C) at (5,-3) { };
   
\end{scope}

\begin{scope} %[>={Stealth[black]},
              every node/.style={fill=white,circle},
              every edge/.style={draw=black,very thick}]
              
    %\path [->] (A) edge [bend left=40] node {} (B);
    %\path [->] (B) edge [bend left=40] node {} (A);
    %\draw [->] (B) to  [out=320,in=40,looseness=9] (B);
    %\path [->] (D) edge node {$3$} (C);
    %\path [->] (A) edge node {$3$} (E);
    %\path [->] (D) edge node {$3$} (E);
    %\path [->] (D) edge node {$3$} (F);
    %\path [->] (C) edge node {$5$} (F);
    %\path [->] (E) edge node {$8$} (F); 
    %\path [->] (B) edge[bend right=60] node {$1$} (E); 

    %\draw[]  (A) node[draw=none][midway,above] {$b_1$} (B);
    \path [->] (A) edge node [midway,above] {$b_1$} (B);
   
    \path [->] (B) edge node [midway,above] {$b_1$} (C);
    \path [->] (C) edge node [midway,left] {$b_1$} (A);
   
    %\draw [->] (A) to  [out=90,in=140,looseness=8] edge node [midway,left] {$c_2$} (A);
    \path [->] (A) edge [out=30, in=150, looseness=15, "$b_2$"] (A);
     \path [->] (B) edge [out=150, in=250, looseness=15, "$b_2$"] (B);
      \path [->] (C) edge [out=310, in=40, looseness=15, "$b_2$"] (C);
\filldraw[black] (3,0) circle (2.5pt) node[anchor=west]{};
\filldraw[black] (1,-3) circle (2.5pt) node[anchor=west]{};
\filldraw[black] (5,-3) circle (2.5pt) node[anchor=west]{};
\end{scope}
\end{tikzpicture}
\subcaption*{Connector (4): ker$(\pi_i)=\langle b_2 \rangle$}
\end{minipage}
\caption{All shapes of connectors.}
\label{figcon}
\end{figure}

Let us fragment $\langle b\rangle$. 
 Denote $\xi=2^{\omega}$, $W_n=\mathsf{X}^{\omega}02^n$. Then Supp$(b)=\bigcup\limits_{n=0}^{\infty}W_n$. Let $P_0=\bigcup\limits_{k=0}^{\infty}W_{4k}$, $P_1=\bigcup\limits_{k=0}^{\infty}W_{4k+1}$, $P_2=\bigcup\limits_{k=0}^{\infty}W_{4k+2}$, and $P_3=\bigcup\limits_{k=0}^{\infty}W_{4k+3}$. Then $\mathcal{P}=\{P_0,P_1,P_2,P_3\}$ forms a partition of Supp$(b)$, and each piece is $b$-invariant accumulating on $\xi$. Let $b_1^{(1)}$ act as $b$ on $P_1\cup P_3$, as $b^2$ on $P_2$, and as identity on $P_0$. Let $b_2^{(1)}$ act as $b$ on $P_0\cup P_1\cup P_2$ and as identity on $P_3$. Define $H_1=\langle b_1^{(1)},b_2^{(1)} \rangle$. Consider the subdirect product $H_1\rightarrow \langle b \rangle^4\cong (\mathbb{Z}/3\mathbb{Z})^4$ defined by 
 \begin{align*}
    \begin{cases}
    b_1^{(1)}\mapsto (0,1,2,1),\\
    b_2^{(1)}\mapsto (1,1,1,0).\\
    \end{cases}
\end{align*}
In other words, we have epimorphisms $\pi_i:H_1\rightarrow \mathbb{Z}/3\mathbb{Z}$ defined by $P_i$, for $i=0,1,2,3$. 

By \cite[Proposition 4.5.1]{JC19}, the image of each element of $H_1$ in $(\mathbb{Z}/3\mathbb{Z})^4$ has a coordinate equal to $0$. Define $H_i=\langle b_1^{(i)},b_2^{(i)} \rangle$, for $i=2,3,4$, where $b_1^{(i)}$ and $b_2^{(i)}$ satisfy
\begin{align}\label{gen1}
    \begin{cases}
    b_1^{(1)}=(Id,Id,b_1^{(2)}),\\
    b_1^{(2)}=(a,Id,b_1^{(3)}),\\
    b_1^{(3)}=(a^2,Id,b_1^{(4)}),\\
    b_1^{(4)}=(a,Id,b_1^{(1)}),
    \end{cases}
\end{align}
and 
\begin{align}\label{gen2}
    \begin{cases}
    b_2^{(1)}=(a,Id,b_2^{(2)}),\\
    b_2^{(2)}=(a,Id,b_2^{(3)}),\\
    b_2^{(3)}=(a,Id,b_2^{(4)}),\\
    b_2^{(4)}=(Id,Id,b_2^{(1)}).
    \end{cases}
\end{align}
Then each $H_i$ has a subdirect product defined similarly as that of $H_1$. By the same reasoning, the image of each element of $H_i$ in $(\mathbb{Z}/3\mathbb{Z})^4$ has a coordinate equal to $0$. Each $H_i$ for $i=1,2,3,4$ is isomorphic to $(\mathbb{Z}/3\mathbb{Z})^2$ with invariant domains of supports.   
Let $G_i$ be the groups generated by $\langle a \rangle\cup H_i$ for $i=1,2,3,4$. Let $$\Omega=\{(0123)^{\omega},(1230)^{\omega},(2301)^{\omega},(3012)^{\omega}\}=\{\tau_1,\tau_2,\tau_3,\tau_4\}$$ and $\sigma:\Omega\rightarrow\Omega$ be the shift map. Then the family $F=(\{G_{\tau_i}\}_{\tau_i\in\Omega},\sigma)$, where $G_{\tau_i}=G_i$ for $i=1,2,3,4$, is a self-similar family of groups defined by the wreath recursions (\ref{gen1}) and (\ref{gen2}) along with $a=(012)$, where the wreath recursions are of the form $G_{\tau_i}\rightarrow S_3\wr G_{\sigma(\tau_i)}$. Note that the fragmentations do not determine a stationary group. However, if we start from level $1$ of $\mathsf{B}$ and pass to levels $1+4k$ for $k\in\mathbb{N}$, then we have a self-similar (stationary) group. 

Tile inflations of the unfragmented $G$ are shown in Figure \ref{fgtile}. Let us describe the tiles of $F$. The point $\xi$ is the purely non-Hausdorff singularity. There are two boundary points on each $\mathcal{T}_n$: $\eta_n=02^{n-1}$ and $\xi_n=2^{n}$. The graph of germs of each $G_i$ at $\xi$, denoted $\widetilde{\Gamma}_{\xi,i}$ is obtained by taking $|H_i|$ ($=9$) copies of $\Gamma_{\xi}$ and connecting them at $\xi$ by the Cayley graph of $H_i$. Fix $i_1\in \{1,2,3,4\}$ and consider the action of $G_{i_1}$ on the rooted tree $\mathsf{X}^{*}$. 
%Then the relationship between the tiles $\mathcal{T}_n$ level and the $\mathcal{T}_{k+1}$ is shown in Figure \ref{fgtile}. 
There is a $3$-sheeted covering $\widetilde{\Gamma}_{\xi,i_1}\rightarrow \Lambda_{j}$ mapping the central part of $\widetilde{\Gamma}_{\xi,i_{1}}$ to the central part of $\Lambda_{j}$, where $\Lambda_{j}$ is the limit graph of the orbital graphs of a set of regular points in $P_j\in\mathcal{P}$, for $j=0,1,2,3$, ($\{2^{4k}1v\}_{k\in \mathbb{N}}$, for example) converging to $2^{\omega}$. Hence $\mathcal{T}_{k+1}$ is isomorphic to a central part of $\Lambda_{j}$. 
Let $w\in S^*$ be a word. A rank $k$ traverse of type $h$ of $w$, for $h\in H_{i_1}$ is a subwalk of $w$ on $\mathcal{T}_n$ starting at the boundary point 
%labeled by $12^{\omega}$
, crossing the connectors labeled by $H_{i_1}$ ending in the boundary point labeled by 
%$12^{\omega}$
, and not touching any boundary points in between, and this subwalk is lifted to a walk in $\widetilde{\Gamma}_{\xi,i_{k}}$ starting at the branch corresponding to $Id$ and ending at the branch corresponding to $h$. All shapes of connectors are shown in Figure \ref{figcon}. Note that the higher indices of the generators are omitted since all the groups $G_i$ will have all these $4$ types of connectors. The connectors appear in order (up to cyclic permutations) of (1)(2)(3)(4) in the tile inflations of $G_{i_1}$. Hence the tiles $\mathcal{T}_{k+3}$ will contain all types of rank $k$ bridges, i.e., $\mathcal{T}_{k+3}$ contains all types of connectors that connect isomorphic copies of $\mathcal{T}_n$, for $k\geq 1$. It follows that a rank $k+3$ traverse of type $h$ will contain rank $k$ traverses that pass through all connectors $(1)(2)(3)(4)$.    

\begin{prop} \label{fgreturn} Let $h\neq Id\in H_{i_1}$ and $k\geq 2$. For each $i\in\{1,2,3,4\}$, a rank $k$ traverse of type $h$ passing through the connector $(i)$ is a rank $k$ traverse of type $Id$ passing through $(i')$ for an $i'\neq i$ and $i'\in\{1,2,3,4\}$. 
\end{prop}
\begin{proof} This follows directly from the fragmentation and \cite[Proposition 4.5.1]{JC19}. However, we would like to give an explicit description of such relations. The group $H_{i_1}$ has order $9$ and the set of non-identity elements is $\{b_1,b_2,b_1b_2,b_1^2,b_2^2,b_1^2b_2,b_1b_2^2,b_1^2b_2^2\}$. Suppose the traverse is passing through $(1)$. Then the non-identity types are $$b_2,b_2^2,b_1b_2,b_1b_2^2,b_1^2b_2,b_1^2b_2^2,$$ since $b_1$ is a loop on $(1)$. Then the types $b_2$ and $b_2^2$ are type $Id$ on $(4)$; the types $b_1b_2$ and $b_1^2b_2^2$ are type $Id$ on $(3)$; the types $b_1^2b_2$ and $b_1b_2^2$ are type $Id$ on $(2)$. Similar arguments hold if we start the traverses on bridges with $(2)(3)(4)$. 
\end{proof}

Denote by $\Theta_{k,h}$ the set of rank $k$ traverse of type $h$ and let $\Theta_n=\bigcup\limits_{h\in H_{i_1}}\Theta_{n,h}$. Let $T_{n,h}=\#\Theta_{n,h}$, and $T_n=\sum\limits_{h\in H_{i_1}\backslash\{Id\}}T_{n,h}$. Let $F_h(t)=\sum\limits_{n=0}^{\infty}T_{n,h}t^n$ and $F(t)=\sum\limits_{n=0}^{\infty}T_nt^n$. By the above observations and Proposition \ref{fgreturn}, we have 
$$T_{n+3}\leq\sum\limits_{h\notin\text{Ker}(\pi_i)}T_{n,h},$$ $$\implies t^{-3}(F(t)-(T_0+T_1t+T_2t^2))\leq \sum\limits_{h\notin\text{Ker}(\pi_i)}F_h(t)$$ $$\implies F(t)\leq |w| \dfrac{g(t)}{4t^{-3}-3}$$ where $g(t)>0$ for $t>0$. Hence if $t\in (0,\eta^{-1})$, where $\eta=\sqrt[3]{\frac{3}{4}}$ is the real positive root of the polynomial equation $4x^3-3=0$, the right-hand side of the above inequality is positive. It follows that the growth function for $G_{i_1}$, denoted $\gamma_{i_1}(R)$, is dominated by $\exp(R^{\alpha})$ for every $\alpha>\frac{\log 3}{\log 3-\log\eta}\approx 0.9197$.

%\def\BState{\State\hskip-\ALG@thistlm}
%\makeatother
%%%%%%%%%%%%%%%%%%%%%%%%%%%
%\let\oldbibitem\bibitem
%\renewcommand{\bibitem}{\setlength{\itemsep}{0pt}\oldbibitem}
%%%%%%%%%%%%%%%%%%%%%%%%%%%%%%%%%%%%%%%%%%%%%%%%%%%%%%%%%%%%%%%

%\newpage
%\let\oldbibitem\bibitem
%\renewcommand{\bibitem}{\setlength{\itemsep}{0pt}\oldbibitem}
%\bibliographystyle{ieeetr}

%\addcontentsline{toc}{section}{REFERENCES}

%\renewcommand{\bibname}{{\normalsize\rm REFERENCES}}

%\bibliography{Growth.bib}
%\printbibliography
%\newpage
%\printbibliography[heading=bibintoc]
% title={References}]
%Prints the entire bibliography with the title "Whole bibliography"
%\clearpage
%Filters bibliography
%\printbibliography[heading=subbibintoc,type=article,title={Articles only}]
%\printbibliography[type=book,title={Books only}]
%\printbibliography[keyword={physics},title={Physics-related only}]
%\printbibliography[keyword={latex},title={\LaTeX-related only}]
%\bibliographystyle{alpha}
%\bibliography{Incompress}
\newcommand{\etalchar}[1]{$^{#1}$}

\end{document}